\pdfoutput=1
\documentclass[11pt]{article}
\usepackage[left=1in,right=1in,top=1in,bottom=1in]{geometry}
\usepackage{times}
\usepackage{expl3}
\usepackage{cite}
\usepackage[table]{xcolor}
\usepackage{multirow}
\usepackage{stackengine} 
\usepackage{hhline}
\usepackage{lipsum}
\usepackage{titlesec}
\usepackage{wrapfig}
\usepackage{enumerate}
\usepackage{epsfig}
\usepackage{amsmath}
\usepackage{tabularx}
\usepackage{array}
\usepackage{booktabs}
\usepackage{enumitem}
\usepackage{bbm}
\usepackage{calc}
\usepackage{graphicx}
\usepackage{amsmath}
\usepackage[title]{appendix}
\usepackage{amssymb}
\usepackage{epstopdf}
\usepackage{boldline}
\usepackage{arydshln}
\usepackage{calligra}
\usepackage{bm}
\usepackage{url}
\usepackage{blindtext}
\usepackage{accents}

\newcommand{\define}{\stackrel{\mbox{\tiny def}}{=}}

\newtheorem{theorem}{Theorem}
\newtheorem{proposition}{Proposition}

\newtheorem{lemma}{Lemma}

\newtheorem{example}{Example}

\usepackage{mathtools}
\usepackage{epstopdf}
\usepackage{balance}
\usepackage{thmtools}
\usepackage{thm-restate}
\usepackage{hyperref}
\usepackage{cleveref}
\usepackage[mathscr]{euscript}

\usepackage[ruled,vlined]{algorithm2e}
\include{pythonlisting}

\newcommand{\ostar}{\mathbin{\mathpalette\make@circled\star}}

\makeatletter
\newcommand{\removelatexerror}{\let\@latex@error\@gobble}
\makeatother
\makeatletter
\newcommand*{\rom}[1]{\expandafter\@slowromancap\romannumeral #1@}
\makeatother

\ExplSyntaxOn
\newcommand\latinabbrev[1]{
  \peek_meaning:NTF. {
    #1\@}%
  { \peek_catcode:NTF a {
      #1.\@ }%
    {#1.\@}}}
\ExplSyntaxOff

\titleclass{\subsubsubsection}{straight}[\subsubsection]

\begin{document}
\vspace{1cm}
\title{Generalized Multiple Operator Integrals and Perturbation Theory for Operators with Continuous Spectra}
\vspace{1.8cm}
\author{Shih-Yu~Chang
\thanks{Shih-Yu Chang is with the Department of Applied Data Science,
San Jose State University, San Jose, CA, U. S. A. (e-mail: {\tt
shihyu.chang@sjsu.edu}). 
           }}

\maketitle

\begin{abstract}
Operators with continuous spectra naturally arise in spectral theory, quantum mechanics, automorphic forms, and noncommutative geometry. However, analyzing such operators—particularly in the non-self-adjoint setting—remains challenging due to spectral instability and the lack of an orthonormal basis. This work advances the theory of Multiple Operator Integrals (MOIs) by developing a unified framework for generalized MOIs (GMOIs) associated with general (non-normal, non-self-adjoint) operators possessing continuous spectra. Building on prior work in Generalized Double Operator Integrals (GDOIs) and finite-dimensional GMOIs, we extend the theory to include: the formulation of GMOIs in the continuous spectrum setting, their algebraic structure, continuity properties, norm and Lipschitz estimates, and a perturbation formula that generalizes classical results. As a key application, we derive a Krein-type spectral shift formula for GDOIs in the continuous spectrum setting and further extend it to arbitrary-order approximations. These contributions provide a foundation for broader developments in spectral theory, operator algebras, noncommutative geometry, and noncommutative analysis.
\end{abstract}


\section{Introduction}\label{sec: Introduction}

Operators with continuous spectrum play a crucial role in various fields, e.g., spectral theory, automorphic forms, noncommutative geometry and quantum physics. Different from operators with discrete spectra, which are often associated with classical linear algebra and finite-dimensional eigenspaces, the continuous spectrum emerges naturally in the analysis of differential operators on noncompact domains, such as the Laplace--Beltrami operator on modular surfaces \cite{iwaniec1997topics, hejhal1976selberg, sarnak1990some}. These spectra encode rich arithmetic structures, including connections to prime number distribution, modular forms, and zeta functions via the Selberg trace formula and Eisenstein series \cite{connes1999trace, connes2017riemann}. Moreover, continuous spectral components are essential in the formulation of scattering theory and quantum mechanics, where self-adjoint operators like the Schrödinger operator exhibit purely continuous spectra for free particles. In this context, spectral resolution requires integration over a continuum of eigenstates, often realized through spectral measures and direct integral decompositions. This feature significantly complicates the definition of operator functions and functional calculus compared to the discrete case.

In noncommutative geometry, as formulated by Connes, continuous spectra arise in spectral triples, where the Dirac operator's spectrum captures the geometric and topological content of noncommutative spaces \cite{connes1994noncommutative}. Similarly, in representation theory, especially of noncompact Lie groups, continuous spectra appear in the decomposition of unitary representations, often necessitating the use of Plancherel-type theorems and intertwining operators. These challenges have motivated the development of generalized functional calculus and operator integration methods beyond the classical spectral theorem over discrete spectrum.

Consequently, there is a growing need for systematic approaches to handle multi-operator integration and perturbation theory in settings involving continuous spectra. Such advancements enable the extension of tools like the spectral shift function, double and multiple operator integrals, and resolvent expansions to broader analytic and algebraic frameworks. From mathematics development perspective, such ongoing development not only deepens our understanding of operator algebras but also fosters applications in mathematical physics, arithmetic geometry, and noncommutative analysis. In this broader context, the development of multi-operator integration frameworks involving continuous-spectrum operators has become an important analytical and algebraic tool. For example, scattering theory, spectral representations of automorphic forms, and noncommutative geometric models frequently invoke simultaneous or interacting operator actions, where each operator contributes a potentially continuous part of the spectrum \cite{langlands1976functional, arthur2005introtrace}.

The investigation of non-self-adjoint operators with continuous spectra has attracted attention in both mathematical science and theoretical physics. Unlike self-adjoint operators, which have a rich spectral theory grounded within the spectral theorem and useful calculus, non-self-adjoint operators present signficant technical difficulties due to the lack of an orthonormal eigenbasis and the potential instability in their spectrum under perturbations~\cite{trefethen2005spectra, davies2007linear}. In quantum mechanics and wave propagation issues, such operators are induced naturally while modeling open structures, dissipative phenomena, or complex potentials—especially in $\mathcal{PT}$-symmetric quantum principle~\cite{bender1998real, bender2007making}.  From a mathematical perspective, non-self-adjoint operators also appear in the study of resonance expansions, non-normal evolution semigroups, and control theory~\cite{engel2000one, zworski2017mathematical}. The presence of continuous spectrum also beings complicaiton to the development of spectral decompositions, requiring more advanced tools such as rigged Hilbert areas, dilation principle, or generalized resolvent strategies. Consequently, the development of operator integration frameworks and perturbation idea for non-self-adjoint continuous-spectrum operators is essential for expanding the limits of functional analysis, spectral theory, and mathematical physics.

However, most comventional  tasks mainly focused on special classes of operators such as self-adjoint, unitary, or normal operators, where the spectral theorem provides a well-established theoretical foundation. In contrast, general (non-normal, non-self-adjoint) operators with continuous spectrum remain relatively underexplored. To address this gap, we have initiated a systematic investigation into multiple-operator integrals (MOIs) and double-operator integrals (DOIs) beyond the classical framework in a series of prior studies. These include the development of DOIs for finite-dimensional general operators~\cite{chang2025GDOIMatrix}, the extension of DOIs to continuous-spectrum general operators~\cite{chang2025GDOICont}, and the formulation of MOIs for finite-dimensional general operators~\cite{chang2025GMOIFinite}. These foundational efforts lay the groundwork for the present work, which further advances MOI theory in the context of continuous spectra.

This work contributes to the ongoing development of operator integration theory by extending the framework of generalized multiple operator integrals (GMOIs) to accommodate general (non-self-adjoint, non-normal) operators with continuous spectrum. Building upon the finite-dimensional theory introduced in~\cite{chang2025GMOIFinite}, we develop a new approach that unifies spectral mapping principles with multi-operator structures arising in noncommutative geometry, singular integral analysis, and representation theory. The key contributions of this work include the formulation of a generalized multi-operator integral (GMOI) framework for operators with continuous spectra, the establishment of its fundamental algebraic properties, a rigorous analysis of its continuity under suitable conditions, GMOI norm and Lipschitz estimates and the derivation of a perturbation formula that extends classical results to this broader setting.  As an application, we first establish a Krein-type spectral shift formula for generalized double operator integrals (GDOIs) associated with continuous spectra, and extending this to the $(n-1)$-th order approximation of the spectral shift formula.

The remainder of this paper is organized as follows. In Section~\ref{sec:MOIs as Special Cases of the Spectral Mapping Theorem for Continuous Spectrum Operators}, we present a perspective that interprets classical multiple operator integrals (MOIs) as special cases of the spectral mapping theorem for continuous spectrum operators. Section~\ref{sec:Generalized Multiple Operator Integrals for Continuous Spectrum Operators and Their Algebraic Properties} develops the notion of generalized multiple operator integrals (GMOIs) for continuous spectrum operators and explores their fundamental algebraic properties. In Section~\ref{sec:Perturbation Formula for GMOIs with Continuous Spectrum Operators}, we derive a perturbation formula involving GMOIs, extending classical results to the continuous spectrum case. Section~\ref{sec:Inequalities Related to GMOI for Continuous Spectrum Operators} establishes key inequalities related to GMOIs, providing analytical tools for further study. The continuity properties of GMOIs are then addressed in Section~\ref{sec:Continuity of GMOI for Continuous Spectrum Operators}, where we prove several convergence results. Finally, in Section~\ref{sec:Applications:Spectral Shift}, we discuss applications of the developed theory to spectral shift functions and trace formulas.

\section{A Spectral Mapping Perspective on MOIs for Operators with Continuous Spectrum}\label{sec:MOIs as Special Cases of the Spectral Mapping Theorem for Continuous
Spectrum Operators}

Let us consider conventional MOI definitions. Given a function $\beta: \mathbb{R}^{\zeta+1} \rightarrow \mathbb{C}$, $\zeta+1$ self-adjoint operators (parameter operators) $\bm{X}_1,\bm{X}_2,\ldots, \bm{X}_{\zeta+1}$, and any $\zeta$ operators (argument operators) $\bm{Y}_1,\bm{Y}_2,\ldots, \bm{Y}_{\zeta}$. From spectral mapping theorem, we have parameter operators
\begin{eqnarray}\label{eq0-1:  conv DOI def}
\bm{X}_1&=&\int\limits_{\lambda_1\in \sigma(\bm{X}_1)}\lambda_1 d\bm{E}_{\bm{X}_1}(\lambda_1);\nonumber \\
\bm{X}_2&=&\int\limits_{\lambda_2\in \sigma(\bm{X}_2)}\lambda_2 d\bm{E}_{\bm{X}_2}(\lambda_2);\nonumber \\
\vdots&&\nonumber \\
\bm{X}_{\zeta+1}&=&\int\limits_{\lambda_{\zeta+1}\in \sigma(\bm{X}_{\zeta+1})}\lambda_{\zeta+1} d\bm{E}_{\bm{X}_{\zeta+1}}(\lambda_{\zeta+1}),
\end{eqnarray}
where $\lambda_1$, $\lambda_2$, $\cdots$, and $\lambda_{\zeta+1}$ are eigenvalues of the operators $\bm{X}_1$, $\bm{X}_2$, $\cdots$, and $\bm{X}_{\zeta+1}$, respectively; and  $\bm{E}_{\bm{X}_1}(\lambda_1), \bm{E}_{\bm{X}_2}(\lambda_2)$, $\cdots$, and $\bm{E}_{\bm{X}_{\zeta+1}}(\lambda_{\zeta+1})$ are spectrum measures with respect to eigenvalues $\lambda_1, \lambda_2$, $\cdots$ and $\lambda_{\zeta+1}$, respectively. 

Similarly, from spectral mapping theorem, we have argument operators
\begin{eqnarray}\label{eq0-2:  conv DOI def}
\bm{Y}_1&=&\int\limits_{\lambda'_{1} \in \sigma(\bm{Y}_1)}\lambda'_{1} d\bm{E}_{\bm{Y}_1}(\lambda'_{1})+
\int\limits_{\lambda'_{1} \in \sigma(\bm{Y}_1)}\left(\bm{Y}_1-\lambda'_{1}\bm{I}\right)d\bm{E}_{\bm{Y}_1}(\lambda'_{1});\nonumber \\
\bm{Y}_2&=&\int\limits_{\lambda'_{2} \in \sigma(\bm{Y}_2)}\lambda'_{2} d\bm{E}_{\bm{Y}_2}(\lambda'_{2})+
\int\limits_{\lambda'_{2} \in \sigma(\bm{Y}_2)}\left(\bm{Y}_2-\lambda'_{2}\bm{I}\right)d\bm{E}_{\bm{Y}_2}(\lambda'_{2});\nonumber \\
\vdots&&\nonumber \\
\bm{Y}_\zeta&=&\int\limits_{\lambda'_{\zeta} \in \sigma(\bm{Y}_\zeta)}\lambda'_{\zeta} d\bm{E}_{\bm{Y}_\zeta}(\lambda'_{\zeta})+
\int\limits_{\lambda'_{\zeta} \in \sigma(\bm{Y}_\zeta)}\left(\bm{Y}_\zeta-\lambda'_{\zeta}\bm{I}\right)d\bm{E}_{\bm{Y}_\zeta}(\lambda'_{\zeta}),
\end{eqnarray}
where $\lambda'_1$, $\lambda'_2$, $\cdots$, and $\lambda'_{\zeta}$ are eigenvalues of the operators $\bm{Y}_1$, $\bm{Y}_2$, $\cdots$, and $\bm{Y}_{\zeta}$, respectively; and  $\bm{E}_{\bm{Y}_1}(\lambda'_1)$, $\bm{E}_{\bm{Y}_2}(\lambda'_2)$, $\cdots$, and $\bm{E}_{\bm{Y}_{\zeta}}(\lambda'_{\zeta})$ are spectrum measures with respect to eigenvalues $\lambda'_1, \lambda'_2$, $\cdots$ and $\lambda'_{\zeta}$, respectively. 

We review Theorem 11 in~\cite{chang2024operatorChar} first. Before presenting this theorem, we review several special ntations related to this Theorem 11 in~\cite{chang2024operatorChar}.

Given $r$ positive integers $q_1, q_2, \ldots, q_r$, we define $\alpha_{\kappa}(q_1, \ldots, q_r)$ to be the selection of these $r$ arguments $q_1, \ldots, q_r$ into $\kappa$ arguments, i.e., we have
\begin{equation}
\alpha_{\kappa}(q_1, \ldots, q_r) = \{ q_{\iota_1}, q_{\iota_2}, \ldots, q_{\iota_\kappa} \}.
\end{equation}

We use $\mathrm{Ind}(\alpha_{\kappa}(q_1, \ldots, q_r))$ to obtain the indices of those $\kappa$ positive integers $\{ q_{\iota_1}, q_{\iota_2}, \ldots, q_{\iota_\kappa} \}$, i.e., we have
\begin{equation}
\mathrm{Ind}(\alpha_{\kappa}(q_1, \ldots, q_r)) = \{ \iota_1, \iota_2, \ldots, \iota_\kappa \}.
\end{equation}

We use $\alpha_{\kappa}(q_1, \ldots, q_r) = 1$ to represent $q_{\iota_1} = 1, q_{\iota_2} = 1, \ldots, q_{\iota_\kappa} = 1$. 

We also use
\[
m_{\lambda_{\mbox{Ind}(\alpha_{\kappa}(q_1,\ldots,q_r))}}-1
\]
to represent
\[
m_{\lambda_{\iota_1}}-1,\ldots,m_{\lambda_{\iota_\kappa}}-1,
\]
where $m_{\lambda_{\iota_j}}$ is the order of the nilpotent  $\left(\bm{X}_{\iota_j}-\lambda_{\iota_j}\bm{I}\right)d\bm{E}_{\bm{X}_{\iota_j}}(\lambda_{\iota_j})$, i.e., 
\[
\left(\bm{X}_{\iota_j}-\lambda_{\iota_j}\bm{I}\right)^{\ell}d\bm{E}_{\bm{X}_{\iota_j}}(\lambda_{\iota_j})= \bm{0}, \quad \text{for } \ell \geq m_{\lambda_{\iota_j}}  \text{ and } j = 1, 2, \ldots, \kappa.
\]

Theorem 11 from~\cite{chang2024operatorChar} is given below. 
\begin{theorem}\label{thm: Spectral Mapping Theorem for r Variables inf}
Given an analytic function $f(z_1,z_2,\ldots,z_r)$ within the domain for $|z_l| < R_l$, and the operator $\bm{X}_l$ decomposed by:
\begin{eqnarray}\label{eq1-1: thm: Spectral Mapping Theorem for r Variables inf}
\bm{X}_l&=&\int\limits_{\lambda_l \in \sigma(\bm{X}_l)}\lambda_l d\bm{E}_{\bm{X}_l}(\lambda_l)+
\int\limits_{\lambda_l \in \sigma(\bm{X}_l)}\left(\bm{X}_l-\lambda_l\bm{I}\right)d\bm{E}_{\bm{X}_l}(\lambda_l),
\end{eqnarray}
where $\left\vert\lambda_{l}\right\vert<R_l$ for $l=1,2,\ldots,r$.

Then, we have
\begin{eqnarray}\label{eq2: thm: Spectral Mapping Theorem for kappa Variables inf}
\lefteqn{f(\bm{X}_1,\ldots,\bm{X}_r)=}\nonumber \\
&&\int\limits_{\lambda_1 \in \sigma(\bm{X}_1)}\cdots\int\limits_{\lambda_r \in \sigma(\bm{X}_r)}
f(\lambda_1,\ldots,\lambda_r)d\bm{E}_{\bm{X}_1}(\lambda_1)\cdots d\bm{E}_{\bm{X}_r}(\lambda_r) \nonumber \\
&&+\int\limits_{\lambda_1 \in \sigma(\bm{X}_1)}\cdots\int\limits_{\lambda_r \in \sigma(\bm{X}_r)}\sum\limits_{\kappa=1}^{r-1}\sum\limits_{\alpha_\kappa(q_1,\ldots,q_r)}\Bigg(\sum\limits_{\alpha_{\kappa}(q_1,\ldots,q_r)=1}^{m_{\lambda_{\mbox{Ind}(\alpha_{\kappa}(q_1,\ldots,q_r))}}-1}\nonumber \\
&&~~~~~ \frac{f^{\alpha_{\kappa}(q_1,\ldots,q_r)}(\lambda_1,\ldots,\lambda_r)}{q_{\iota_1}!q_{\iota_2}!\ldots q_{\iota_\kappa}!}\times \prod\limits_{\substack{\beta' =\mbox{Ind}(\alpha_{\kappa}(q_1,\ldots,q_r)), \bm{Y}=\left(\bm{X}_{\beta'} - \lambda_{\beta'}\bm{I}\right)^{q_{\beta'}}d\bm{E}_{\bm{X}_{\beta'}}(\lambda_{\beta'}) \\ \beta' \neq \mbox{Ind}(\alpha_{\kappa}(q_1,\ldots,q_r)), \bm{Y}=d\bm{E}_{\bm{X}_{\beta'}}(\lambda_{\beta'})}
}^{r} \bm{Y}\Bigg) 
\nonumber \\
&&+\int\limits_{\lambda_1 \in \sigma(\bm{X}_1)}\cdots\int\limits_{\lambda_r \in \sigma(\bm{X}_r)}\sum\limits_{q_1=\ldots=q_r=1}^{m_{\lambda_1}-1,\ldots,m_{\lambda_r}-1}
\frac{f^{(q_1,\ldots,q_r)}(\lambda_1,\ldots,\lambda_r)}{q_1!\cdots q_r!}\nonumber \\
&&\times \left(\bm{X}_1 - \lambda_1\bm{I}\right)^{q_1}d\bm{E}_{\bm{X}_1}(\lambda_1) \left(\bm{X}_2 - \lambda_2\bm{I}\right)^{q_2}d\bm{E}_{\bm{X}_2}(\lambda_2)\cdots \left(\bm{X}_r - \lambda_r\bm{I}\right)^{q_r}d\bm{E}_{\bm{X}_r}(\lambda_r),
\end{eqnarray}
where we have
\begin{itemize}
\item $\sum\limits_{\alpha_\kappa(q_1,\ldots,q_r)}$ runs over all $\kappa$ selections of $q_1,\ldots,q_r$ by $\alpha_\kappa(q_1,\ldots,q_r)$;
\item $m_{\lambda_{\mbox{Ind}(\alpha_{\kappa}(q_1,\ldots,q_r))}}-1 =$ $m_{\lambda_{\iota_1}}-1$,$\ldots,m_{\lambda_{\iota_\kappa}}-1$;
\item $f^{\alpha_{\kappa}(q_1,\ldots,q_r)}(\lambda_1,\ldots,\lambda_r)$ represents the partial derivatives with respect to variables with indices $\iota_1,\iota_2,\ldots,\iota_\kappa$ and the orders of derivatives given by $q_{\iota_1},q_{\iota_2},\ldots,q_{\iota_\kappa}$.
\end{itemize}
\end{theorem}

The conventional MOI is an operator, denoted by $T_{\beta}^{\bm{X}_1,\ldots,\bm{X}_{\zeta+1}}(\bm{Y}_1,\ldots,\bm{Y}_{\zeta})$, which can be expressed as~\cite{skripka2019multilinear}:
\begin{eqnarray}\label{eq1:  conv MOI def}
\lefteqn{T_{\beta}^{\bm{X}_1,\ldots,\bm{X}_{\zeta+1}}(\bm{Y}_1,\ldots,\bm{Y}_{\zeta})\define}\nonumber \\&=&\int\limits_{\lambda_1 \in \sigma(\bm{X}_1)}\cdots\int\limits_{\lambda_{\zeta+1} \in \sigma(\bm{X}_{\zeta+1})} d\bm{E}_{\bm{X}_1}(\lambda_1) \bm{Y}_1 d\bm{E}_{\bm{X}_2}(\lambda_2) \ldots d\bm{E}_{\bm{X}_\zeta}(\lambda_\zeta) \bm{Y}_\zeta d\bm{E}_{\bm{X}_{\zeta+1}}(\lambda_{\zeta+1})
\end{eqnarray}
From the decomposition of the operator $\bm{Y}_j$ given by Eq.~\eqref{eq0-2:  conv DOI def}, Eq.~\eqref{eq1:  conv MOI def} can further be expressed as
\begin{eqnarray}\label{eq2:  conv MOI def}
\lefteqn{T_{\beta}^{\bm{X}_1,\ldots,\bm{X}_{\zeta+1}}(\bm{Y}_1,\ldots,\bm{Y}_{\zeta})=\int\limits_{\lambda_1 \in \sigma(\bm{X}_1)}\cdots\int\limits_{\lambda_{\zeta+1}\in \sigma(\bm{X}_{\zeta+1})}\beta(\lambda_{1},\ldots, \lambda_{\zeta+1})}\nonumber \\
&& \times d\bm{E}_{\bm{X}_1}(\lambda_1) \left(\int\limits_{\lambda'_{1} \in \sigma(\bm{Y}_1)}\lambda'_{1} d\bm{E}_{\bm{Y}_1}(\lambda'_{1})+
\int\limits_{\lambda'_{1} \in \sigma(\bm{Y}_1)}\left(\bm{Y}_1-\lambda'_{1}\bm{I}\right)d\bm{E}_{\bm{Y}_1}(\lambda'_{1})\right)d\bm{E}_{\bm{X}_2}(\lambda_2) \ldots \nonumber \\
&& \times d\bm{E}_{\bm{X}_\zeta}(\lambda_\zeta) \left(\int\limits_{\lambda'_{\zeta} \in \sigma(\bm{Y}_\zeta)}\lambda'_{\zeta} d\bm{E}_{\bm{Y}_\zeta}(\lambda'_{\zeta})+
\int\limits_{\lambda'_{\zeta} \in \sigma(\bm{Y}_\zeta)}\left(\bm{Y}_\zeta-\lambda'_{\zeta}\bm{I}\right)d\bm{E}_{\bm{Y}_\zeta}(\lambda'_{\zeta}) \right)d\bm{E}_{\bm{X}_{\zeta+1}}(\lambda_{\zeta+1}) \nonumber \\
&=& \int\limits_{\lambda_1 \in \sigma(\bm{X}_1)}\cdots\int\limits_{\lambda_{\zeta+1}\in \sigma(\bm{X}_{\zeta+1})}  \int\limits_{\lambda'_1 \in \sigma(\bm{Y}_1)}\cdots\int\limits_{\lambda'_{\zeta}\in \sigma(\bm{Y}_{\zeta})}
\beta(\lambda_{1},\ldots,\lambda_{\zeta+1})\lambda'_{1}\cdots\lambda'_{\zeta}\nonumber \\
&& \times d\bm{E}_{\bm{X}_1}(\lambda_1)d\bm{E}_{\bm{Y}_1}(\lambda'_1)\cdots d\bm{E}_{\bm{Y}_\zeta}(\lambda'_\zeta)d\bm{E}_{\bm{X}_{\zeta+1}}(\lambda_{\zeta+1})\nonumber \\
&&+\int\limits_{\lambda_1 \in \sigma(\bm{X}_1)}\cdots\int\limits_{\lambda_{\zeta+1}\in \sigma(\bm{X}_{\zeta+1})}  \int\limits_{\lambda'_1 \in \sigma(\bm{Y}_1)}\cdots\int\limits_{\lambda'_{\zeta}\in \sigma(\bm{Y}_{\zeta})} \beta(\lambda_1,\ldots,\lambda_{k_{\zeta+1}})\sum\limits_{\kappa=1}^{\zeta-1}\sum\limits_{\varrho_\kappa(1,\ldots,\zeta)}\nonumber \\
&& \Bigg(\prod\limits_{\substack{\varsigma =\mbox{Ind}(\varrho_{\kappa}(1,\ldots,\zeta)), \bm{Z}=\left(\bm{Y}_{\varsigma '} - \lambda'_{\varsigma }\bm{I}\right)d\bm{E}_{\bm{Y}_{\varsigma }}(\lambda'_{\varsigma})  \\ \varsigma \neq \mbox{Ind}(\varrho_{\kappa}(1,\ldots,\zeta)), \bm{Z}=\lambda'_\varsigma d\bm{E}_{\bm{Y}_{\varsigma }}(\lambda'_{\varsigma})}}^{\zeta} \bm{Z}\Bigg) 
\nonumber \\
&&+\int\limits_{\lambda_1 \in \sigma(\bm{X}_1)}\cdots\int\limits_{\lambda_{\zeta+1}\in \sigma(\bm{X}_{\zeta+1})}  \int\limits_{\lambda'_1 \in \sigma(\bm{Y}_1)}\cdots\int\limits_{\lambda'_{\zeta}\in \sigma(\bm{Y}_{\zeta})} 
\beta(\lambda_{1},\ldots,\lambda_{\zeta+1})\nonumber \\
&& \times d\bm{E}_{\bm{X}_1}(\lambda_1)\left(\bm{Y}_{1} - \lambda'_{1}\bm{I}\right)d\bm{E}_{\bm{Y}_{1}}(\lambda'_{1})d\bm{E}_{\bm{X}_2}(\lambda_2) \cdots \nonumber \\
&& \times d\bm{E}_{\bm{X}_\zeta}(\lambda_\zeta) \left(\bm{Y}_{\zeta} - \lambda'_{\zeta}\bm{I}\right)d\bm{E}_{\bm{Y}_{\zeta}}(\lambda'_{\zeta})d\bm{E}_{\bm{X}_{\zeta+1}}(\lambda_{\zeta+1}),
\end{eqnarray}
where $\varrho_{\kappa}(1,\ldots,\zeta)$ is the selection of $\kappa$ indices from the indices set $\{1,2,\ldots,\zeta\}$, and $\mbox{Ind}(\varrho_{\kappa}(1,\ldots,\zeta))$ represents those $\kappa$ indices being selected. 

If we set the function $f$ in Theorem~\ref{thm: Spectral Mapping Theorem for r Variables inf} by
\begin{eqnarray}
f(z_1,z_2, \ldots, z_{2\zeta+1}) &=& \beta(z_1, z_3, \ldots,z_{2\zeta+1})z_2 z_4 \ldots z_{2\zeta}, 
\end{eqnarray}
and set 
\begin{eqnarray}
z_1&=&\lambda_1,~~z_3~=~\lambda_2,~~z_5~=~\lambda_3, ~~\cdots,~~z_{2\zeta+1}~=~\lambda_{k_{\zeta+1}}, \nonumber \\  
z_2&=&\lambda_{k'_1},~~z_4~=~\lambda_{k'_2},~~z_6~=~\lambda_{k'_3}, ~~\cdots,~~z_{2\zeta}~=~\lambda_{k'_\zeta},
\end{eqnarray}
then, by applying Theorem~\ref{thm: Spectral Mapping Theorem for r Variables inf}, we will obtain Eq.~\eqref{eq2:  conv MOI def}.

Following Example~\ref{exp:TOI is a special case} shows that the conventional MOI is a special case of spectrum mapping theorem when $\zeta=2$.
\begin{example}\label{exp:TOI is a special case}
Let us consider a triple operator integral $T_\beta^{\bm{X}_1,\bm{X}_2,\bm{X}_3}(\bm{Y}_1, \bm{Y}_2)$ with self-adjoint operators $\bm{X}_1,  \bm{X}_2$ and $\bm{X}_3$. Then, we can express $T_\beta^{\bm{X}_1,\bm{X}_2,\bm{X}_3}(\bm{Y}_1, \bm{Y}_2)$ as 
\begin{eqnarray}\label{eq1:exp:TOI is a special case}
\lefteqn{T_{\beta}^{\bm{X}_1,\bm{X}_2,\bm{X}_3}(\bm{Y}_1, \bm{Y}_2)}\nonumber\\
&\define&\int\limits_{\lambda_1 \in \sigma(\bm{X}_1)}\int\limits_{\lambda_2 \in \sigma(\bm{X}_2)}\int\limits_{\lambda_3 \in \sigma(\bm{X}_3)}  \beta(\lambda_{1}, \lambda_{2}, \lambda_{3}) d\bm{E}_{\bm{X}_1}(\lambda_1) \bm{Y}_1 d\bm{E}_{\bm{X}_2}(\lambda_2)\bm{Y}_2 d\bm{E}_{\bm{X}_3}(\lambda_3)\nonumber \\
&=&\int\limits_{\lambda_1 \in \sigma(\bm{X}_1)}\int\limits_{\lambda_2 \in \sigma(\bm{X}_2)}\int\limits_{\lambda_3 \in \sigma(\bm{X}_3)} \beta(\lambda_{1}, \lambda_{2}, \lambda_{3})d\bm{E}_{\bm{X}_1}(\lambda_1) \nonumber \\
&& \times \left( \int\limits_{\lambda'_{1} \in \sigma(\bm{Y}_1)}\lambda'_{1} d\bm{E}_{\bm{Y}_1}(\lambda'_{1})+
\int\limits_{\lambda'_{1} \in \sigma(\bm{Y}_1)}\left(\bm{Y}_1-\lambda'_{1}\bm{I}\right)d\bm{E}_{\bm{Y}_1}(\lambda'_{1}) \right) d\bm{E}_{\bm{X}_2}(\lambda_2) \nonumber \\
&& \times \left( \int\limits_{\lambda'_{2} \in \sigma(\bm{Y}_2)}\lambda'_{2} d\bm{E}_{\bm{Y}_2}(\lambda'_{2})+
\int\limits_{\lambda'_{2} \in \sigma(\bm{Y}_2)}\left(\bm{Y}_2-\lambda'_{2}\bm{I}\right)d\bm{E}_{\bm{Y}_2}(\lambda'_{2}) \right) d\bm{E}_{\bm{X}_3}(\lambda_3) \nonumber \\
&=&\int\limits_{\lambda_1 \in \sigma(\bm{X}_1)}\int\limits_{\lambda_2 \in \sigma(\bm{X}_2)}\int\limits_{\lambda_3 \in \sigma(\bm{X}_3)}  \beta(\lambda_{1}, \lambda_{2}, \lambda_{3}) d\bm{E}_{\bm{X}_1}(\lambda_1)\nonumber \\
&& \times \left( \int\limits_{\lambda'_1 \in \sigma(\bm{Y}_1)} \lambda'_1 d\bm{E}_{\bm{Y}_1}(\lambda'_1) \right)d\bm{E}_{\bm{X}_2}(\lambda_2) \left( \int\limits_{\lambda'_2 \in \sigma(\bm{Y}_2)} \lambda'_2 d\bm{E}_{\bm{Y}_2}(\lambda'_2) \right) d\bm{E}_{\bm{X}_3}(\lambda_3) \nonumber \\
&&+\int\limits_{\lambda_1 \in \sigma(\bm{X}_1)}\int\limits_{\lambda_2 \in \sigma(\bm{X}_2)}\int\limits_{\lambda_3 \in \sigma(\bm{X}_3)}\beta(\lambda_1, \lambda_2, \lambda_3)d\bm{E}_{\bm{X}_1}(\lambda_1)\nonumber \\
&& \times \left( \int\limits_{\lambda'_1 \in \sigma(\bm{Y}_1)} \lambda'_1 d\bm{E}_{\bm{Y}_1}(\lambda'_1)  \right)d\bm{E}_{\bm{X}_2}(\lambda_2)\left( \int\limits_{\lambda'_{2} \in \sigma(\bm{Y}_2)}\left(\bm{Y}_2-\lambda'_{2}\bm{I}\right)d\bm{E}_{\bm{Y}_2}(\lambda'_{2})\right)d\bm{E}_{\bm{X}_3}(\lambda_3)  \nonumber \\
&& +\int\limits_{\lambda_1 \in \sigma(\bm{X}_1)}\int\limits_{\lambda_2 \in \sigma(\bm{X}_2)}\int\limits_{\lambda_3 \in \sigma(\bm{X}_3)}\beta(\lambda_1, \lambda_2, \lambda_3)d\bm{E}_{\bm{X}_1}(\lambda_1)\nonumber \\
&& \times \left( \int\limits_{\lambda'_{1} \in \sigma(\bm{Y}_1)}\left(\bm{Y}_1-\lambda'_{1}\bm{I}\right)d\bm{E}_{\bm{Y}_1}(\lambda'_{1}) \right)d\bm{E}_{\bm{X}_2}(\lambda_2)\left( \int\limits_{\lambda'_2 \in \sigma(\bm{Y}_2)} \lambda'_2 d\bm{E}_{\bm{Y}_2}(\lambda'_2) \right)d\bm{E}_{\bm{X}_3}(\lambda_3) \nonumber \\
&&+\int\limits_{\lambda_1 \in \sigma(\bm{X}_1)}\int\limits_{\lambda_2 \in \sigma(\bm{X}_2)}\int\limits_{\lambda_3 \in \sigma(\bm{X}_3)}\beta(\lambda_1, \lambda_2, \lambda_3)d\bm{E}_{\bm{X}_1}(\lambda_1)\nonumber \\
&& \times \left( \int\limits_{\lambda'_{1} \in \sigma(\bm{Y}_1)}\left(\bm{Y}_1-\lambda'_{1}\bm{I}\right)d\bm{E}_{\bm{Y}_1}(\lambda'_{1}) \right)d\bm{E}_{\bm{X}_2}(\lambda_2) \nonumber \\
&&\times \left( \int\limits_{\lambda'_{2} \in \sigma(\bm{Y}_2)}\left(\bm{Y}_2-\lambda'_{2}\bm{I}\right)d\bm{E}_{\bm{Y}_2}(\lambda'_{2}) \right)d\bm{E}_{\bm{X}_3}(\lambda_3)
\end{eqnarray}

If we set the function $f$ in Theorem~\ref{thm: Spectral Mapping Theorem for r Variables inf} by
\begin{eqnarray}
f(z_1,z_2,z_3,z_4, z_5) &=& \beta(z_1, z_3, \,z_5)z_2 z_4, 
\end{eqnarray}
and set 
\begin{eqnarray}
z_1&=&\lambda_1,~~z_3~=~\lambda_2,~~z_5~=~\lambda_3, \nonumber \\  
z_2&=&\lambda_{k'_1},~~z_4~=~\lambda_{k'_2},
\end{eqnarray}
then, by applying Theorem~\ref{thm: Spectral Mapping Theorem for r Variables inf}, we have
\begin{eqnarray}
\frac{\partial f(z_1,z_2,z_3,z_4, z_5)}{\partial z_2} &=&  \beta(z_1, z_3, \,z_5)z_4;\nonumber \\
\frac{\partial f(z_1,z_2,z_3,z_4, z_5)}{\partial z_4} &=&  \beta(z_1, z_3, \,z_5)z_2;\nonumber \\
\frac{\partial f(z_1,z_2,z_3,z_4, z_5)}{\partial z_2 \partial z_4} &=&  \beta(z_1, z_3, \,z_5). 
\end{eqnarray}
Therefore, we obtain:
\begin{eqnarray}\label{eq2:exp:TOI is a special case}
\lefteqn{f(\bm{X}_1,\bm{Y}_1,\bm{X}_2,\bm{Y}_2,\bm{X}_3)}\nonumber \\
&=&\int\limits_{\lambda_1 \in \sigma(\bm{X}_1)}\int\limits_{\lambda_2 \in \sigma(\bm{X}_2)}\int\limits_{\lambda_3 \in \sigma(\bm{X}_3)}\beta(\lambda_{1}, \lambda_{2}, \lambda_{3})\lambda'_{1}\lambda'_{2}d\bm{E}_{\bm{X}_1}(\lambda_1)\nonumber \\
&& \times \left(  \int\limits_{\lambda'_{1} \in \sigma(\bm{Y}_1)}\lambda'_{1} d\bm{E}_{\bm{Y}_1}(\lambda'_{1}) \right)d\bm{E}_{\bm{X}_2}(\lambda_2)\left(  \int\limits_{\lambda'_{2} \in \sigma(\bm{Y}_2)}\lambda'_{2} d\bm{E}_{\bm{Y}_2}(\lambda'_{2}) \right)d\bm{E}_{\bm{X}_3}(\lambda_3) \nonumber \\
&&+\int\limits_{\lambda_1 \in \sigma(\bm{X}_1)}\int\limits_{\lambda_2 \in \sigma(\bm{X}_2)}\int\limits_{\lambda_3 \in \sigma(\bm{X}_3)}\beta(\lambda_{1}, \lambda_{2}, \lambda_{3})\lambda'_{1}d\bm{E}_{\bm{X}_1}(\lambda_1)\nonumber \\
&& \times \left(  \int\limits_{\lambda'_{1} \in \sigma(\bm{Y}_1)}\lambda'_{1} d\bm{E}_{\bm{Y}_1}(\lambda'_{1}) \right)d\bm{E}_{\bm{X}_2}(\lambda_2)\left(  \int\limits_{\lambda'_{2} \in \sigma(\bm{Y}_2)}\left(\bm{Y}_2-\lambda'_{2}\bm{I}\right)d\bm{E}_{\bm{Y}_2}(\lambda'_{2}) \right) d\bm{E}_{\bm{X}_3}(\lambda_3)   \nonumber \\
&& +\int\limits_{\lambda_1 \in \sigma(\bm{X}_1)}\int\limits_{\lambda_2 \in \sigma(\bm{X}_2)}\int\limits_{\lambda_3 \in \sigma(\bm{X}_3)}\beta(\lambda_{1}, \lambda_{2}, \lambda_{3})\lambda'_{2}d\bm{E}_{\bm{X}_1}(\lambda_1)\nonumber \\
&& \times \left(  \int\limits_{\lambda'_{1} \in \sigma(\bm{Y}_1)}\left(\bm{Y}_1-\lambda'_{1}\bm{I}\right)d\bm{E}_{\bm{Y}_1}(\lambda'_{1}) \right)d\bm{E}_{\bm{X}_2}(\lambda_2)\left(  \int\limits_{\lambda'_{2} \in \sigma(\bm{Y}_2)}\lambda'_{2} d\bm{E}_{\bm{Y}_2}(\lambda'_{2}) \right)d\bm{E}_{\bm{X}_3}(\lambda_3) \nonumber \\
&&+\int\limits_{\lambda_1 \in \sigma(\bm{X}_1)}\int\limits_{\lambda_2 \in \sigma(\bm{X}_2)}\int\limits_{\lambda_3 \in \sigma(\bm{X}_3)}\beta(\lambda_{1}, \lambda_{2}, \lambda_{3})d\bm{E}_{\bm{X}_1}(\lambda_1)\nonumber \\
&& \times \left(  \int\limits_{\lambda'_{1} \in \sigma(\bm{Y}_1)}\left(\bm{Y}_1-\lambda'_{1}\bm{I}\right)d\bm{E}_{\bm{Y}_1}(\lambda'_{1}) \right)d\bm{E}_{\bm{X}_2}(\lambda_2)\nonumber \\
&& \times \left(  \int\limits_{\lambda'_{2} \in \sigma(\bm{Y}_2)}\left(\bm{Y}_2-\lambda'_{2}\bm{I}\right)d\bm{E}_{\bm{Y}_2}(\lambda'_{2})\right)d\bm{E}_{\bm{X}_3}(\lambda_3). 
\end{eqnarray}
By comparing Eq.~\eqref{eq1:exp:TOI is a special case} and Eq.~\eqref{eq2:exp:TOI is a special case}, we have
\begin{eqnarray}
f(\bm{X}_1,\bm{Y}_1,\bm{X}_2,\bm{Y}_2,\bm{X}_3)&=&T_{\beta}^{\bm{X}_1,\bm{X}_2,\bm{X}_3}(\bm{Y}_1, \bm{Y}_2).
\end{eqnarray}
\end{example}

\section{Generalized Multiple Operator Integrals for Continuous Spectrum Operators and Their Algebraic Properties}\label{sec:Generalized Multiple Operator Integrals for Continuous Spectrum Operators and Their Algebraic Properties}

The definition of GMOIs for continuous spectrum is provided by Section~\ref{sec: Generalized Multiple Operator Integrals for Continuous Spectrum Operators} first. Then, the algebraic properties of GMOIs are given by Section~\ref{sec: Algebraic Properties of Generalized Multiple Operator Integrals for Continuous Spectrum Operators}

\subsection{Generalized Multiple Operator Integrals for Continuous Spectrum Operators}\label{sec: Generalized Multiple Operator Integrals for Continuous Spectrum Operators}

According to Theorem~\ref{thm: Spectral Mapping Theorem for r Variables inf}, GMOIs can be defined as follows.
\begin{eqnarray}\label{eq1:GMOI def}
\lefteqn{T_{\beta}^{\bm{X}_1,\ldots,\bm{X}_{\zeta+1}}(\bm{Y}_1,\ldots,\bm{Y}_{\zeta})}\nonumber \\
&\define& \int\limits_{\lambda_1 \in \sigma(\bm{X}_1)}\cdots\int\limits_{\lambda_{\zeta+1}\in \sigma(\bm{X}_{\zeta+1})}\beta(\lambda_{1},\ldots,\lambda_{\zeta+1})d\bm{E}_{\bm{X}_1}(\lambda_1)\bm{Y}_1d\bm{E}_{\bm{X}_2}(\lambda_2)\bm{Y}_2\ldots\bm{Y}_{\zeta}d\bm{E}_{\bm{X}_{\zeta+1}}(\lambda_{\zeta+1})\nonumber \\
&&+ \int\limits_{\lambda_1 \in \sigma(\bm{X}_1)}\cdots\int\limits_{\lambda_{\zeta+1}\in \sigma(\bm{X}_{\zeta+1})}\sum\limits_{\kappa=1}^{\zeta}\sum\limits_{\varrho_\kappa(q_1,\ldots,q_{\zeta+1})}\Bigg(\sum\limits_{\varrho_{\kappa}(q_1,\ldots,q_{\zeta+1})=1}^{m_{\lambda_{\mbox{Ind}(\varrho_{\kappa}(q_1,\ldots,q_{\zeta+1}))}}-1}\nonumber \\
&&~~~~~ \frac{\beta^{\varrho_{\kappa}(q_1,\ldots,q_{\zeta+1})}(\lambda_1,\ldots,\lambda_{k_{\zeta+1}})}{q_{\iota_1}!q_{\iota_2}!\ldots q_{\iota_\kappa}!}\times \prod\limits_{\substack{\varsigma =\mbox{Ind}(\varrho_{\kappa}(q_1,\ldots,q_r)), \bm{Z}_\varsigma=\left(\bm{X}_{\varsigma} - \lambda_{\varsigma}\bm{I}\right)^{q_{\varsigma}}d\bm{E}_{\bm{X}_{\varsigma}}(\lambda_{\varsigma}) \bm{Y}_{\varsigma} \\ \varsigma \neq \mbox{Ind}(\varrho_{\kappa}(q_1,\ldots,q_r)), \bm{Z}_\varsigma=d\bm{E}_{\bm{X}_{\varsigma}}(\lambda_{\varsigma})\bm{Y}_{\varsigma}}
}^{\zeta+1} \bm{Z}_\varsigma\Bigg) 
\nonumber \\
&&+\int\limits_{\lambda_1 \in \sigma(\bm{X}_1)}\cdots\int\limits_{\lambda_{\zeta+1}\in \sigma(\bm{X}_{\zeta+1})}\sum\limits_{q_1=\ldots=q_{\zeta+1}=1}^{m_{\lambda_1}-1,\ldots,m_{\lambda_{\zeta+1}}-1}
\frac{\beta^{(q_1,\ldots,q_{\zeta+1})}(\lambda_1,\ldots,\lambda_{k_{\zeta+1}})}{q_1!\cdots q_{\zeta+1}!}\nonumber \\
&&
~~\times \left(\bm{X}_{1} - \lambda_{1}\bm{I}\right)^{q_{1}}d\bm{E}_{\bm{X}_{1}}(\lambda_{1}) \bm{Y}_1 \left(\bm{X}_{2} - \lambda_{2}\bm{I}\right)^{q_{2}}d\bm{E}_{\bm{X}_{2}}(\lambda_{2}) \bm{Y}_2 \times \nonumber \\
&& \ldots \times \bm{Y}_{\zeta}\left(\bm{X}_{\zeta+1} - \lambda_{\zeta+1}\bm{I}\right)^{q_{\zeta+1}}d\bm{E}_{\bm{X}_{\zeta+1}}(\lambda_{\zeta+1}).
\end{eqnarray}
where $\bm{Y}_{\zeta+1}$ is set as identity operator.

\subsection{Algebraic Properties of Generalized Multiple Operator Integrals for Continuous Spectrum Operators}\label{sec: Algebraic Properties of Generalized Multiple Operator Integrals for Continuous Spectrum Operators}

In this section, we will establish the algebraic properties of the operator $T_{\beta}^{\bm{X}_1,\ldots,\bm{X}_{\zeta+1}}(\bm{Y}_1,\ldots,\bm{Y}_{\zeta})$ defined by Eq.~\eqref{eq1:GMOI def}. 

If the operators $\bm{X}_p$ for $p=1,2,\ldots,\zeta+1$ are decomposed as:
\begin{eqnarray}\label{eq:X1 decomp p and n parts}
\bm{X}_p&=&\int\limits_{\lambda_p \in \sigma(\bm{X}_p)} \lambda_{p}d\bm{E}_{\bm{X}_{p}}(\lambda_p)+
\int\limits_{\lambda_p \in \sigma(\bm{X}_p)}(\bm{X}_p - \lambda_{p}\bm{I})d\bm{E}_{\bm{X}_{p}}(\lambda_p)\nonumber \\
&\define&\bm{X}_{p,P}+\bm{X}_{p,N},
\end{eqnarray}
then, from the definition of $T_{\beta}^{\bm{X}_1,\ldots,\bm{X}_{\zeta+1}}(\bm{Y}_1,\ldots,\bm{Y}_{\zeta})$ provided in Eq.~\eqref{eq1:GMOI def}, we have the following decomposition proposition with respect to parameters operators $\bm{X}_p$ immediately.
\begin{proposition}\label{prop:GMOI decomp by parameters X P X N}
If the operators  $\bm{X}_p$ are decomposed as Eq.~\eqref{eq:X1 decomp p and n parts} for $p=1,2,\ldots,\zeta+1$, we have
\begin{eqnarray}
T_{\beta}^{\bm{X}_1,\ldots,\bm{X}_{\zeta+1}}(\bm{Y}_1,\ldots,\bm{Y}_{\zeta})&=&\sum\limits_{i=1}^{2^{\zeta+1}}T_{\beta}^{[\bm{X}]_{\xi(i)}}(\bm{Y}_1,\ldots,\bm{Y}_{\zeta}),
\end{eqnarray}
where $[\bm{X}]_{\xi(i)}$ is an array with $\zeta+1$ entries such that its $j$-th entry is expressed by
\begin{eqnarray}
([\bm{X}]_{\xi(i)})_j&=& \begin{cases}
     \bm{X}_{j,P}, & \text{if} \left\lfloor \frac{i-1}{2^{j-1}} \right\rfloor \equiv 0 \pmod{2}; \\
     \bm{X}_{j,N}, & \text{otherwise}.
   \end{cases}
\end{eqnarray}
\end{proposition}

For example, if we consider $\zeta=2$, we have GMOI as below:
\begin{eqnarray}\label{eq1:  GTOI def}
\lefteqn{T_{\beta}^{\bm{X}_1,\bm{X}_2,\bm{X}_3}(\bm{Y}_1, \bm{Y}_2)\define\int\limits_{\lambda_1 \in \sigma(\bm{X}_1)}\int\limits_{\lambda_2 \in \sigma(\bm{X}_2)}\int\limits_{\lambda_3 \in \sigma(\bm{X}_3)}\beta(\lambda_{1}, \lambda_{2}, \lambda_{3})d\bm{E}_{\bm{X}_1}(\lambda_1)\bm{Y}_1d\bm{E}_{\bm{X}_2}(\lambda_2)\bm{Y}_2 d\bm{E}_{\bm{X}_3}(\lambda_3)}\nonumber\\
&&+\int\limits_{\lambda_1 \in \sigma(\bm{X}_1)}\int\limits_{\lambda_2 \in \sigma(\bm{X}_2)}\int\limits_{\lambda_3 \in \sigma(\bm{X}_3)}\sum_{q_3=1}^{m_{\lambda_3}-1}\frac{\beta^{(-,-,q_3)}(\lambda_1,\lambda_2,\lambda_3)}{q_3!} \nonumber \\
&&\times d\bm{E}_{\bm{X}_1}(\lambda_1)\bm{Y}_1d\bm{E}_{\bm{X}_2}(\lambda_2)\bm{Y}_2(\bm{X}_3-\lambda_3\bm{I})^{q_3}d\bm{E}_{\bm{X}_3}(\lambda_3)\nonumber \\
&&+\int\limits_{\lambda_1 \in \sigma(\bm{X}_1)}\int\limits_{\lambda_2 \in \sigma(\bm{X}_2)}\int\limits_{\lambda_3 \in \sigma(\bm{X}_3)}\sum_{q_2=1}^{m_{\lambda_2}-1}\frac{\beta^{(-,q_2,-)}(\lambda_1,\lambda_2,\lambda_3)}{q_2!}d\bm{E}_{\bm{X}_1}(\lambda_1)\bm{Y}_1(\bm{X}_2-\lambda_2\bm{I})^{q_2}\nonumber \\
&&\times d\bm{E}_{\bm{X}_2}(\lambda_2)\bm{Y}_2d\bm{E}_{\bm{X}_3}(\lambda_3) \nonumber \\
&&+\int\limits_{\lambda_1 \in \sigma(\bm{X}_1)}\int\limits_{\lambda_2 \in \sigma(\bm{X}_2)}\int\limits_{\lambda_3 \in \sigma(\bm{X}_3)}\sum_{q_1=1}^{m_{\lambda_1}-1}\frac{\beta^{(q_1,-,-)}(\lambda_1,\lambda_2,\lambda_3)}{q_1!}(\bm{X}_1-\lambda_1\bm{I})^{q_1}d\bm{E}_{\bm{X}_1}(\lambda_1)\bm{Y}_1 \nonumber \\
&&\times d\bm{E}_{\bm{X}_2}(\lambda_2)\bm{Y}_2d\bm{E}_{\bm{X}_3}(\lambda_3) \nonumber \\
&&+\int\limits_{\lambda_1 \in \sigma(\bm{X}_1)}\int\limits_{\lambda_2 \in \sigma(\bm{X}_2)}\int\limits_{\lambda_3 \in \sigma(\bm{X}_3)}\sum_{q_2=1}^{m_{\lambda_2}-1}\sum_{q_3=1}^{m_{\lambda_3}-1}\frac{\beta^{(-,q_2,q_3)}(\lambda_1,\lambda_2,\lambda_3)}{q_2! q_3!}d\bm{E}_{\bm{X}_1}(\lambda_1)\bm{Y}_1(\bm{X}_2-\lambda_2\bm{I})^{q_2}\nonumber \\
&&\times d\bm{E}_{\bm{X}_2}(\lambda_2)\bm{Y}_2(\bm{X}_3-\lambda_3\bm{I})^{q_3}d\bm{E}_{\bm{X}_3}(\lambda_3)\nonumber \\
&&+\int\limits_{\lambda_1 \in \sigma(\bm{X}_1)}\int\limits_{\lambda_2 \in \sigma(\bm{X}_2)}\int\limits_{\lambda_3 \in \sigma(\bm{X}_3)}\sum_{q_1=1}^{m_{\lambda_1}-1}\sum_{q_3=1}^{m_{\lambda_3}-1}\frac{\beta^{(q_1,-,q_3)}(\lambda_1,\lambda_2,\lambda_3)}{q_1!q_3!}(\bm{X}_1-\lambda_1\bm{I})^{q_1}d\bm{E}_{\bm{X}_1}(\lambda_1)\bm{Y}_1 \nonumber \\
&&\times d\bm{E}_{\bm{X}_2}(\lambda_2)\bm{Y}_2(\bm{X}_3-\lambda_3\bm{I})^{q_3}d\bm{E}_{\bm{X}_3}(\lambda_3) \nonumber \\
&&+\int\limits_{\lambda_1 \in \sigma(\bm{X}_1)}\int\limits_{\lambda_2 \in \sigma(\bm{X}_2)}\int\limits_{\lambda_3 \in \sigma(\bm{X}_3)}\sum_{q_1=1}^{m_{\lambda_1}-1}\sum_{q_2=1}^{m_{\lambda_2}-1}\frac{\beta^{(q_1,q_2,-)}(\lambda_1,\lambda_2,\lambda_3)}{q_1!q_2!} \nonumber \\
&&\times (\bm{X}_1-\lambda_1\bm{I})^{q_1}d\bm{E}_{\bm{X}_1}(\lambda_1)\bm{Y}_1(\bm{X}_2-\lambda_2\bm{I})^{q_2}d\bm{E}_{\bm{X}_2}(\lambda_2)\bm{Y}_2d\bm{E}_{\bm{X}_3}(\lambda_3) \nonumber \\
&&+\int\limits_{\lambda_1 \in \sigma(\bm{X}_1)}\int\limits_{\lambda_2 \in \sigma(\bm{X}_2)}\int\limits_{\lambda_3 \in \sigma(\bm{X}_3)}\sum_{q_1=1}^{m_{\lambda_1}-1}\sum_{q_2=1}^{m_{\lambda_2}-1}\sum_{q_3=1}^{m_{\lambda_3}-1}\frac{\beta^{(q_1,q_2,q_3)}(\lambda_1,\lambda_2,\lambda_3)}{q_1!q_2!q_3!}\nonumber \\
&&~~\times(\bm{X}_1-\lambda_1\bm{I})^{q_1}d\bm{E}_{\bm{X}_1}(\lambda_1)\bm{Y}_1(\bm{X}_2-\lambda_2\bm{I})^{q_2}d\bm{E}_{\bm{X}_2}(\lambda_2)\bm{Y}_2(\bm{X}_3-\lambda_3\bm{I})^{q_3}d\bm{E}_{\bm{X}_3}(\lambda_3) \nonumber \\
&=_1& \sum\limits_{i=1}^{2^3}T_{\beta}^{[\bm{X}]_{\xi(i)}}(\bm{Y}_1,\bm{Y}_2),
\end{eqnarray}
where $=_1$ comes from Proposition~\ref{prop:GMOI decomp by parameters X P X N}.

In GDOI, the composition of two GDOIs results in another GDOI. Nevertheless, in the case of GMOI, composing a GMOI with multiple GMOIs corresponding to operators with continuous spectra alters the number of parameter operators in the resulting composed GMOI. We state this behavior formally in Proposition~\ref{prop:comp GMOIs}.
\begin{proposition}\label{prop:comp GMOIs}
Given the following GMOIs for operators with continuous spectra : 
\begin{eqnarray}\label{eq1:prop:comp GMOIs}
T_{f}^{\bm{X}_1,\ldots,\bm{X}_{\zeta+1}}(\bm{Y}_1,\ldots,\bm{Y}_{\zeta}),
\end{eqnarray}
and
\begin{eqnarray}\label{eq2:prop:comp GMOIs}
T_{\beta_i}^{\bm{X}_1,\ldots,\bm{X}_{\zeta+1}}(\bm{Y}_1,\ldots,\bm{Y}_{\zeta}),
\end{eqnarray}
where $i=1,2,\ldots,\zeta$ and $[\bm{Y}]^{\zeta}_1$ is an array of operators given by $[\bm{Y}]_1^{\zeta} \define \bm{Y}_1,\bm{Y}_2,\ldots,\bm{Y}_{\zeta}$. If the operators $\bm{X}_p$ for $p=1,2,\ldots,\zeta+1$ are decomposed as:
\begin{eqnarray}\label{eq2.1:prop:comp GMOIs}
\bm{X}_p&=&\int\limits_{\lambda_p \in \sigma(\bm{X}_p)} \lambda_{p}d\bm{E}_{\bm{X}_{p}}(\lambda_p)+
\int\limits_{\lambda_p \in \sigma(\bm{X}_p)}(\bm{X}_p - \lambda_{p}\bm{I})d\bm{E}_{\bm{X}_{p}}(\lambda_p).
\end{eqnarray}

Then, we have
\begin{eqnarray}\label{eq3:prop:comp GMOIs}
\lefteqn{T_{f \prod\limits_{i=1}^{\zeta}\beta_i}^{\bm{X}_1,[\bm{X}]_1^{\zeta+1},\bm{X}_2,[\bm{X}]_1^{\zeta+1},\ldots,[\bm{X}]_1^{\zeta+1},\bm{X}_{\zeta+1}}( \overbrace{\bm{I},[\bm{Y}]^{\zeta}_1, \bm{I},\ldots,\bm{I},[\bm{Y}]^{\zeta}_1,\bm{I}}^{\mbox{there are $\zeta$ terms of $\bm{I},[\bm{Y}]^{\zeta}_1, \bm{I}$}} )}\nonumber \\
&=&T_{f}^{\bm{X}_1,\ldots,\bm{X}_{\zeta+1}}(T_{\beta_1}^{\bm{X}_1,\ldots,\bm{X}_{\zeta+1}}(\bm{Y}_1,\ldots,\bm{Y}_{\zeta}),\ldots,T_{\beta_\zeta}^{\bm{X}_1,\ldots,\bm{X}_{\zeta+1}}(\bm{Y}_1,\ldots,\bm{Y}_{\zeta})),
\end{eqnarray}
where $[\bm{X}]_1^{\zeta+1}$ is an array of operators given by $[\bm{X}]_1^{\zeta+1} \define \bm{X}_1,\bm{X}_2,\ldots,\bm{X}_{\zeta+1}$. Note that there are $(\zeta+1)^2$ parameter operators in GMOI given by L.H.S. of Eq.~\eqref{eq3:prop:comp GMOIs}
\end{proposition}
\textbf{Proof:}
From GMOI for operators with continuous spectra definition given by Eq.~\eqref{eq1:GMOI def}, we have 
\begin{eqnarray}\label{eq4:prop:comp GMOIs}
\lefteqn{T_{f}^{\bm{X}_1,\ldots,\bm{X}_{\zeta+1}}(T_{\beta_1}^{\bm{X}_1,\ldots,\bm{X}_{\zeta+1}}(\bm{Y}_1,\ldots,\bm{Y}_{\zeta}),\ldots,T_{\beta_\zeta}^{\bm{X}_1,\ldots,\bm{X}_{\zeta+1}}(\bm{Y}_1,\ldots,\bm{Y}_{\zeta}))}\nonumber \\
&=&  \int\limits_{\lambda_1 \in \sigma(\bm{X}_1)}\cdots\int\limits_{\lambda_{\zeta+1}\in \sigma(\bm{X}_{\zeta+1})}
f(\lambda_1,\ldots,\lambda_{k_{\zeta+1}})d\bm{E}_{\bm{X}_1}(\lambda_1)T_{\beta_1}^{\bm{X}_1,\ldots,\bm{X}_{\zeta+1}}(\bm{Y}_1,\ldots,\bm{Y}_{\zeta})\nonumber \\
&&~~\times d\bm{E}_{\bm{X}_2}(\lambda_2)T_{\beta_2}^{\bm{X}_1,\ldots,\bm{X}_{\zeta+1}}(\bm{Y}_1,\ldots,\bm{Y}_{\zeta}) \ldots T_{\beta_\zeta}^{\bm{X}_1,\ldots,\bm{X}_{\zeta+1}}(\bm{Y}_1,\ldots,\bm{Y}_{\zeta}) d\bm{E}_{\bm{X}_{\zeta+1}}(\lambda_{\zeta+1})\nonumber \\
&&+\int\limits_{\lambda_1 \in \sigma(\bm{X}_1)}\cdots\int\limits_{\lambda_{\zeta+1}\in \sigma(\bm{X}_{\zeta+1})}\sum\limits_{\kappa=1}^{\zeta}\sum\limits_{\varrho_\kappa(q_1,\ldots,q_{\zeta+1})}\Bigg(\sum\limits_{\varrho_{\kappa}(q_1,\ldots,q_{\zeta+1})=1}^{m_{\lambda_{\mbox{Ind}(\varrho_{\kappa}(q_1,\ldots,q_{\zeta+1}))}}-1}\nonumber \\
&&~~~~~ \frac{f^{\varrho_{\kappa}(q_1,\ldots,q_{\zeta+1})}(\lambda_1,\ldots,\lambda_{k_{\zeta+1}})}{q_{\iota_1}!q_{\iota_2}!\ldots q_{\iota_\kappa}!} \nonumber \\
&& \times \prod\limits_{\substack{\varsigma =\mbox{Ind}(\varrho_{\kappa}(q_1,\ldots,q_r)), \bm{Z}_\varsigma= \left(\bm{X}_{\varsigma} - \lambda_{\varsigma}\bm{I}\right)^{q_{\varsigma}}d\bm{E}_{\bm{X}_{\varsigma}}(\lambda_{\varsigma}) T_{\beta_\varsigma}^{\bm{X}_1,\ldots,\bm{X}_{\zeta+1}}(\bm{Y}_1,\ldots,\bm{Y}_{\zeta})\\ \varsigma \neq \mbox{Ind}(\varrho_{\kappa}(q_1,\ldots,q_r)), \bm{Z}_\varsigma=d \bm{E}_{\bm{X}_\varsigma}(\lambda_{\varsigma}) T_{\beta_\varsigma}^{\bm{X}_1,\ldots,\bm{X}_{\zeta+1}}(\bm{Y}_1,\ldots,\bm{Y}_{\zeta})}
}^{\zeta+1} \bm{Z}_\varsigma\Bigg) 
\nonumber \\
&&+\int\limits_{\lambda_1 \in \sigma(\bm{X}_1)}\cdots\int\limits_{\lambda_{\zeta+1}\in \sigma(\bm{X}_{\zeta+1})}  \sum\limits_{q_1=\ldots=q_{\zeta+1}=1}^{m_{\lambda_1}-1,\ldots,m_{\lambda_{\zeta+1}}-1}
\nonumber \\
&&
\frac{f^{(q_1,\ldots,q_{\zeta+1})}(\lambda_1,\ldots,\lambda_{\zeta+1})}{q_1!\cdots q_{\zeta+1}!}\left(\bm{X}_{1} - \lambda_{1}\bm{I}\right)^{q_{1}}d\bm{E}_{\bm{X}_{1}}(\lambda_{1}) \nonumber \\
&& ~~\times  T_{\beta_1}^{\bm{X}_1,\ldots,\bm{X}_{\zeta+1}}(\bm{Y}_1,\ldots,\bm{Y}_{\zeta}) \left(\bm{X}_{2} - \lambda_{2}\bm{I}\right)^{q_{2}}d\bm{E}_{\bm{X}_{2}}(\lambda_{2}) T_{\beta_2}^{\bm{X}_1,\ldots,\bm{X}_{\zeta+1}}(\bm{Y}_1,\ldots,\bm{Y}_{\zeta})\ldots\nonumber \\
&& ~~ \times \left(\bm{X}_{\zeta} - \lambda_{\zeta}\bm{I}\right)^{q_{\zeta}}d\bm{E}_{\bm{X}_{\zeta}}(\lambda_{\zeta}) T_{\beta_\zeta}^{\bm{X}_1,\ldots,\bm{X}_{\zeta+1}}(\bm{Y}_1,\ldots,\bm{Y}_{\zeta}) \left(\bm{X}_{\zeta+1} - \lambda_{\zeta+1}\bm{I}\right)^{q_{\zeta+1}}d\bm{E}_{\bm{X}_{\zeta+1}}(\lambda_{\zeta+1}),
\end{eqnarray}
where we set $T_{\beta_{\zeta+1}}^{\bm{X}_1,\ldots,\bm{X}_{\zeta+1}}(\bm{Y}_1,\ldots,\bm{Y}_{\zeta})$ as the identity operator.

By observing each term in Eq.~\eqref{eq4:prop:comp GMOIs}, we have the following format regarding to parameter operators
\begin{eqnarray}\label{eq5:prop:comp GMOIs}
\left\{
\begin{array}{l}
d\bm{E}_{\bm{X}_1}(\lambda_1) \\
\left(\bm{X}_{1} - \lambda_{1}\bm{I}\right)^{q_{1}}d\bm{E}_{\bm{X}_{1}}(\lambda_{1})
\end{array}
\right\}, [\bm{X}]_1^{\zeta+1}, \left\{
\begin{array}{l}
d\bm{E}_{\bm{X}_2}(\lambda_2) \\
\left(\bm{X}_{2} - \lambda_{2}\bm{I}\right)^{q_{2}}d\bm{E}_{\bm{X}_{2}}(\lambda_{2})
\end{array}
\right\}, [\bm{X}]_1^{\zeta+1} \nonumber \\
\ldots  [\bm{X}]_1^{\zeta+1}, 
\left\{
\begin{array}{l}
d\bm{E}_{\bm{X}_{\zeta+1}}(\lambda_{\zeta+1}) \\
\left(\bm{X}_{\zeta+1} - \lambda_{\zeta+1}\bm{I}\right)^{q_{\zeta+1}}d\bm{E}_{\bm{X}_{\zeta+1}}(\lambda_{\zeta+1})
\end{array}
\right\},
\end{eqnarray}
where each $\left\{
\begin{array}{l}
d\bm{E}_{\bm{X}_{p}}(\lambda_{p}) \\
\left(\bm{X}_{p} - \lambda_{p}\bm{I}\right)^{q_{p}}d\bm{E}_{\bm{X}_{p}}(\lambda_{p})
\end{array}
\right\}$
is obtained by the projection or nilpotent parts from the operator $\bm{X}_p$. From Eq.\eqref{eq5:prop:comp GMOIs}, we have the argument operators $[\bm{Y}]^{\zeta}_1$ corresponding to each occurrence of $[\bm{X}]_1^{\zeta+1}$. To align Eq.\eqref{eq5:prop:comp GMOIs} with each summand term in Eq.\eqref{eq4:prop:comp GMOIs}, we insert an identity operator both before and after each $[\bm{X}]_1^{\zeta+1}$. 

Accordingly, the full sequence of argument operators can be given by:
\begin{eqnarray}
\overbrace{\bm{I},[\bm{Y}]^{\zeta}_1, \bm{I},\ldots,\bm{I},[\bm{Y}]^{\zeta}_1,\bm{I}}^{\mbox{there are $\zeta$ terms of $\bm{I},[\bm{Y}]^{\zeta}_1, \bm{I}$}}.
\end{eqnarray}
Finally, since the complex-valued functions $f$ and $\beta_i$ are commutative, Eq.\eqref{eq4:prop:comp GMOIs} becomes:
\begin{eqnarray}
T_{f \prod\limits_{i=1}^{\zeta}\beta_i}^{\bm{X}_1,[\bm{X}]_1^{\zeta+1},\bm{X}_2,[\bm{X}]_1^{\zeta+1},\ldots,[\bm{X}]1^{\zeta+1},\bm{X}{\zeta+1}}( \overbrace{\bm{I},[\bm{Y}]^{\zeta}_1, \bm{I},\ldots,\bm{I},[\bm{Y}]^{\zeta}_1,\bm{I}}^{\mbox{there are $\zeta$ terms of $\bm{I},[\bm{Y}]^{\zeta}_1, \bm{I}$}}).
\end{eqnarray}
$\hfill\Box$

\section{Perturbation Formula for GMOIs with Continuous Spectrum Operators}\label{sec:Perturbation Formula for GMOIs with Continuous Spectrum Operators}

In this section, the perturbation formula of GMOI with Continuous Spectrum Operators is derived.  We begin by defining several new notations for later perturbation formula simpler presentation. We have 

\begin{eqnarray}\label{eq1: pert nimpler notation}
\bm{S}_{c}&=&\begin{cases}
d\bm{E}_{\bm{C}}(\lambda_{c}) \\
\left(\bm{C} - \lambda_{c}\bm{I}\right)^{q_{c}}d\bm{E}_{\bm{C}}(\lambda_{c})
\end{cases};\nonumber \\
\bm{S}_{d}&=&\begin{cases}
d\bm{E}_{\bm{D}}(\lambda_{d}) \\
\left(\bm{D} - \lambda_{d}\bm{I}\right)^{q_{d}}d\bm{E}_{\bm{D}}(\lambda_{d})
\end{cases}.
\end{eqnarray}
where $c,d$ are indices for the operators $\bm{C}$ and $\bm{D}$. Also, we set 
\begin{eqnarray}\label{eq2: pert nimpler notation}
\bm{S}_{p}&=&\begin{cases}
d\bm{E}_{\bm{X}_{p}}(\lambda_{p}) \\
\left(\bm{X}_{p} - \lambda_{p}\bm{I}\right)^{q_{p}}d\bm{E}_{\bm{X}_{p}}(\lambda_{p})
\end{cases};\nonumber \\
\bm{S}_{p'}&=&\begin{cases}
d\bm{E}_{\bm{X}_{p'}}(\lambda_{p'}) \\
\left(\bm{X}_{p'} - \lambda_{p'}\bm{I}\right)^{q_{p'}}d\bm{E}_{\bm{X}_{p'}}(\lambda_{p'})
\end{cases},
\end{eqnarray}
where $p,p' \in \mathbb{N}$ are indices for the operators $\bm{X}_p$ and $\bm{X}_{p'}$. We use the following summation notation:
\begin{eqnarray}\label{eq3: pert nimpler notation}
\int\limits_{\bm{S}_{p}} &=&
\begin{cases}
\int\limits_{\lambda_p \in \sigma(\bm{X}_p)}d\bm{E}_{\bm{X}_{p}}(\lambda_{p}), \mbox{~~if $\bm{S}_{p}=d\bm{E}_{\bm{X}_{p}}(\lambda_{p})$}, \\
\int\limits_{\lambda_p \in \sigma(\bm{X}_p)}\sum\limits_{q_p=1}^{m_{\lambda_p}}\left(\bm{X}_{p} - \lambda_{p}\bm{I}\right)^{q_{p}}d\bm{E}_{\bm{X}_{p}}(\lambda_{p}), \mbox{~~if $\bm{S}_{p}=\left(\bm{X}_{p} - \lambda_{p}\bm{I}\right)^{q_{p}}d\bm{E}_{\bm{X}_{p}}(\lambda_{p})$}, 
\end{cases}
\end{eqnarray}
therefore, if we write $\int\limits_{\bm{S}_{p}, p \in [p_1, p_2]}$, which can be expressed by
\begin{eqnarray}\label{eq4: pert nimpler notation}
\int\limits_{\bm{S}_{p}, p \in [p_1, p_2]} &=&
\int\limits_{\bm{S}_{p_1}}\int\limits_{\bm{S}_{p_1 + 1}}\ldots\int\limits_{\bm{S}_{p_2}},
\end{eqnarray}
where $[p_1, p_2]$ denotes the range of positive integers from $p_1, p_1 + 1, \ldots, p_2$. Note that the expression in Eq.~\eqref{eq4: pert nimpler notation} includes the summation results over all binary combinations of $\bm{S}_{p}$. Therefore, Eq.~\eqref{eq4: pert nimpler notation} yields a total of $2^{p_2 - p_1 + 1}$ summation terms.

Given any differentiable function $f: \mathbb{C}^{\zeta} \rightarrow \mathbb{C}$, we use the following notation to represent partial derivatives with respect to different arguments of the function $f$ and its normalization with respect to differentiation orders:
\begin{eqnarray}\label{eq5: pert nimpler notation}
\frac{f^{([q]_1^{q_\zeta})}([\lambda]_1^{\zeta})}{([q!]_1^{q_\zeta})}&\define&\frac{f^{(q_1,q_2,\ldots,q_\zeta)}(\lambda_1, \lambda_2, \ldots, \lambda_{\zeta})}{q_1 ! q_2 ! \ldots,q_\zeta! }.
\end{eqnarray}

The Lemma~\ref{lma:first-order divided difference identity} (adopted from Lemma 2 in~\cite{chang2025GMOIFinite}) will be given to show that the $k$-order divided difference identity will be valid for the function and its partial derivatives.
\begin{lemma}\label{lma:first-order divided difference identity}
Let \( f : \mathbb{C}^n \to \mathbb{C} \) be a function of class \( C^{k+1+r} \). Let
\[
D^{(r)} f = \frac{\partial^r f}{\partial x_{i_1} \cdots \partial x_{i_r}}
\]
be an \( r \)-th order partial derivative of \( f \). Then the divided difference identity
\[
D^{(r)} f^{[k]}(\dots, \lambda_{\alpha_j}, \dots) - D^{(r)} f^{[k]}(\dots, \mu_{\alpha_j}, \dots) = (\lambda_{\alpha_j} - \mu_{\alpha_j}) D^{(r)} f^{[k+1]}(\dots, \lambda_{\alpha_j}, \mu_{\alpha_j}, \dots)
\]
holds provided that the divided difference is taken with respect to a variable independent of \( x_{i_1}, \dots, x_{i_r} \).
\end{lemma}
\textbf{Proof:}
The proof of this Lemma can be found at Lemma 2 in~\cite{chang2025GMOIFinite}
$\hfill\Box$

We are ready to present Theorem~\ref{thm:GMOI Perturbation Formula} about  perturbation formula for GMOI with continuous spectra.
\begin{theorem}\label{thm:GMOI Perturbation Formula}
Given operators $\bm{X}_p$ decomposed as:
\begin{eqnarray}\label{eq1:thm:GMOI Perturbation Formula}
\bm{X}_p&=&\int\limits_{\lambda_p \in \sigma(\bm{X}_p)} \lambda_{p}d\bm{E}_{\bm{X}_{p}}(\lambda_p)+
\int\limits_{\lambda_p \in \sigma(\bm{X}_p)}(\bm{X}_p - \lambda_{p}\bm{I})d\bm{E}_{\bm{X}_{p}}(\lambda_p).
\end{eqnarray}
where $p=1,2,\ldots,\zeta$, and operators $\bm{C}, \bm{D}$ decomposed as:
\begin{eqnarray}\label{eq1:thm:GMOI Perturbation Formula}
\bm{C}&=&\int\limits_{\lambda_c \in \sigma(\bm{X}_c)} \lambda_{c}d\bm{E}_{\bm{C}}(\lambda_c)+
\int\limits_{\lambda_c \in \sigma(\bm{X}_c)}(\bm{X}_c - \lambda_{c}\bm{I})d\bm{E}_{\bm{C}}(\lambda_c), \nonumber \\ 
\bm{D}&=&\int\limits_{\lambda_d \in \sigma(\bm{X}_d)} \lambda_{d}d\bm{E}_{\bm{D}}(\lambda_d)+
\int\limits_{\lambda_d \in \sigma(\bm{X}_d)}(\bm{X}_d - \lambda_{d}\bm{I})d\bm{E}_{\bm{D}}(\lambda_d).
\end{eqnarray}

We have the following formula for GMOI:
{\tiny
\begin{eqnarray}\label{eq1:thm:GMOI Perturbation Formula}
\lefteqn{T_{\beta^{[\zeta+1]}}^{[\bm{X}]_1^{j-1},\bm{C},\bm{D},[\bm{X}]_j^{\zeta}}([\bm{Y}]_1^{j-1},\bm{C}-\bm{D},[\bm{Y}]_j^{\zeta})=}\nonumber \\
&& \Bigg(T_{\beta^{[\zeta]}}^{[\bm{X}]_1^{j-2},\bm{X}_{j-1,P},\bm{C}_P,\bm{X}_{j,P},[\bm{X}]_{j+1}^{\zeta}}([\bm{Y}]_1^{j-1},[\bm{Y}]_j^{\zeta}) \nonumber \\
&-& T_{\beta^{[\zeta]}}^{[\bm{X}]_1^{j-2},\bm{X}_{j-1,P},\bm{D}_P,\bm{X}_{j,P},[\bm{X}]_{j+1}^{\zeta}}([\bm{Y}]_1^{j-1},[\bm{Y}]_j^{\zeta})\Bigg)_I + \mathfrak{X}^{\bm{X}_{j-1,P},\bm{C}_P,\bm{D}_P,\bm{X}_{j,P}} \nonumber \\
&+&\Bigg( T_{\beta^{[\zeta]}}^{[\bm{X}]_1^{j-2},\bm{X}_{j-1,P},\bm{C}_N,\bm{X}_{j,P},[\bm{X}]_{j+1}^{\zeta}}([\bm{Y}]_1^{j-1},[\bm{Y}]_j^{\zeta} \nonumber \\
&-& T_{\beta^{[\zeta]}}^{[\bm{X}]_1^{j-2},\bm{X}_{j-1,P},\bm{D}_N,\bm{X}_{j,P},[\bm{X}]_{j+1}^{\zeta}}([\bm{Y}]_1^{j-1},[\bm{Y}]_j^{\zeta}))\Bigg)_{II} + \mathfrak{X}^{\bm{X}_{j-1,P},\bm{C}_N,\bm{D}_N,\bm{X}_{j,P}} \nonumber \\
&+& \Bigg( T_{\beta^{[\zeta]}}^{[\bm{X}]_1^{j-2},\bm{X}_{j-1,P},\bm{C}_P,\bm{X}_{j,N},[\bm{X}]_{j+1}^{\zeta}}([\bm{Y}]_1^{j-1},\bm{C},[\bm{Y}]_j^{\zeta}) \nonumber \\
&-& T_{\beta^{[\zeta]}}^{[\bm{X}]_1^{j-2},\bm{X}_{j-1,P},\bm{D}_P,\bm{X}_{j,N},[\bm{X}]_{j+1}^{\zeta}}([\bm{Y}]_1^{j-1},\bm{D},[\bm{Y}]_j^{\zeta})\Bigg)_{III}+ \mathfrak{X}^{\bm{X}_{j-1,P},\bm{C}_P,\bm{D}_P,\bm{X}_{j,N}} \nonumber \\
&+&\Bigg( T_{\beta^{[\zeta]}}^{[\bm{X}]_1^{j-2},\bm{X}_{j-1,P},\bm{C}_N,\bm{X}_{j,N},[\bm{X}]_{j+1}^{\zeta}}([\bm{Y}]_1^{j-1},[\bm{Y}]_j^{\zeta}) \nonumber \\
&-& T_{\beta^{[\zeta]}}^{[\bm{X}]_1^{j-2},\bm{X}_{j-1,P},\bm{D}_N,\bm{X}_{j,N},[\bm{X}]_{j+1}^{\zeta}}([\bm{Y}]_1^{j-1},[\bm{Y}]_j^{\zeta})\Bigg)_{IV}+ \mathfrak{X}^{\bm{X}_{j-1,P},\bm{C}_N,\bm{D}_N,\bm{X}_{j,N}} \nonumber \\
&+& \Bigg( T_{\beta^{[\zeta]}}^{[\bm{X}]_1^{j-2},\bm{X}_{j-1,N},\bm{C}_P,\bm{X}_{j,P},[\bm{X}]_{j+1}^{\zeta}}([\bm{Y}]_1^{j-1},[\bm{Y}]_j^{\zeta}) \nonumber \\
&-& T_{\beta^{[\zeta]}}^{[\bm{X}]_1^{j-2},\bm{X}_{j-1,N},\bm{D}_P,\bm{X}_{j,P},[\bm{X}]_{j+1}^{\zeta}}([\bm{Y}]_1^{j-1},[\bm{Y}]_j^{\zeta})\Bigg)_{V}+ \mathfrak{X}^{\bm{X}_{j-1,N},\bm{C}_P,\bm{D}_P,\bm{X}_{j,P}} \nonumber \\
&+&\Bigg(T_{\beta^{[\zeta]}}^{[\bm{X}]_1^{j-2},\bm{X}_{j-1,N},\bm{C}_N,\bm{X}_{j,P},[\bm{X}]_{j+1}^{\zeta}}([\bm{Y}]_1^{j-1},[\bm{Y}]_j^{\zeta}) \nonumber \\
&-& T_{\beta^{[\zeta]}}^{[\bm{X}]_1^{j-2},\bm{X}_{j-1,N},\bm{D}_N,\bm{X}_{j,P},[\bm{X}]_{j+1}^{\zeta}}([\bm{Y}]_1^{j-1},[\bm{Y}]_j^{\zeta})\Bigg)_{VI} + \mathfrak{X}^{\bm{X}_{j-1,N},\bm{C}_N,\bm{D}_N,\bm{X}_{j,P}} \nonumber \\
&+& \Bigg(T_{\beta^{[\zeta]}}^{[\bm{X}]_1^{j-2},\bm{X}_{j-1,N},\bm{C}_P,\bm{X}_{j,N},[\bm{X}]_{j+1}^{\zeta}}([\bm{Y}]_1^{j-1},[\bm{Y}]_j^{\zeta}) \nonumber \\
&-& T_{\beta^{[\zeta]}}^{[\bm{X}]_1^{j-2},\bm{X}_{j-1,N},\bm{D}_P,\bm{X}_{j,N},[\bm{X}]_{j+1}^{\zeta}}([\bm{Y}]_1^{j-1},[\bm{Y}]_j^{\zeta})\Bigg)_{VII} + \mathfrak{X}^{\bm{X}_{j-1,N},\bm{C}_P,\bm{D}_P,\bm{X}_{j,N}} \nonumber \\
&+&\Bigg(T_{\beta^{[\zeta]}}^{[\bm{X}]_1^{j-2},\bm{X}_{j-1,N},\bm{C}_N,\bm{X}_{j,N},[\bm{X}]_{j+1}^{\zeta}}([\bm{Y}]_1^{j-1},[\bm{Y}]_j^{\zeta}) \nonumber \\
&-& T_{\beta^{[\zeta]}}^{[\bm{X}]_1^{j-2},\bm{X}_{j-1,N},\bm{D}_N,\bm{X}_{j,N},[\bm{X}]_{j+1}^{\zeta}}([\bm{Y}]_1^{j-1},[\bm{Y}]_j^{\zeta})\Bigg)_{VIII} + \mathfrak{X}^{\bm{X}_{j-1,N},\bm{C}_N,\bm{D}_N,\bm{X}_{j,N}}\nonumber \\
&+& T_{\beta^{[\zeta+1]}}^{[\bm{X}]_1^{j-2},\bm{X}_{j-1,P},\bm{C}_N,\bm{D}_N,\bm{X}_{j,P},[\bm{X}]_{j+1}^{\zeta}}([\bm{Y}]_1^{j-1},\bm{C}-\bm{D},[\bm{Y}]_j^{\zeta})\nonumber \\
&+& T_{\beta^{[\zeta+1]}}^{[\bm{X}]_1^{j-2},\bm{X}_{j-1,P},\bm{C}_N,\bm{D}_N,\bm{X}_{j,N},[\bm{X}]_{j+1}^{\zeta}}([\bm{Y}]_1^{j-1},\bm{C}-\bm{D},[\bm{Y}]_j^{\zeta})\nonumber \\
&+& T_{\beta^{[\zeta+1]}}^{[\bm{X}]_1^{j-2},\bm{X}_{j-1,N},\bm{C}_N,\bm{D}_N,\bm{X}_{j,P},[\bm{X}]_{j+1}^{\zeta}}([\bm{Y}]_1^{j-1},\bm{C}-\bm{D},[\bm{Y}]_j^{\zeta})\nonumber \\
&+& T_{\beta^{[\zeta+1]}}^{[\bm{X}]_1^{j-2},\bm{X}_{j-1,N},\bm{C}_N,\bm{D}_N,\bm{X}_{j,N},[\bm{X}]_{j+1}^{\zeta}}([\bm{Y}]_1^{j-1},\bm{C}-\bm{D},[\bm{Y}]_j^{\zeta})\nonumber \\
&=&T_{\beta^{[\zeta]}}^{[\bm{X}]_1^{j-1},\bm{C},[\bm{X}]_{j}^{\zeta}}([\bm{Y}]_1^{j-1},[\bm{Y}]_j^{\zeta}) - T_{\beta^{[\zeta]}}^{[\bm{X}]_1^{j-1},\bm{D},[\bm{X}]_{j}^{\zeta}}([\bm{Y}]_1^{j-1},[\bm{Y}]_j^{\zeta}) + \bar{\mathfrak{X}}([\bm{X}]_1^{j-1},\bm{C},\bm{D},[\bm{X}]_j^{\zeta}) \nonumber \\
&+& T_{\beta^{[\zeta+1]}}^{[\bm{X}]_1^{j-2},\bm{X}_{j-1,P},\bm{C}_N,\bm{D}_N,\bm{X}_{j,P},[\bm{X}]_{j+1}^{\zeta}}([\bm{Y}]_1^{j-1},\bm{C}-\bm{D},[\bm{Y}]_j^{\zeta})\nonumber \\
&+& T_{\beta^{[\zeta+1]}}^{[\bm{X}]_1^{j-2},\bm{X}_{j-1,P},\bm{C}_N,\bm{D}_N,\bm{X}_{j,N},[\bm{X}]_{j+1}^{\zeta}}([\bm{Y}]_1^{j-1},\bm{C}-\bm{D},[\bm{Y}]_j^{\zeta})\nonumber \\
&+& T_{\beta^{[\zeta+1]}}^{[\bm{X}]_1^{j-2},\bm{X}_{j-1,N},\bm{C}_N,\bm{D}_N,\bm{X}_{j,P},[\bm{X}]_{j+1}^{\zeta}}([\bm{Y}]_1^{j-1},\bm{C}-\bm{D},[\bm{Y}]_j^{\zeta})\nonumber \\
&+& T_{\beta^{[\zeta+1]}}^{[\bm{X}]_1^{j-2},\bm{X}_{j-1,N},\bm{C}_N,\bm{D}_N,\bm{X}_{j,N},[\bm{X}]_{j+1}^{\zeta}}([\bm{Y}]_1^{j-1},\bm{C}-\bm{D},[\bm{Y}]_j^{\zeta})
\end{eqnarray}
}
where $\beta^{[\ell]}$ is the $\ell$-th order divided difference of a differentiable function $\beta$,  $2 \leq \zeta \in \mathbb{N}$, and $\mathfrak{X}^{\bm{Z}_1,\bm{Z}_2,\bm{Z}_3,\bm{Z}_4}$ are perturbation correction terms with respect to 
\begin{eqnarray}\label{eq2:thm:GMOI Perturbation Formula}
\lefteqn{T_{\beta^{[\zeta+1]}}^{[\bm{X}]_1^{j-2},\bm{Z}_1, \bm{Z}_2,\bm{Z}_3,\bm{Z}_4, [\bm{X}]_{j+1}^{\zeta}}([\bm{Y}]_1^{j-1},\bm{Z}_2-\bm{Z}_3,[\bm{Y}]_j^{\zeta})=}\nonumber \\
&&T_{\beta^{[\zeta]}}^{[\bm{X}]_1^{j-2},\bm{Z}_1,\bm{Z}_2,\bm{Z}_4,[\bm{X}]_{j+1}^{\zeta}}([\bm{Y}]_1^{j-1},[\bm{Y}]_j^{\zeta}) \nonumber \\
&&-T_{\beta^{[\zeta]}}^{[\bm{X}]_1^{j-2},\bm{Z}_1,\bm{Z}_3,\bm{Z}_4,[\bm{X}]_{j+1}^{\zeta}}([\bm{Y}]_1^{j-1},[\bm{Y}]_j^{\zeta})+ \mathfrak{X}^{\bm{Z}_1,\bm{Z}_2,\bm{Z}_3,\bm{Z}_4},
\end{eqnarray}
and
\begin{eqnarray}\label{eq2:thm:GMOI Perturbation Formula}
\bar{\mathfrak{X}}([\bm{X}]_1^{j-1},\bm{C},\bm{D},[\bm{X}]_j^{\zeta})&=& \mathfrak{X}^{\bm{X}_{j-1,P},\bm{C}_P,\bm{D}_P,\bm{X}_{j,P}} + \mathfrak{X}^{\bm{X}_{j-1,P},\bm{C}_N,\bm{D}_N,\bm{X}_{j,P}} \nonumber \\
&+& \mathfrak{X}^{\bm{X}_{j-1,P},\bm{C}_P,\bm{D}_P,\bm{X}_{j,N}} + \mathfrak{X}^{\bm{X}_{j-1,P},\bm{C}_N,\bm{D}_N,\bm{X}_{j,N}} \nonumber \\
&+&  \mathfrak{X}^{\bm{X}_{j-1,N},\bm{C}_P,\bm{D}_P,\bm{X}_{j,P}} + \mathfrak{X}^{\bm{X}_{j-1,N},\bm{C}_N,\bm{D}_N,\bm{X}_{j,P}} \nonumber \\
&+& \mathfrak{X}^{\bm{X}_{j-1,N},\bm{C}_P,\bm{D}_P,\bm{X}_{j,N}} + \mathfrak{X}^{\bm{X}_{j-1,N},\bm{C}_N,\bm{D}_N,\bm{X}_{j,N}}. 
\end{eqnarray}

Each term related to $\mathfrak{X}^{\bm{Z}_1,\bm{Z}_2,\bm{Z}_3,\bm{Z}_4}$ can be expressed by
\begin{eqnarray}\label{eq1:cross:thm:GMOI Perturbation Formula}
\lefteqn{\mathfrak{X}^{\bm{X}_{j-1,P},\bm{C}_P,\bm{D}_P,\bm{X}_{j,P}}=}\nonumber \\
&&\int\limits_{\substack{S_{p}, p \in [1,j-2], \\
S_{p'}, p' \in [j+1,\zeta]}}\int\limits_{\lambda_{j-1},\lambda_c, \lambda_d, \lambda_j}\sum\limits_{q_c=1,q_d=1}^{m_{\lambda_c},m_{\lambda_d}} \frac{(\beta^{[\zeta+1]}([\lambda_k]_1^{j-2},\lambda_{j-1},\lambda_{c},\lambda_{d},\lambda_{j},[\lambda_k]_{j+1}^{\zeta}))^{([q]_1^{j-2}[0,0,0,0][q]_{j+1}^{\zeta})}}{ ([q!]_1^{j-2}[0!,0!,0!,0!][q!]_{j+1}^{\zeta})  } \nonumber \\
&& \times \left(\prod\limits_{p=1}^{j-2}S_{p}\bm{Y}_p\right)d\bm{E}_{\bm{X}_{j-1}}(\lambda_{j-1})\bm{Y}_{j-1} d\bm{E}_{\bm{C}}(\lambda_{c}) \nonumber \\
&&\times (\left(\bm{C} - \lambda_{c}\bm{I}\right)^{q_{c}}d\bm{E}_{\bm{C}}(\lambda_{c})-\left(\bm{C} - \lambda_{d}\bm{I}\right)^{q_{d}}d\bm{E}_{\bm{D}}(\lambda_{d})) d\bm{E}_{\bm{D}}(\lambda_{d}) \bm{Y}_j \nonumber \\
&&\times d\bm{E}_{\bm{X}_{j}}(\lambda_{j}) \left(\prod\limits_{p'=j+1}^{\zeta}\bm{Y}_{p'}S_{p'}\right),
\end{eqnarray}
where $[\lambda_k]_{\ell_1}^{\ell_2}$ for $\ell_1, \ell_2 \in \mathbb{N}$ represents eigenvalues of $\lambda_{\ell_1}$, $\lambda_{\ell_1 + 1}$, $\cdots$, $\lambda_{\ell_2}$.  
{\small
\begin{eqnarray}\label{eq2:cross:thm:GMOI Perturbation Formula}
\lefteqn{\mathfrak{X}^{\bm{X}_{j-1,P},\bm{C}_N,\bm{D}_N,\bm{X}_{j,P}}=}\nonumber \\
&& \int\limits_{\substack{S_{p}, p \in [1,j-2], \\
S_{p'}, p' \in [j+1,\zeta]}}\int\limits_{\lambda_{j-1},\lambda_c, \lambda_d, \lambda_j}\sum\limits_{q_c=1,q_d=1}^{m_{\lambda_c},m_{\lambda_d}} \frac{(\beta^{[\zeta+1]}([\lambda_k]_1^{j-2},\lambda_{k_{j-1}},\lambda_{c},\lambda_{d},\lambda_{j},[\lambda_k]_{j+1}^{\zeta}))^{([q]_1^{j-2}[0,0,q_d,0][q]_{j+1}^{\zeta})}}{ ([q!]_1^{j-2}[0!,0!,q_d!,0!][q!]_{j+1}^{\zeta})  } \nonumber \\
&& \times \left(\prod\limits_{p=1}^{j-2}S_{p}\bm{Y}_p\right) d\bm{E}_{\bm{X}_{j-1}}(\lambda_{j-1})  \bm{Y}_{j-1} d\bm{E}_{\bm{C}}(\lambda_{c}) \nonumber \\
&& \times (\left(\bm{C} - \lambda_{c}\bm{I}\right)^{q_{c}}d\bm{E}_{\bm{C}}(\lambda_{c}) - \left(\bm{D} - \lambda_{d}\bm{I}\right)^{q_{d}}d\bm{E}_{\bm{D}}(\lambda_{d})) \nonumber \\
&& \times \left(\bm{D} - \lambda_{d}\bm{I}\right)^{q_{d}} d\bm{E}_{\bm{D}}(\lambda_{d}) \bm{Y}_j d\bm{E}_{\bm{X}_{j}}(\lambda_{j}) \left(\prod\limits_{p'=j+1}^{\zeta}\bm{Y}_{p'}S_{p'}\right) \nonumber \\
&+&  \int\limits_{\substack{S_{p}, p \in [1,j-2], \\
S_{p'}, p' \in [j+1,\zeta]}}\int\limits_{\lambda_{j-1},\lambda_c, \lambda_d, \lambda_j}\sum\limits_{q_c=1,q_d=1}^{m_{\lambda_c},m_{\lambda_d}} \frac{(\beta^{[\zeta+1]}([\lambda_k]_1^{j-2},\lambda_{j-1},\lambda_{c},\lambda_{d},\lambda_{j},[\lambda_k]_{j+1}^{\zeta}))^{([q]_1^{j-2}[0,q_c,0,0][q]_{j+1}^{\zeta})}}{ ([q!]_1^{j-2}[0!,q_c!,0!,0!][q!]_{j+1}^{\zeta})  } \nonumber \\
&& \times \left(\prod\limits_{p=1}^{j-2}S_{p}\bm{Y}_p\right) d\bm{E}_{\bm{X}_{j-1}}(\lambda_{j-1}) \bm{Y}_{j-1} \left(\bm{C} - \lambda_{c}\bm{I}\right)^{q_{c}}d\bm{E}_{\bm{C}}(\lambda_{c}) \nonumber \\
&& \times (\left(\bm{C} - \lambda_{c}\bm{I}\right)^{q_{c}}d\bm{E}_{\bm{C}}(\lambda_{c})- \left(\bm{D} - \lambda_{d}\bm{I}\right)^{q_{d}}d\bm{E}_{\bm{D}}(\lambda_{d}) ) \nonumber \\
&& \times d\bm{E}_{\bm{D}}(\lambda_{d})  \bm{Y}_j d\bm{E}_{\bm{X}_{j}}(\lambda_{j}) \left(\prod\limits_{p'=j+1}^{\zeta}\bm{Y}_{p'}S_{p'}\right) \nonumber \\
&+&  \int\limits_{\substack{S_{p}, p \in [1,j-2], \\
S_{p'}, p' \in [j+1,\zeta]}}\int\limits_{\lambda_{j-1},\lambda_c, \lambda_d, \lambda_j}\sum\limits_{q_d=1}^{m_{\lambda_d}}\frac{(\beta^{[\zeta]}([\lambda_k]_1^{j-2},\lambda_{j-1},\lambda_{c},\lambda_{j},[\lambda_k]_{j+1}^{\zeta}))^{([q]_1^{j-2}[0,q_d,0][q]_{j+1}^{\zeta})}}{ ([q!]_1^{j-2}[0!,q_d!,0!][q!]_{j+1}^{\zeta})  } \nonumber \\
&& \times \left(\prod\limits_{p=1}^{j-2}S_{p}\bm{Y}_p\right)d\bm{E}_{\bm{X}_{j-1}}(\lambda_{j-1}) \bm{Y}_{j-1}d\bm{E}_{\bm{C}}(\lambda_{c}) \left(\bm{D} - \lambda_{d}\bm{I}\right)^{q_{d}}\nonumber \\
&& \times d\bm{E}_{\bm{D}}(\lambda_{d}) \bm{Y}_jd\bm{E}_{\bm{X}_{j}}(\lambda_{j})\left(\prod\limits_{p'=j+1}^{\zeta}\bm{Y}_{p'}S_{p'}\right) \nonumber \\  
&-&  \int\limits_{\substack{S_{p}, p \in [1,j-2], \\
S_{p'}, p' \in [j+1,\zeta]}}\int\limits_{\lambda_{j-1},\lambda_c,\lambda_d,\lambda_j}\sum\limits_{q_c=1}^{m_{\lambda_c}} \frac{(\beta^{[\zeta]}([\lambda_k]_1^{j-2},\lambda_{j-1},\lambda_{d},\lambda_{j},[\lambda_k]_{j+1}^{\zeta}))^{([q]_1^{j-2}[0,q_c,0][q]_{j+1}^{\zeta})}}{ ([q!]_1^{j-2}[0!,q_c!,0!][q!]_{j+1}^{\zeta})  } \nonumber \\
&& \times \left(\prod\limits_{p=1}^{j-2}S_{p}\bm{Y}_p\right)d\bm{E}_{\bm{X}_{j-1}}(\lambda_{j-1})\bm{Y}_{j-1} \left(\bm{C} - \lambda_{c}\bm{I}\right)^{q_{c}}d\bm{E}_{\bm{C}}(\lambda_{c}) \nonumber \\
&& \times d\bm{E}_{\bm{D}}(\lambda_{d})\bm{Y}_jd\bm{E}_{\bm{X}_{j}}(\lambda_{j})\left(\prod\limits_{p'=j+1}^{\zeta}\bm{Y}_{p'}S_{p'}\right),
\end{eqnarray}
}
\begin{eqnarray}\label{eq3:cross:thm:GMOI Perturbation Formula}
\lefteqn{\mathfrak{X}^{\bm{X}_{j-1,P},\bm{C}_P,\bm{D}_P,\bm{X}_{j,N}}=}\nonumber \\
&&  \int\limits_{\substack{S_{p}, p \in [1,j-2], \\
S_{p'}, p' \in [j+1,\zeta]}}\int\limits_{\lambda_{j-1},\lambda_c,\lambda_d,\lambda_j}\sum\limits_{q_c=1,q_d=1,q_j=1}^{m_{\lambda_c},m_{\lambda_d},m_{\lambda_j}}\frac{(\beta^{[\zeta+1]}([\lambda_k]_1^{j-2},\lambda_{j-1},\lambda_{c},\lambda_{d},\lambda_{j},[\lambda_k]_{j+1}^{\zeta}))^{([q]_1^{j-2}[0,0,0,q_j][q]_{j+1}^{\zeta})}}{ ([q!]_1^{j-2}[0!,0!,0!,q_j !][q!]_{j+1}^{\zeta})  } \nonumber \\
&& \times \left(\prod\limits_{p=1}^{j-2}S_{p}\bm{Y}_p\right) d\bm{E}_{\bm{X}_{j-1}}(\lambda_{j-1}) \bm{Y}_{j-1}d\bm{E}_{\bm{C}}(\lambda_{c})  \nonumber \\
&& \times (\left(\bm{C} - \lambda_{c}\bm{I}\right)^{q_{c}}d\bm{E}_{\bm{C}}(\lambda_{c}) -\left(\bm{D} - \lambda_{d}\bm{I}\right)^{q_{d}}d\bm{E}_{\bm{D}}(\lambda_{d}))  \nonumber \\
&& \times d\bm{E}_{\bm{D}}(\lambda_{d}) \bm{Y}_j \left(\bm{X}_{j} - \lambda_{j}\bm{I}\right)^{q_{j}}d\bm{E}_{\bm{X}_{j}}(\lambda_{j}) \left(\prod\limits_{p'=j+1}^{\zeta}\bm{Y}_{p'}S_{p'}\right),
\end{eqnarray}
{\small
\begin{eqnarray}\label{eq4:cross:thm:GMOI Perturbation Formula}
\lefteqn{\mathfrak{X}^{\bm{X}_{j-1,P},\bm{C}_N,\bm{D}_N,\bm{X}_{j,N}}=}\nonumber \\
&&  \int\limits_{\substack{S_{p}, p \in [1,j-2], \\
S_{p'}, p' \in [j+1,\zeta]}}\int\limits_{\lambda_{j-1},\lambda_c,\lambda_d,\lambda_j}\sum\limits_{q_c=1,q_d=1,q_j=1}^{m_{\lambda_c},m_{\lambda_d},m_{\lambda_j}}\frac{(\beta^{[\zeta+1]}([\lambda_k]_1^{j-2},\lambda_{j-1},\lambda_{c},\lambda_{d},\lambda_{j},[\lambda_k]_{j+1}^{\zeta}))^{([q]_1^{j-2}[0,0,q_d,q_j][q]_{j+1}^{\zeta})}}{ ([q!]_1^{j-2}[0!,0!,q_d!,q_j!][q!]_{j+1}^{\zeta})  } \nonumber \\
&& \times \left(\prod\limits_{p=1}^{j-2}S_{p}\bm{Y}_p\right) d\bm{E}_{\bm{X}_{j-1}}(\lambda_{j-1}) \bm{Y}_{j-1} d\bm{E}_{\bm{C}}(\lambda_{c}) \nonumber \\
&& \times  (\left(\bm{C} - \lambda_{c}\bm{I}\right)^{q_{c}}d\bm{E}_{\bm{C}}(\lambda_{c})- \left(\bm{D} - \lambda_{d}\bm{I}\right)^{q_{d}}d\bm{E}_{\bm{D}}(\lambda_{d}) ) \left(\bm{D} - \lambda_{d}\bm{I}\right)^{q_{d}}\nonumber \\
&& \times d\bm{E}_{\bm{D}}(\lambda_{d}) \bm{Y}_j \left(\bm{X}_{j} - \lambda_{j}\bm{I}\right)^{q_{j}}d\bm{E}_{\bm{X}_{j}}(\lambda_{j}) \left(\prod\limits_{p'=j+1}^{\zeta}\bm{Y}_{p'}S_{p'}\right) \nonumber \\
&+&  \int\limits_{\substack{S_{p}, p \in [1,j-2], \\
S_{p'}, p' \in [j+1,\zeta]}}\int\limits_{\lambda_{j-1},\lambda_c,\lambda_d,\lambda_j}\sum\limits_{q_c=1,q_d=1,q_j=1}^{m_{\lambda_c},m_{\lambda_d},m_{\lambda_j}}\frac{(\beta^{[\zeta+1]}([\lambda_k]_1^{j-2},\lambda_{j-1},\lambda_{c},\lambda_{d},\lambda_{j},[\lambda_k]_{j+1}^{\zeta}))^{([q]_1^{j-2}[0,q_c,0,q_j][q]_{j+1}^{\zeta})}}{ ([q!]_1^{j-2}[0!,q_c!,0!,q_j!][q!]_{j+1}^{\zeta})  } \nonumber \\
&& \times \left(\prod\limits_{p=1}^{j-2}S_{p}\bm{Y}_p\right) d\bm{E}_{\bm{X}_{j-1}}(\lambda_{j-1}) \bm{Y}_{j-1} \left(\bm{C} - \lambda_{c}\bm{I}\right)^{q_{c}}d\bm{E}_{\bm{C}}(\lambda_{c}) \nonumber \\
&& \times  ( \left(\bm{C} - \lambda_{c}\bm{I}\right)^{q_{c}}d\bm{E}_{\bm{C}}(\lambda_{c}) - \left(\bm{D} - \lambda_{d}\bm{I}\right)^{q_{d}}d\bm{E}_{\bm{D}}(\lambda_{d}) )\nonumber \\
&& \times d\bm{E}_{\bm{D}}(\lambda_{d}) \bm{Y}_j \left(\bm{X}_{j} - \lambda_{j}\bm{I}\right)^{q_{j}}d\bm{E}_{\bm{X}_{j}}(\lambda_{j}) \left(\prod\limits_{p'=j+1}^{\zeta}\bm{Y}_{p'}S_{p'}\right) \nonumber \\
&+&  \int\limits_{\substack{S_{p}, p \in [1,j-2], \\
S_{p'}, p' \in [j+1,\zeta]}}\int\limits_{\lambda_{j-1},\lambda_c,\lambda_j}\sum\limits_{q_d=1,q_j=1}^{m_{\lambda_d},m_{\lambda_j}} \frac{(\beta^{[\zeta]}([\lambda_k]_1^{j-2},\lambda_{j-1},\lambda_{c},\lambda_{j},[\lambda_k]_{j+1}^{\zeta}))^{([q]_1^{j-2}[0,q_d,q_j][q]_{j+1}^{\zeta})}}{ ([q!]_1^{j-2}[0!,q_d!,q_j!][q!]_{j+1}^{\zeta})  } \nonumber \\
&& \times \left(\prod\limits_{p=1}^{j-2}S_{p}\bm{Y}_p\right) d\bm{E}_{\bm{X}_{j-1}}(\lambda_{j-1})  \bm{Y}_{j-1} d\bm{E}_{\bm{C}}(\lambda_{c}) \left(\bm{D} - \lambda_{d}\bm{I}\right)^{q_{d}}  \nonumber \\
&& \times d\bm{E}_{\bm{D}}(\lambda_{d}) \bm{Y}_j \left(\bm{X}_{j} - \lambda_{j}\bm{I}\right)^{q_{j}}d\bm{E}_{\bm{X}_{j}}(\lambda_{j}) \left(\prod\limits_{p'=j+1}^{\zeta}\bm{Y}_{p'}S_{p'}\right) \nonumber \\  
&-&  \int\limits_{\substack{S_{p}, p \in [1,j-2], \\
S_{p'}, p' \in [j+1,\zeta]}}\int\limits_{\lambda_{j-1},\lambda_d,\lambda_j}\sum\limits_{q_c=1,q_j=1}^{m_{\lambda_c},m_{\lambda_j}} \frac{(\beta^{[\zeta]}([\lambda_k]_1^{j-2},\lambda_{j-1},\lambda_{d},\lambda_{j},[\lambda_k]_{j+1}^{\zeta}))^{([q]_1^{j-2}[0,q_c,q_j][q]_{j+1}^{\zeta})}}{ ([q!]_1^{j-2}[0!,q_c!,q_j!][q!]_{j+1}^{\zeta})  } \nonumber \\
&& \times \left(\prod\limits_{p=1}^{j-2}S_{p}\bm{Y}_p\right) d\bm{E}_{\bm{X}_{j-1}}(\lambda_{j-1})  \bm{Y}_{j-1} \left(\bm{C} - \lambda_{c}\bm{I}\right)^{q_{c}}d\bm{E}_{\bm{C}}(\lambda_{c})  \nonumber \\
&& \times d\bm{E}_{\bm{D}}(\lambda_{d}) \bm{Y}_j \left(\bm{X}_{j} - \lambda_{j}\bm{I}\right)^{q_{j}}d\bm{E}_{\bm{X}_{j}}(\lambda_{j}) \left(\prod\limits_{p'=j+1}^{\zeta}\bm{Y}_{p'}S_{p'}\right),
\end{eqnarray}
}
\begin{eqnarray}\label{eq5:cross:thm:GMOI Perturbation Formula}
\lefteqn{\mathfrak{X}^{\bm{X}_{j-1,N},\bm{C}_P,\bm{D}_P,\bm{X}_{j,P}}=}\nonumber \\
&& \int\limits_{\substack{S_{p}, p \in [1,j-2], \\
S_{p'}, p' \in [j+1,\zeta]}}\int\limits_{\lambda_{j-1},\lambda_c,\lambda_d,\lambda_j}\sum\limits_{q_{j-1}=1,q_c=1,q_d=1}^{m_{\lambda_{j-1}}, m_{\lambda_c},m_{\lambda_d}}\frac{(\beta^{[\zeta+1]}([\lambda_k]_1^{j-2},\lambda_{j-1},\lambda_{c},\lambda_{d},\lambda_{j},[\lambda_k]_{j+1}^{\zeta}))^{([q]_1^{j-2}[q_{j-1},0,0,0][q]_{j+1}^{\zeta})}}{ ([q!]_1^{j-2}[q_{j-1}!,0!,0!,0!][q!]_{j+1}^{\zeta})  } \nonumber \\
&& \times \left(\prod\limits_{p=1}^{j-2}S_{p}\bm{Y}_p\right) \left(\bm{X}_{j-1} - \lambda_{j-1}\bm{I}\right)^{q_{j-1}}d\bm{E}_{\bm{X}_{j-1}}(\lambda_{j-1}) \bm{Y}_{j-1}d\bm{E}_{\bm{C}}(\lambda_{c})  \nonumber \\
&& \times (\left(\bm{C} - \lambda_{c}\bm{I}\right)^{q_{c}}d\bm{E}_{\bm{C}}(\lambda_{c})- \left(\bm{X}_{p} - \lambda_{d}\bm{I}\right)^{q_{d}}d\bm{E}_{\bm{D}}(\lambda_{d}))  \nonumber \\
&& \times d\bm{E}_{\bm{D}}(\lambda_{d}) \bm{Y}_j d\bm{E}_{\bm{X}_{j}}(\lambda_{j}) \left(\prod\limits_{p'=j+1}^{\zeta}\bm{Y}_{p'}S_{p'}\right),
\end{eqnarray}
{\small
\begin{eqnarray}\label{eq6:cross:thm:GMOI Perturbation Formula}
\lefteqn{\mathfrak{X}^{\bm{X}_{j-1,N},\bm{C}_N,\bm{D}_N,\bm{X}_{j,P}}=}\nonumber \\
&&  \int\limits_{\substack{S_{p}, p \in [1,j-2], \\
S_{p'}, p' \in [j+1,\zeta]}}\int\limits_{\lambda_{j-1},\lambda_c,\lambda_d,\lambda_j}\sum\limits_{q_{j-1}=1,q_c=1,q_d=1}^{m_{\lambda_{j-1}}, m_{\lambda_c},m_{\lambda_d}} \frac{(\beta^{[\zeta+1]}([\lambda_k]_1^{j-2},\lambda_{j-1},\lambda_{c},\lambda_{d},\lambda_{j},[\lambda_k]_{j+1}^{\zeta}))^{([q]_1^{j-2}[q_{j-1},0,q_d,0][q]_{j+1}^{\zeta})}}{ ([q!]_1^{j-2}[q_{j-1}!,0!,q_d!,0!][q!]_{j+1}^{\zeta})  } \nonumber \\
&& \times \left(\prod\limits_{p=1}^{j-2}S_{p}\bm{Y}_p\right) \left(\bm{X}_{j-1} - \lambda_{j-1}\bm{I}\right)^{q_{j-1}}d\bm{E}_{\bm{X}_{j-1}}(\lambda_{j-1}) \bm{Y}_{j-1} d\bm{E}_{\bm{C}}(\lambda_{c})  \nonumber \\
&& \times ( \left(\bm{C} - \lambda_{c}\bm{I}\right)^{q_{c}}d\bm{E}_{\bm{C}}(\lambda_{c}) - \left(\bm{D} - \lambda_{d}\bm{I}\right)^{q_{d}}d\bm{E}_{\bm{D}}(\lambda_{d}) )  \nonumber \\
&& \times \left(\bm{D} - \lambda_{d}\bm{I}\right)^{q_{d}}d\bm{E}_{\bm{D}}(\lambda_{d}) \bm{Y}_j d\bm{E}_{\bm{X}_{j}}(\lambda_{j})\left(\prod\limits_{p'=j+1}^{\zeta}\bm{Y}_{p'}S_{p'}\right) \nonumber \\
&+&  \int\limits_{\substack{S_{p}, p \in [1,j-2], \\
S_{p'}, p' \in [j+1,\zeta]}}\int\limits_{\lambda_{j-1},\lambda_c,\lambda_d,\lambda_j}\sum\limits_{q_{j-1}=1,q_c=1,q_d=1}^{m_{\lambda_{j-1}}, m_{\lambda_c},m_{\lambda_d}} \frac{(\beta^{[\zeta+1]}([\lambda_k]_1^{j-2},\lambda_{j-1},\lambda_{c},\lambda_{d},\lambda_{j},[\lambda_k]_{j+1}^{\zeta}))^{([q]_1^{j-2}[q_{j-1},q_c,0,0][q]_{j+1}^{\zeta})}}{ ([q!]_1^{j-2}[q_{j-1}!,q_c!,0!,0!][q!]_{j+1}^{\zeta})  } \nonumber \\
&& \times \left(\prod\limits_{p=1}^{j-2}S_{p}\bm{Y}_p\right) \left(\bm{X}_{j-1} - \lambda_{j-1}\bm{I}\right)^{q_{j-1}}d\bm{E}_{\bm{X}_{j-1}}(\lambda_{j-1}) \bm{Y}_{j-1}\left(\bm{C} - \lambda_{c}\bm{I}\right)^{q_{c}}d\bm{E}_{\bm{C}}(\lambda_{c}) \nonumber \\
&& \times  (\left(\bm{C} - \lambda_{c}\bm{I}\right)^{q_{c}}d\bm{E}_{\bm{C}}(\lambda_{c}) - \left(\bm{D} - \lambda_{d}\bm{I}\right)^{q_{d}}d\bm{E}_{\bm{D}}(\lambda_{d})) \nonumber \\
&& \times d\bm{E}_{\bm{X}_d}(\lambda_d) \bm{Y}_j  d\bm{E}_{\bm{X}_j}(\lambda_j) \left(\prod\limits_{p'=j+1}^{\zeta}\bm{Y}_{p'}S_{p'}\right) \nonumber \\
&+&  \int\limits_{\substack{S_{p}, p \in [1,j-2], \\
S_{p'}, p' \in [j+1,\zeta]}}\int\limits_{\lambda_{j-1},\lambda_c,\lambda_d,\lambda_j}\sum\limits_{q_{j-1}=1,q_d=1}^{m_{\lambda_{j-1}},m_{\lambda_d}} \frac{(\beta^{[\zeta]}([\lambda_k]_1^{j-2},\lambda_{j-1},\lambda_{c},\lambda_{j},[\lambda_k]_{j+1}^{\zeta}))^{([q]_1^{j-2}[q_{j-1},q_d,0][q]_{j+1}^{\zeta})}}{ ([q!]_1^{j-2}[q_{j-1}!,q_d!,0!][q!]_{j+1}^{\zeta})  } \nonumber \\
&& \times \left(\prod\limits_{p=1}^{j-2}S_{p}\bm{Y}_p\right) \left(\bm{X}_{j-1} - \lambda_{j-1}\bm{I}\right)^{q_{j-1}}d\bm{E}_{\bm{X}_{j-1}}(\lambda_{j-1}) \bm{Y}_{j-1} d\bm{E}_{\bm{C}}(\lambda_{c}) \nonumber \\
&& \times \left(\bm{D} - \lambda_{d}\bm{I}\right)^{q_{d}}d\bm{E}_{\bm{D}}(\lambda_{d})  \bm{Y}_j d\bm{E}_{\bm{X}_{j}}(\lambda_{j})\left(\prod\limits_{p'=j+1}^{\zeta}\bm{Y}_{p'}S_{p'}\right) \nonumber \\  
&-&  \int\limits_{\substack{S_{p}, p \in [1,j-2], \\
S_{p'}, p' \in [j+1,\zeta]}}\int\limits_{\lambda_{j-1},\lambda_c,\lambda_d,\lambda_j}\sum\limits_{q_{j-1}=1,q_c=1}^{m_{\lambda_{j-1}}, m_{\lambda_c}} \frac{(\beta^{[\zeta]}([\lambda_k]_1^{j-2},\lambda_{j-1},\lambda_{d},\lambda_{j},[\lambda_k]_{j+1}^{\zeta}))^{([q]_1^{j-2}[q_{j-1},q_c,0][q]_{j+1}^{\zeta})}}{ ([q!]_1^{j-2}[q_{j-1}!,q_c!,0!][q!]_{j+1}^{\zeta})  } \nonumber \\
&& \times \left(\prod\limits_{p=1}^{j-2}S_{p}\bm{Y}_p\right) \left(\bm{X}_{j-1} - \lambda_{j-1}\bm{I}\right)^{q_{j-1}}d\bm{E}_{\bm{X}_{j-1}}(\lambda_{j-1}) \bm{Y}_{j-1} \nonumber \\
&& \times \left(\bm{C} - \lambda_{c}\bm{I}\right)^{q_{c}}d\bm{E}_{\bm{C}}(\lambda_{c}) d\bm{E}_{\bm{D}}(\lambda_{d} \bm{Y}_j d\bm{E}_{\bm{X}_{j}}(\lambda_{j} \left(\prod\limits_{p'=j+1}^{\zeta}\bm{Y}_{p'}S_{p'}\right),
\end{eqnarray}
}
\begin{eqnarray}\label{eq7:cross:thm:GMOI Perturbation Formula}
\lefteqn{\mathfrak{X}^{\bm{X}_{j-1,N},\bm{C}_P,\bm{D}_P,\bm{X}_{j,N}}=}\nonumber \\
&& \int\limits_{\substack{S_{p}, p \in [1,j-2], \\
S_{p'}, p' \in [j+1,\zeta]}}\int\limits_{\lambda_{j-1},\lambda_c,\lambda_d,\lambda_j}\sum\limits_{q_{j-1}=1,q_c=1,q_d=1,q_j=1}^{m_{\lambda_{j-1}}, m_{\lambda_c},m_{\lambda_d},m_{\lambda_j}} \nonumber \\
&& \frac{(\beta^{[\zeta+1]}([\lambda_k]_1^{j-2},\lambda_{k_{j-1}},\lambda_{c},\lambda_{d},\lambda_{j},[\lambda_k]_{j+1}^{\zeta}))^{([q]_1^{j-2}[q_{j-1},0,0,q_j][q]_{j+1}^{\zeta})}}{ ([q!]_1^{j-2}[q_{j-1}!,0!,0!,q_j!][q!]_{j+1}^{\zeta})  } \nonumber \\
&& \times \left(\prod\limits_{p=1}^{j-2}S_{p}\bm{Y}_p\right) \left(\bm{X}_{j-1} - \lambda_{j-1}\bm{I}\right)^{q_{j-1}}d\bm{E}_{\bm{X}_{j-1}}(\lambda_{j-1}) \bm{Y}_{j-1} d\bm{E}_{\bm{C}}(\lambda_{c})  \nonumber \\
&& \times (\left(\bm{C} - \lambda_{c}\bm{I}\right)^{q_{c}}d\bm{E}_{\bm{C}}(\lambda_{c})- \left(\bm{D} - \lambda_{d}\bm{I}\right)^{q_{d}}d\bm{E}_{\bm{D}}(\lambda_{d}))\nonumber \\
&& \times d\bm{E}_{\bm{D}}(\lambda_{d}) \bm{Y}_j \left(\bm{X}_{j} - \lambda_{p}\bm{I}\right)^{q_{j}}d\bm{E}_{\bm{X}_{j}}(\lambda_{j}) \left(\prod\limits_{p'=j+1}^{\zeta}\bm{Y}_{p'}S_{p'}\right),
\end{eqnarray}
{\tiny
\begin{eqnarray}\label{eq8:cross:thm:GMOI Perturbation Formula}
\lefteqn{\mathfrak{X}^{\bm{X}_{j-1,N},\bm{C}_N,\bm{D}_N,\bm{X}_{j,N}}=}\nonumber \\
&&   \int\limits_{\substack{S_{p}, p \in [1,j-2], \\
S_{p'}, p' \in [j+1,\zeta]}}\int\limits_{\lambda_{j-1},\lambda_c,\lambda_d,\lambda_j}\sum\limits_{q_{j-1}=1,q_c=1,q_d=1,q_j=1}^{m_{\lambda_{j-1}}, m_{\lambda_c},m_{\lambda_d},m_{\lambda_j}}\nonumber \\
&& \frac{(\beta^{[\zeta+1]}([\lambda_k]_1^{j-2},\lambda_{j-1},\lambda_{c},\lambda_{d},\lambda_{j},[\lambda_k]_{j+1}^{\zeta}))^{([q]_1^{j-2}[q_{j-1},0,q_d,q_j][q]_{j+1}^{\zeta})}}{ ([q!]_1^{j-2}[q_{j-1}!,0!,q_d!,q_j!][q!]_{j+1}^{\zeta})  } \nonumber \\
&& \times \left(\prod\limits_{p=1}^{j-2}S_{p}\bm{Y}_p\right) \left(\bm{X}_{j-1} - \lambda_{j-1}\bm{I}\right)^{q_{j-1}}d\bm{E}_{\bm{X}_{j-1}}(\lambda_{j-1}) \bm{Y}_{j-1} d\bm{E}_{\bm{C}}(\lambda_{c})  \nonumber \\
&& \times ( \left(\bm{C} - \lambda_{c}\bm{I}\right)^{q_{c}}d\bm{E}_{\bm{C}}(\lambda_{c})- \left(\bm{X}_{p} - \lambda_{d}\bm{I}\right)^{q_{d}}d\bm{E}_{\bm{D}}(\lambda_{d})) \left(\bm{D} - \lambda_{d}\bm{I}\right)^{q_{d}}  \nonumber \\
&& \times d\bm{E}_{\bm{D}}(\lambda_{d}) \bm{Y}_j \left(\bm{X}_{j} - \lambda_{j}\bm{I}\right)^{q_{j}}d\bm{E}_{\bm{X}_{j}}(\lambda_{j}) \left(\prod\limits_{p'=j+1}^{\zeta}\bm{Y}_{p'}S_{p'}\right) \nonumber \\
&+&  \int\limits_{\substack{S_{p}, p \in [1,j-2], \\
S_{p'}, p' \in [j+1,\zeta]}}\int\limits_{\lambda_{j-1},\lambda_c,\lambda_d,\lambda_j}\sum\limits_{q_{j-1}=1,q_c=1,q_d=1,q_j=1}^{m_{\lambda_{j-1}}, m_{\lambda_c},m_{\lambda_d},m_{\lambda_j}} \nonumber \\
&&
\frac{(\beta^{[\zeta+1]}([\lambda_k]_1^{j-2},\lambda_{j-1},\lambda_{c},\lambda_{d},\lambda_{j},[\lambda_k]_{j+1}^{\zeta}))^{([q]_1^{j-2}[q_{j-1},q_c,0,q_j][q]_{j+1}^{\zeta})}}{ ([q!]_1^{j-2}[q_{j-1}!,q_c!,0!,q_j!][q!]_{j+1}^{\zeta})  } \nonumber \\
&& \times \left(\prod\limits_{p=1}^{j-2}S_{p}\bm{Y}_p\right) \left(\bm{X}_{j-1} - \lambda_{j-1}\bm{I}\right)^{q_{j-1}}d\bm{E}_{\bm{X}_{j-1}}(\lambda_{j-1}) \bm{Y}_{j-1} \left(\bm{C} - \lambda_{c}\bm{I}\right)^{q_{c}}d\bm{E}_{\bm{C}}(\lambda_{c}) \nonumber \\
&& \times  ( \left(\bm{C} - \lambda_{c}\bm{I}\right)^{q_{c}}d\bm{E}_{\bm{C}}(\lambda_{c}) - \left(\bm{D} - \lambda_{d}\bm{I}\right)^{q_{d}}d\bm{E}_{\bm{D}}(\lambda_{d}) ) \nonumber \\
&& \times  d\bm{E}_{\bm{D}}(\lambda_{d}) \bm{Y}_j \left(\bm{X}_{j} - \lambda_{j}\bm{I}\right)^{q_{j}}d\bm{E}_{\bm{X}_{j}}(\lambda_{j}) \left(\prod\limits_{p'=j+1}^{\zeta}\bm{Y}_{p'}S_{p'}\right) \nonumber \\
&+&  \int\limits_{\substack{S_{p}, p \in [1,j-2], \\
S_{p'}, p' \in [j+1,\zeta]}}\int\limits_{\lambda_{j-1},\lambda_c,\lambda_d,\lambda_j}\sum\limits_{q_{j-1}=1,q_d=1,q_j=1}^{m_{\lambda_{j-1}},m_{\lambda_d},m_{\lambda_j}} \nonumber \\
&&
\frac{(\beta^{[\zeta]}([\lambda_k]_1^{j-2},\lambda_{j-1},\lambda_{c},\lambda_{j},[\lambda_k]_{j+1}^{\zeta}))^{([q]_1^{j-2}[q_{j-1},q_d,q_j][q]_{j+1}^{\zeta})}}{ ([q!]_1^{j-2}[q_{j-1}!,q_d!,q_j!][q!]_{j+1}^{\zeta})  } \nonumber \\
&& \times \left(\prod\limits_{p=1}^{j-2}S_{p}\bm{Y}_p\right) \left(\bm{X}_{j-1} - \lambda_{j-1}\bm{I}\right)^{q_{j-1}}d\bm{E}_{\bm{X}_{j-1}}(\lambda_{j-1}) \bm{Y}_{j-1} d\bm{E}_{\bm{C}}(\lambda_{c})  \nonumber \\
&& \times \left(\bm{D} - \lambda_{d}\bm{I}\right)^{q_{d}}d\bm{E}_{\bm{D}}(\lambda_{d}) \bm{Y}_j \left(\bm{X}_{j} - \lambda_{j}\bm{I}\right)^{q_{j}}d\bm{E}_{\bm{X}_{j}}(\lambda_{j}) \left(\prod\limits_{p'=j+1}^{\zeta}\bm{Y}_{p'}S_{p'}\right) \nonumber \\  
&-& \int\limits_{\substack{S_{p}, p \in [1,j-2], \\
S_{p'}, p' \in [j+1,\zeta]}}\int\limits_{\lambda_{j-1},\lambda_c,\lambda_d,\lambda_j}\sum\limits_{q_{j-1}=1,q_c=1,q_j=1}^{m_{\lambda_{j-1}}, m_{\lambda_c},m_{\lambda_j}} \nonumber \\
&&
\frac{(\beta^{[\zeta]}([\lambda_k]_1^{j-2},\lambda_{j-1},\lambda_{d},\lambda_{j},[\lambda_k]_{j+1}^{\zeta}))^{([q]_1^{j-2}[q_{j-1},q_c,q_j][q]_{j+1}^{\zeta})}}{ ([q!]_1^{j-2}[q_{j-1}!,q_c!,q_j!][q!]_{j+1}^{\zeta})  } \nonumber \\
&& \times \left(\prod\limits_{p=1}^{j-2}S_{p}\bm{Y}_p\right) \left(\bm{X}_{j-1} - \lambda_{j-1}\bm{I}\right)^{q_{j-1}}d\bm{E}_{\bm{X}_{j-1}}(\lambda_{j-1}) \bm{Y}_{j-1} \left(\bm{C} - \lambda_{c}\bm{I}\right)^{q_{c}}  \nonumber \\
&& \times d\bm{E}_{\bm{C}}(\lambda_{c}) d\bm{E}_{\bm{D}}(\lambda_{d}) \bm{Y}_j \left(\bm{X}_{p} - \lambda_{j}\bm{I}\right)^{q_{j}}d\bm{E}_{\bm{X}_{j}}(\lambda_{j}) \left(\prod\limits_{p'=j+1}^{\zeta}\bm{Y}_{p'}S_{p'}\right).
\end{eqnarray}
}
\end{theorem}
\textbf{Proof:}
The proof idea follows the proof in Theorem 3 from~\cite{chang2025GMOIFinite}. From Proposition~\ref{prop:GMOI decomp by parameters X P X N}, we have 
{\small
\begin{eqnarray}\label{eq3:thm:GMOI Perturbation Formula}
\lefteqn{T_{\beta^{[\zeta+1]}}^{[\bm{X}]_1^{j-1},\bm{C},\bm{D},[\bm{X}]_j^{\zeta}}([\bm{Y}]_1^{j-1},\bm{C}-\bm{D},[\bm{Y}]_j^{\zeta})=\bm{O}}\nonumber \\
&+_1& T_{\beta^{[\zeta+1]}}^{[\bm{X}]_1^{j-2},\bm{X}_{j-1,P},\bm{C}_P,\bm{D}_P,\bm{X}_{j,P},[\bm{X}]_{j+1}^{\zeta}}([\bm{Y}]_1^{j-1},\bm{C}-\bm{D},[\bm{Y}]_j^{\zeta})\nonumber \\
&+_2& T_{\beta^{[\zeta+1]}}^{[\bm{X}]_1^{j-2},\bm{X}_{j-1,P},\bm{C}_N,\bm{D}_P,\bm{X}_{j,P},[\bm{X}]_{j+1}^{\zeta}}([\bm{Y}]_1^{j-1},\bm{C}-\bm{D},[\bm{Y}]_j^{\zeta})\nonumber \\
&+_3& T_{\beta^{[\zeta+1]}}^{[\bm{X}]_1^{j-2},\bm{X}_{j-1,P},\bm{C}_P,\bm{D}_N,\bm{X}_{j,P},[\bm{X}]_{j+1}^{\zeta}}([\bm{Y}]_1^{j-1},\bm{C}-\bm{D},[\bm{Y}]_j^{\zeta})\nonumber \\
&+_4& T_{\beta^{[\zeta+1]}}^{[\bm{X}]_1^{j-2},\bm{X}_{j-1,P},\bm{C}_P,\bm{D}_P,\bm{X}_{j,N},[\bm{X}]_{j+1}^{\zeta}}([\bm{Y}]_1^{j-1},\bm{C}-\bm{D},[\bm{Y}]_j^{\zeta}) \nonumber \\
&+_5& T_{\beta^{[\zeta+1]}}^{[\bm{X}]_1^{j-2},\bm{X}_{j-1,P},\bm{C}_N,\bm{D}_P,\bm{X}_{j,N},[\bm{X}]_{j+1}^{\zeta}}([\bm{Y}]_1^{j-1},\bm{C}-\bm{D},[\bm{Y}]_j^{\zeta})\nonumber \\
&+_6& T_{\beta^{[\zeta+1]}}^{[\bm{X}]_1^{j-2},\bm{X}_{j-1,P},\bm{C}_P,\bm{D}_N,\bm{X}_{j,N},[\bm{X}]_{j+1}^{\zeta}}([\bm{Y}]_1^{j-1},\bm{C}-\bm{D},[\bm{Y}]_j^{\zeta})\nonumber \\
&+_7& T_{\beta^{[\zeta+1]}}^{[\bm{X}]_1^{j-2},\bm{X}_{j-1,N},\bm{C}_P,\bm{D}_P,\bm{X}_{j,P},[\bm{X}]_{j+1}^{\zeta}}([\bm{Y}]_1^{j-1},\bm{C}-\bm{D},[\bm{Y}]_j^{\zeta})\nonumber \\
&+_8& T_{\beta^{[\zeta+1]}}^{[\bm{X}]_1^{j-2},\bm{X}_{j-1,N},\bm{C}_N,\bm{D}_P,\bm{X}_{j,P},[\bm{X}]_{j+1}^{\zeta}}([\bm{Y}]_1^{j-1},\bm{C}-\bm{D},[\bm{Y}]_j^{\zeta})\nonumber \\
&+_9& T_{\beta^{[\zeta+1]}}^{[\bm{X}]_1^{j-2},\bm{X}_{j-1,N},\bm{C}_P,\bm{D}_N,\bm{X}_{j,P},[\bm{X}]_{j+1}^{\zeta}}([\bm{Y}]_1^{j-1},\bm{C}-\bm{D},[\bm{Y}]_j^{\zeta})\nonumber \\
&+_{10}&T_{\beta^{[\zeta+1]}}^{[\bm{X}]_1^{j-2},\bm{X}_{j-1,N},\bm{C}_P,\bm{D}_P,\bm{X}_{j,N},[\bm{X}]_{j+1}^{\zeta}}([\bm{Y}]_1^{j-1},\bm{C}-\bm{D},[\bm{Y}]_j^{\zeta})\nonumber \\
&+_{11}& T_{\beta^{[\zeta+1]}}^{[\bm{X}]_1^{j-2},\bm{X}_{j-1,N},\bm{C}_N,\bm{D}_P,\bm{X}_{j,N},[\bm{X}]_{j+1}^{\zeta}}([\bm{Y}]_1^{j-1},\bm{C}-\bm{D},[\bm{Y}]_j^{\zeta})\nonumber \\
&+_{12}& T_{\beta^{[\zeta+1]}}^{[\bm{X}]_1^{j-2},\bm{X}_{j-1,N},\bm{C}_P,\bm{D}_N,\bm{X}_{j,N},[\bm{X}]_{j+1}^{\zeta}}([\bm{Y}]_1^{j-1},\bm{C}-\bm{D},[\bm{Y}]_j^{\zeta})\nonumber \\
&+_{}& T_{\beta^{[\zeta+1]}}^{[\bm{X}]_1^{j-2},\bm{X}_{j-1,P},\bm{C}_N,\bm{D}_N,\bm{X}_{j,P},[\bm{X}]_{j+1}^{\zeta}}([\bm{Y}]_1^{j-1},\bm{C}-\bm{D},[\bm{Y}]_j^{\zeta})\nonumber \\
&+& T_{\beta^{[\zeta+1]}}^{[\bm{X}]_1^{j-2},\bm{X}_{j-1,P},\bm{C}_N,\bm{D}_N,\bm{X}_{j,N},[\bm{X}]_{j+1}^{\zeta}}([\bm{Y}]_1^{j-1},\bm{C}-\bm{D},[\bm{Y}]_j^{\zeta})\nonumber \\
&+& T_{\beta^{[\zeta+1]}}^{[\bm{X}]_1^{j-2},\bm{X}_{j-1,N},\bm{C}_N,\bm{D}_N,\bm{X}_{j,P},[\bm{X}]_{j+1}^{\zeta}}([\bm{Y}]_1^{j-1},\bm{C}-\bm{D},[\bm{Y}]_j^{\zeta})\nonumber \\
&+& T_{\beta^{[\zeta+1]}}^{[\bm{X}]_1^{j-2},\bm{X}_{j-1,N},\bm{C}_N,\bm{D}_N,\bm{X}_{j,N},[\bm{X}]_{j+1}^{\zeta}}([\bm{Y}]_1^{j-1},\bm{C}-\bm{D},[\bm{Y}]_j^{\zeta}), 
\end{eqnarray}
}

By comparing Eq.~\eqref{eq1:thm:GMOI Perturbation Formula} and Eq.~\eqref{eq3:thm:GMOI Perturbation Formula}, we claim we have the following pairs of identitoes:
\begin{eqnarray}\label{eq4-1:thm:GMOI Perturbation Formula}
\mbox{GMOI with continuous spectra after $+_1$} &=&I+\mathfrak{X}^{\bm{X}_{j-1,P},\bm{C}_P,\bm{D}_P,\bm{X}_{j,P}}, 
\end{eqnarray}
\begin{eqnarray}\label{eq4-2:thm:GMOI Perturbation Formula}
 \mbox{GMOIs with continuous spectra after $+_2$ and $+_3$} &=&II+\mathfrak{X}^{\bm{X}_{j-1,P},\bm{C}_N,\bm{D}_N,\bm{X}_{j,P}}, 
\end{eqnarray}
\begin{eqnarray}\label{eq4-3:thm:GMOI Perturbation Formula}
\mbox{GMOI with continuous spectra after $+_4$} &=&III+\mathfrak{X}^{\bm{X}_{j-1,P},\bm{C}_P,\bm{D}_P,\bm{X}_{j,N}},  
\end{eqnarray}
\begin{eqnarray}\label{eq4-4:thm:GMOI Perturbation Formula}
\mbox{GMOIs with continuous spectra after $+_5$ and $+_6$} &=&IV+\mathfrak{X}^{\bm{X}_{j-1,P},\bm{C}_N,\bm{D}_N,\bm{X}_{j,N}},  
\end{eqnarray}
\begin{eqnarray}\label{eq4-5:thm:GMOI Perturbation Formula}
\mbox{GMOI with continuous spectra after $+_7$} &=&V+\mathfrak{X}^{\bm{X}_{j-1,N},\bm{C}_P,\bm{D}_P,\bm{X}_{j,P}},  
\end{eqnarray}
\begin{eqnarray}\label{eq4-6:thm:GMOI Perturbation Formula}
\mbox{GMOIs with continuous spectra  after $+_8$ and $+_9$} &=&VI+\mathfrak{X}^{\bm{X}_{j-1,N},\bm{C}_N,\bm{D}_N,\bm{X}_{j,P}},  
\end{eqnarray}
\begin{eqnarray}\label{eq4-7:thm:GMOI Perturbation Formula}
\mbox{GMOI with continuous spectra after $+_{10}$} &=&VII+\mathfrak{X}^{\bm{X}_{j-1,N},\bm{C}_P,\bm{D}_P,\bm{X}_{j,N}},  
\end{eqnarray}
\begin{eqnarray}\label{eq4-8:thm:GMOI Perturbation Formula}
\mbox{GMOIs with continuous spectra after $+_{11}$ and $+_{12}$} &=&VIII+\mathfrak{X}^{\bm{X}_{j-1,N},\bm{C}_N,\bm{D}_N,\bm{X}_{j,N}},
\end{eqnarray}
We will provide detailed proof steps for identity given by Eq.~\eqref{eq4-1:thm:GMOI Perturbation Formula} and for identity given by Eq.~\eqref{eq4-2:thm:GMOI Perturbation Formula} as other identities from Eq.~\eqref{eq4-3:thm:GMOI Perturbation Formula} to Eq.~\eqref{eq4-8:thm:GMOI Perturbation Formula} can be proved similarly. 

\textbf{Identity  given by Eq.~\eqref{eq4-1:thm:GMOI Perturbation Formula} Proof}

Because we have 
{\tiny
\begin{eqnarray}\label{eq1:eq4-1:thm:GMOI Perturbation Formula Proof}
\lefteqn{T_{\beta^{[\zeta]}}^{[\bm{X}]_1^{j-2},\bm{X}_{j-1,P},\bm{C}_P,\bm{X}_{j,P},[\bm{X}]_{j+1}^{\zeta}}([\bm{Y}]_1^{j-1},[\bm{Y}]_j^{\zeta}) - T_{\beta^{[\zeta]}}^{[\bm{X}]_1^{j-2},\bm{X}_{j-1,P},\bm{D}_P,\bm{X}_{j,P},[\bm{X}]_{j+1}^{\zeta}}([\bm{Y}]_1^{j-1},[\bm{Y}]_j^{\zeta})=}\nonumber \\
&& \int\limits_{\substack{S_{p}, p \in [1,j-2], \\
S_{p'}, p' \in [j+1,\zeta]}}\int\limits_{\lambda_{j-1},\lambda_c,\lambda_j} \frac{(\beta^{[\zeta]}([\lambda_k]_1^{j-2},\lambda_{j-1},\lambda_{c},\lambda_{j},[\lambda_k]_{j+1}^{\zeta}))^{([q]_1^{j-2}[0,0,0][q]_{j+1}^{\zeta})}}{ ([q!]_1^{j-2}[0!,0!,0!][q!]_{j+1}^{\zeta})  } \nonumber \\
&& \times \left(\prod\limits_{p=1}^{j-2}S_{p}\bm{Y}_p\right)d\bm{E}_{\bm{X}_{j-1}}(\lambda_{j-1})\bm{Y}_{j-1}d\bm{E}_{\bm{C}}(\lambda_{c})\bm{Y}_jd\bm{E}_{\bm{X}_{j}}(\lambda_{j})\left(\prod\limits_{p'=j+1}^{\zeta}\bm{Y}_{p'}S_{p'}\right) \nonumber \\
&-&  \int\limits_{\substack{S_{p}, p \in [1,j-2], \\
S_{p'}, p' \in [j+1,\zeta]}}\int\limits_{\lambda_{j-1},\lambda_d,\lambda_j} \frac{(\beta^{[\zeta]}([\lambda_k]_1^{j-2},\lambda_{j-1},\lambda_{d},\lambda_{j},[\lambda_k]_{j+1}^{\zeta}))^{([q]_1^{j-2}[0,0,0][q]_{j+1}^{\zeta})}}{ ([q!]_1^{j-2}[0!,0!,0!][q!]_{j+1}^{\zeta})  } \nonumber \\
&& \times \left(\prod\limits_{p=1}^{j-2}S_{p}\bm{Y}_p\right)d\bm{E}_{\bm{X}_{j-1}}(\lambda_{j-1})\bm{Y}_{j-1}d\bm{E}_{\bm{D}}(\lambda_{d})\bm{Y}_jd\bm{E}_{\bm{X}_{j}}(\lambda_{j})\left(\prod\limits_{p'=j+1}^{\zeta}\bm{Y}_{p'}S_{p'}\right) \nonumber \\
&=_1& \int\limits_{\substack{S_{p}, p \in [1,j-2], \\
S_{p'}, p' \in [j+1,\zeta]}}\int\limits_{\lambda_{j-1},\lambda_c,\lambda_j}  \frac{(\beta^{[\zeta]}([\lambda_k]_1^{j-2},\lambda_{j-1},\lambda_{c},\lambda_{j},[\lambda_k]_{j+1}^{\zeta}))^{([q]_1^{j-2}[0,0,0][q]_{j+1}^{\zeta})}}{ ([q!]_1^{j-2}[0!,0!,0!][q!]_{j+1}^{\zeta})  } \nonumber \\
&& \times \left(\prod\limits_{p=1}^{j-2}S_{p}\bm{Y}_p\right)d\bm{E}_{\bm{X}_{j-1}}(\lambda_{j-1})\bm{Y}_{j-1}d\bm{E}_{\bm{C}}(\lambda_{c})d\bm{E}_{\bm{D}}(\lambda_{d})\bm{Y}_jd\bm{E}_{\bm{X}_{j}}(\lambda_{j})\left(\prod\limits_{p'=j+1}^{\zeta}\bm{Y}_{p'}S_{p'}\right) \nonumber \\
&-& \int\limits_{\substack{S_{p}, p \in [1,j-2], \\
S_{p'}, p' \in [j+1,\zeta]}}\int\limits_{\lambda_{j-1},\lambda_d,\lambda_j}  \frac{(\beta^{[\zeta]}([\lambda_k]_1^{j-2},\lambda_{j-1},\lambda_{d},\lambda_{j},[\lambda_k]_{j+1}^{\zeta}))^{([q]_1^{j-2}[0,0,0][q]_{j+1}^{\zeta})}}{ ([q!]_1^{j-2}[0!,0!,0!][q!]_{j+1}^{\zeta})  } \nonumber \\
&& \times \left(\prod\limits_{p=1}^{j-2}S_{p}\bm{Y}_p\right)d\bm{E}_{\bm{X}_{j-1}}(\lambda_{j-1})\bm{Y}_{j-1}d\bm{E}_{\bm{C}}(\lambda_{c})d\bm{E}_{\bm{D}}(\lambda_{d})\bm{Y}_jd\bm{E}_{\bm{X}_{j}}(\lambda_{j})\left(\prod\limits_{p'=j+1}^{\zeta}\bm{Y}_{p'}S_{p'}\right) \nonumber \\
&=_2&  \int\limits_{\substack{S_{p}, p \in [1,j-2], \\
S_{p'}, p' \in [j+1,\zeta]}}\int\limits_{\lambda_{j-1},\lambda_c,\lambda_d,\lambda_j} \frac{(\beta^{[\zeta+1]}([\lambda_k]_1^{j-2},\lambda_{j-1},\lambda_{c},\lambda_{d},\lambda_{j},[\lambda_k]_{j+1}^{\zeta}))^{([q]_1^{j-2}[0,0,0,0][q]_{j+1}^{\zeta})}}{ ([q!]_1^{j-2}[0!,0!,0!,0!][q!]_{j+1}^{\zeta})  } \nonumber \\
&& \times \left(\prod\limits_{p=1}^{j-2}S_{p}\bm{Y}_p\right)d\bm{E}_{\bm{X}_{j-1}}(\lambda_{j-1})\bm{Y}_{j-1}d\bm{E}_{\bm{C}}(\lambda_{c})(d\bm{E}_{\bm{C}}(\lambda_{c})-d\bm{E}_{\bm{D}}(\lambda_{d}))d\bm{E}_{\bm{D}}(\lambda_{d})\bm{Y}_jd\bm{E}_{\bm{X}_{j}}(\lambda_{j})\left(\prod\limits_{p'=j+1}^{\zeta}\bm{Y}_{p'}S_{p'}\right)\nonumber \\
&=&  T_{\beta^{[\zeta+1]}}^{[\bm{X}]_1^{j-2},\bm{X}_{j-1,P},\bm{C}_P,\bm{D}_P,\bm{X}_{j,P},[\bm{X}]_{j+1}^{\zeta}}([\bm{Y}]_1^{j-1},\bm{C}-\bm{D},[\bm{Y}]_j^{\zeta}) \nonumber \\
&-& \int\limits_{\substack{S_{p}, p \in [1,j-2], \\
S_{p'}, p' \in [j+1,\zeta]}}\int\limits_{\lambda_{j-1},\lambda_c,\lambda_d,\lambda_j} \sum\limits_{q_c=1,q_d=1}^{m_{\lambda_c},m_{\lambda_d}}\nonumber \\
&& \frac{(\beta^{[\zeta+1]}([\lambda_k]_1^{j-2},\lambda_{j-1},\lambda_{c},\lambda_{d},\lambda_{j},[\lambda_k]_{j+1}^{\zeta}))^{([q]_1^{j-2}[0,0,0,0][q]_{j+1}^{\zeta})}}{ ([q!]_1^{j-2}[0!,0!,0!,0!][q!]_{j+1}^{\zeta})  } \nonumber \\
&& \times \left(\prod\limits_{p=1}^{j-2}S_{p}\bm{Y}_p\right)d\bm{E}_{\bm{X}_{j-1}}(\lambda_{j-1})\bm{Y}_{j-1}d\bm{E}_{\bm{C}}(\lambda_{c})\nonumber \\
&& \times (\left(\bm{C} - \lambda_{c}\bm{I}\right)^{q_{c}}d\bm{E}_{\bm{C}}(\lambda_{c})- \left(\bm{D} - \lambda_{d}\bm{I}\right)^{q_{d}}d\bm{E}_{\bm{D}}(\lambda_{d}))d\bm{E}_{\bm{D}}(\lambda_{d})\bm{Y}_jd\bm{E}_{\bm{X}_{j}}(\lambda_{j})\left(\prod\limits_{p'=j+1}^{\zeta}\bm{Y}_{p'}S_{p'}\right),
\end{eqnarray}
}
where we apply $\int\limits_{\lambda_c \in \sigma(\bm{X}_c)}d\bm{E}_{\bm{C}}(\lambda_{c})   = \int\limits_{\lambda_d \in \sigma(\bm{X}_d)}d\bm{E}_{\bm{D}}(\lambda_{d}) = \bm{I}$ in $=_1$,  apply Lemma~\ref{lma:first-order divided difference identity} in $=_2$. Moreover, the perturbation correct term $\mathfrak{X}^{\bm{X}_{j-1,P},\bm{C}_P,\bm{D}_P,\bm{X}_{j,P}}$ can be expressed by
\begin{eqnarray}\label{eq1:eq4-1:thm:GMOI Perturbation Formula Proof}
\lefteqn{\mathfrak{X}^{\bm{X}_{j-1,P},\bm{C}_P,\bm{D}_P,\bm{X}_{j,P}}=}\nonumber \\
&&  \int\limits_{\substack{S_{p}, p \in [1,j-2], \\
S_{p'}, p' \in [j+1,\zeta]}}\int\limits_{\lambda_{j-1},\lambda_c,\lambda_d,\lambda_j} \sum\limits_{q_c=1,q_d=1}^{m_{\lambda_c},m_{\lambda_d}}\nonumber \\
&& 
\frac{(\beta^{[\zeta+1]}([\lambda_k]_1^{j-2},\lambda_{j-1},\lambda_{c},\lambda_{d},\lambda_{j},[\lambda_k]_{j+1}^{\zeta}))^{([q]_1^{j-2}[0,0,0,0][q]_{j+1}^{\zeta})}}{ ([q!]_1^{j-2}[0!,0!,0!,0!][q!]_{j+1}^{\zeta})  } \nonumber \\
&& \times \left(\prod\limits_{p=1}^{j-2}S_{p}\bm{Y}_p\right)d\bm{E}_{\bm{X}_{j-1}}(\lambda_{j-1})\bm{Y}_{j-1}d\bm{E}_{\bm{C}}(\lambda_{c})(\left(\bm{C} - \lambda_{c}\bm{I}\right)^{q_{c}}d\bm{E}_{\bm{C}}(\lambda_{c}) \nonumber \\
&&  - \left(\bm{D} - \lambda_{d}\bm{I}\right)^{q_{d}}d\bm{E}_{\bm{D}}(\lambda_{d})) d\bm{E}_{\bm{D}}(\lambda_{d})\bm{Y}_jd\bm{E}_{\bm{X}_{j}}(\lambda_{j})\left(\prod\limits_{p'=j+1}^{\zeta}\bm{Y}_{p'}S_{p'}\right).
\end{eqnarray}

\textbf{Identity  given by Eq.~\eqref{eq4-2:thm:GMOI Perturbation Formula} Proof}

Because we have 
{\tiny
\begin{eqnarray}\label{eq1:eq4-2:thm:GMOI Perturbation Formula Proof}
\lefteqn{T_{\beta^{[\zeta]}}^{[\bm{X}]_1^{j-2},\bm{X}_{j-1,P},\bm{C}_N,\bm{X}_{j,P},[\bm{X}]_{j+1}^{\zeta}}([\bm{Y}]_1^{j-1},[\bm{Y}]_j^{\zeta}) - T_{\beta^{[\zeta]}}^{[\bm{X}]_1^{j-2},\bm{X}_{j-1,P},\bm{D}_N,\bm{X}_{j,P},[\bm{X}]_{j+1}^{\zeta}}([\bm{Y}]_1^{j-1},[\bm{Y}]_j^{\zeta})=}\nonumber \\
&&\int\limits_{\substack{S_{p}, p \in [1,j-2], \\
S_{p'}, p' \in [j+1,\zeta]}}\int\limits_{\lambda_{j-1},\lambda_c,\lambda_j} \sum\limits_{q_c=1}^{m_{\lambda_c}} \nonumber \\
&& \frac{(\beta^{[\zeta]}([\lambda_k]_1^{j-2},\lambda_{j-1},\lambda_{c},\lambda_{j},[\lambda_k]_{j+1}^{\zeta}))^{([q]_1^{j-2}[0,q_c,0][q]_{j+1}^{\zeta})}}{ ([q!]_1^{j-2}[0!,q_c!,0!][q!]_{j+1}^{\zeta})  } \nonumber \\
&& \times \left(\prod\limits_{p=1}^{j-2}S_{p}\bm{Y}_p\right)d\bm{E}_{\bm{X}_{j-1}}(\lambda_{j-1})\bm{Y}_{j-1} \left(\bm{C} - \lambda_{c}\bm{I}\right)^{q_{c}}d\bm{E}_{\bm{C}}(\lambda_{c}) \bm{Y}_jd\bm{E}_{\bm{X}_{j}}(\lambda_{j})\left(\prod\limits_{p'=j+1}^{\zeta}\bm{Y}_{p'}S_{p'}\right) \nonumber \\
&-&\int\limits_{\substack{S_{p}, p \in [1,j-2], \\
S_{p'}, p' \in [j+1,\zeta]}}\int\limits_{\lambda_{j-1},\lambda_d,\lambda_j} \sum\limits_{q_d=1}^{m_{\lambda_d}} \frac{(\beta^{[\zeta]}([\lambda_k]_1^{j-2},\lambda_{j-1},\lambda_{d},\lambda_{j},[\lambda_k]_{j+1}^{\zeta}))^{([q]_1^{j-2}[0,q_d,0][q]_{j+1}^{\zeta})}}{ ([q!]_1^{j-2}[0!,q_d!,0!][q!]_{j+1}^{\zeta})  } \nonumber \\
&& \times \left(\prod\limits_{p=1}^{j-2}S_{p}\bm{Y}_p\right)d\bm{E}_{\bm{X}_{j-1}}(\lambda_{j-1})\bm{Y}_{j-1} \left(\bm{D} - \lambda_{d}\bm{I}\right)^{q_{d}}d\bm{E}_{\bm{D}}(\lambda_{d}) \bm{Y}_jd\bm{E}_{\bm{X}_{j}}(\lambda_{j})\left(\prod\limits_{p'=j+1}^{\zeta}\bm{Y}_{p'}S_{p'}\right) \nonumber \\
&=_1&\int\limits_{\substack{S_{p}, p \in [1,j-2], \\
S_{p'}, p' \in [j+1,\zeta]}}\int\limits_{\lambda_{j-1},\lambda_c,\lambda_d,\lambda_j} \sum\limits_{q_c=1}^{m_{\lambda_c}} \frac{(\beta^{[\zeta]}([\lambda_k]_1^{j-2},\lambda_{j-1},\lambda_{c},\lambda_{j},[\lambda_k]_{j+1}^{\zeta}))^{([q]_1^{j-2}[0,q_c,0][q]_{j+1}^{\zeta})}}{ ([q!]_1^{j-2}[0!,q_c!,0!][q!]_{j+1}^{\zeta})  } \nonumber \\
&& \times \left(\prod\limits_{p=1}^{j-2}S_{p}\bm{Y}_p\right)d\bm{E}_{\bm{X}_{j-1}}(\lambda_{j-1})\bm{Y}_{j-1}\left(\bm{C} - \lambda_{c}\bm{I}\right)^{q_{c}}d\bm{E}_{\bm{C}}(\lambda_{c})\nonumber \\
&& d\bm{E}_{\bm{D}}(\lambda_{d})\bm{Y}_jd\bm{E}_{\bm{X}_{j}}(\lambda_{j})\left(\prod\limits_{p'=j+1}^{\zeta}\bm{Y}_{p'}S_{p'}\right) \nonumber \\
&-& \int\limits_{\substack{S_{p}, p \in [1,j-2], \\
S_{p'}, p' \in [j+1,\zeta]}}\int\limits_{\lambda_{j-1},\lambda_c,\lambda_d,\lambda_j} \sum\limits_{q_d=1}^{m_{\lambda_d}} \frac{(\beta^{[\zeta]}([\lambda_k]_1^{j-2},\lambda_{j-1},\lambda_{d},\lambda_{j},[\lambda_k]_{j+1}^{\zeta}))^{([q]_1^{j-2}[0,q_d,0][q]_{j+1}^{\zeta})}}{ ([q!]_1^{j-2}[0!,q_d!,0!][q!]_{j+1}^{\zeta})  } \nonumber \\
&& \times \left(\prod\limits_{p=1}^{j-2}S_{p}\bm{Y}_p\right)d\bm{E}_{\bm{X}_{j-1}}(\lambda_{j-1})\bm{Y}_{j-1}d\bm{E}_{\bm{C}}(\lambda_{c})\left(\bm{D} - \lambda_{d}\bm{I}\right)^{q_{d}}d\bm{E}_{\bm{D}}(\lambda_{d})\bm{Y}_jd\bm{E}_{\bm{X}_{j}}(\lambda_{j})\left(\prod\limits_{p'=j+1}^{\zeta}\bm{Y}_{p'}S_{p'}\right),
\end{eqnarray}
}
where we apply $\int\limits_{\lambda_c \in \sigma(\bm{X}_c)}d\bm{E}_{\bm{C}}(\lambda_{c}) = \int\limits_{\lambda_d \in \sigma(\bm{X}_d)}d\bm{E}_{\bm{D}}(\lambda_{d}) = \bm{I}$ in $=_1$.

Different proof in Eq.~\eqref{eq4-1:thm:GMOI Perturbation Formula}, the operators production part of the first and the second terms are different, i.e., $\left(\bm{C} - \lambda_{c}\bm{I}\right)^{q_{c}}d\bm{E}_{\bm{C}}(\lambda_{p}) d\bm{E}_{\bm{D}}(\lambda_{d}) \neq d\bm{E}_{\bm{C}}(\lambda_{c}) \left(\bm{D} - \lambda_{d}\bm{I}\right)^{q_{d}}d\bm{E}_{\bm{D}}(\lambda_{d})$.  Therefore, we will add and subtract the same auxillary terms to make two GMOI terems agree with GMOIs terms after $+_2$ and $+_3$ gievn by Eq.~\eqref{eq3:thm:GMOI Perturbation Formula}. Continue from Eq.~\eqref{eq1:eq4-2:thm:GMOI Perturbation Formula Proof}, we have 
{\tiny
\begin{eqnarray}\label{eq2:eq4-2:thm:GMOI Perturbation Formula Proof}
\lefteqn{T_{\beta^{[\zeta]}}^{[\bm{X}]_1^{j-2},\bm{X}_{j-1,P},\bm{C}_N,\bm{X}_{j,P},[\bm{X}]_{j+1}^{\zeta}}([\bm{Y}]_1^{j-1},[\bm{Y}]_j^{\zeta}) - T_{\beta^{[\zeta]}}^{[\bm{X}]_1^{j-2},\bm{X}_{j-1,P},\bm{D}_N,\bm{X}_{j,P},[\bm{X}]_{j+1}^{\zeta}}([\bm{Y}]_1^{j-1},[\bm{Y}]_j^{\zeta})=}\nonumber \\
&&\Bigg[ \int\limits_{\substack{S_{p}, p \in [1,j-2], \\
S_{p'}, p' \in [j+1,\zeta]}}\int\limits_{\lambda_{j-1},\lambda_c,\lambda_d,\lambda_j} \sum\limits_{q_c=1}^{m_{\lambda_c}} \frac{(\beta^{[\zeta]}([\lambda_k]_1^{j-2},\lambda_{j-1},\lambda_{c},\lambda_{j},[\lambda_k]_{j+1}^{\zeta}))^{([q]_1^{j-2}[0,q_c,0][q]_{j+1}^{\zeta})}}{ ([q!]_1^{j-2}[0!,q_c!,0!][q!]_{j+1}^{\zeta})  } \nonumber \\
&& \times \left(\prod\limits_{p=1}^{j-2}S_{p}\bm{Y}_p\right)d\bm{E}_{\bm{X}_{j-1}}(\lambda_{j-1})\bm{Y}_{j-1}\left(\bm{C} - \lambda_{c}\bm{I}\right)^{q_{c}}d\bm{E}_{\bm{C}}(\lambda_{c})d\bm{E}_{\bm{D}}(\lambda_{d})\bm{Y}_jd\bm{E}_{\bm{X}_{j}}(\lambda_{j})\left(\prod\limits_{p'=j+1}^{\zeta}\bm{Y}_{p'}S_{p'}\right) \nonumber \\
&-& \int\limits_{\substack{S_{p}, p \in [1,j-2], \\
S_{p'}, p' \in [j+1,\zeta]}}\int\limits_{\lambda_{j-1},\lambda_c,\lambda_d,\lambda_j} \sum\limits_{q_c=1}^{m_{\lambda_c}} \frac{(\beta^{[\zeta]}([\lambda_k]_1^{j-2},\lambda_{j-1},\lambda_{d},\lambda_{j},[\lambda_k]_{j+1}^{\zeta}))^{([q]_1^{j-2}[0,q_c,0][q]_{j+1}^{\zeta})}}{ ([q!]_1^{j-2}[0!,q_c!,0!][q!]_{j+1}^{\zeta})  } \nonumber \\
&& \times \left(\prod\limits_{p=1}^{j-2}S_{p}\bm{Y}_p\right)d\bm{E}_{\bm{X}_{j-1}}(\lambda_{j-1})\bm{Y}_{j-1} \left(\bm{C} - \lambda_{c}\bm{I}\right)^{q_{c}}d\bm{E}_{\bm{C}}(\lambda_{c}) d\bm{E}_{\bm{D}}(\lambda_{d})\bm{Y}_jd\bm{E}_{\bm{X}_{j}}(\lambda_{j})\left(\prod\limits_{p'=j+1}^{\zeta}\bm{Y}_{p'}S_{p'}\right) \Bigg]\nonumber \\
&+& \int\limits_{\substack{S_{p}, p \in [1,j-2], \\
S_{p'}, p' \in [j+1,\zeta]}}\int\limits_{\lambda_{j-1},\lambda_c,\lambda_d,\lambda_j} \sum\limits_{q_c=1}^{m_{\lambda_c}}\frac{(\beta^{[\zeta]}([\lambda_k]_1^{j-2},\lambda_{k_{j-1}},\lambda_{d},\lambda_{k_j},[\lambda_k]_{j+1}^{\zeta}))^{([q]_1^{j-2}[0,q_c,0][q]_{j+1}^{\zeta})}}{ ([q!]_1^{j-2}[0!,q_c!,0!][q!]_{j+1}^{\zeta})  } \nonumber \\
&& \times \left(\prod\limits_{p=1}^{j-2}S_{p}\bm{Y}_p\right)d\bm{E}_{\bm{X}_{j-1}}(\lambda_{j-1})\bm{Y}_{j-1} \left(\bm{C} - \lambda_{c}\bm{I}\right)^{q_{c}}d\bm{E}_{\bm{C}}(\lambda_{c}) d\bm{E}_{\bm{D}}(\lambda_{d})\bm{Y}_jd\bm{E}_{\bm{X}_{j}}(\lambda_{j})\left(\prod\limits_{p'=j+1}^{\zeta}\bm{Y}_{p'}S_{p'}\right) \nonumber \\
&-&\bigg[\int\limits_{\substack{S_{p}, p \in [1,j-2], \\
S_{p'}, p' \in [j+1,\zeta]}}\int\limits_{\lambda_{j-1},\lambda_c,\lambda_d,\lambda_j} \sum\limits_{q_c=1}^{m_{\lambda_c}} \frac{(\beta^{[\zeta]}([\lambda_k]_1^{j-2},\lambda_{k_{j-1}},\lambda_{d},\lambda_{k_j},[\lambda_k]_{j+1}^{\zeta}))^{([q]_1^{j-2}[0,q_d,0][q]_{j+1}^{\zeta})}}{ ([q!]_1^{j-2}[0!,q_d!,0!][q!]_{j+1}^{\zeta})  } \nonumber \\
&& \times \left(\prod\limits_{p=1}^{j-2}S_{p}\bm{Y}_p\right)d\bm{E}_{\bm{X}_{j-1}}(\lambda_{j-1})\bm{Y}_{j-1}d\bm{E}_{\bm{C}}(\lambda_{c}) \left(\bm{D} - \lambda_{d}\bm{I}\right)^{q_{d}}d\bm{E}_{\bm{D}}(\lambda_{d})  \bm{Y}_jd\bm{E}_{\bm{X}_{j}}(\lambda_{j})\left(\prod\limits_{p'=j+1}^{\zeta}\bm{Y}_{p'}S_{p'}\right)\nonumber \\
&+&\int\limits_{\substack{S_{p}, p \in [1,j-2], \\
S_{p'}, p' \in [j+1,\zeta]}}\int\limits_{\lambda_{j-1},\lambda_c,\lambda_d,\lambda_j} \sum\limits_{q_c=1}^{m_{\lambda_c}} \frac{(\beta^{[\zeta]}([\lambda_k]_1^{j-2},\lambda_{j-1},\lambda_{c},\lambda_{j},[\lambda_k]_{j+1}^{\zeta}))^{([q]_1^{j-2}[0,q_d,0][q]_{j+1}^{\zeta})}}{ ([q!]_1^{j-2}[0!,q_d!,0!][q!]_{j+1}^{\zeta})  } \nonumber \\
&& \times \left(\prod\limits_{p=1}^{j-2}S_{p}\bm{Y}_p\right)d\bm{E}_{\bm{X}_{j-1}}(\lambda_{j-1})\bm{Y}_{j-1}d\bm{E}_{\bm{C}}(\lambda_{c})\left(\bm{D} - \lambda_{d}\bm{I}\right)^{q_{d}}d\bm{E}_{\bm{D}}(\lambda_{d})\bm{Y}_jd\bm{E}_{\bm{X}_{j}}(\lambda_{j})\left(\prod\limits_{p'=j+1}^{\zeta}\bm{Y}_{p'}S_{p'}\right)\Bigg]\nonumber \\
&-&\int\limits_{\substack{S_{p}, p \in [1,j-2], \\
S_{p'}, p' \in [j+1,\zeta]}}\int\limits_{\lambda_{j-1},\lambda_c,\lambda_d,\lambda_j} \sum\limits_{q_c=1}^{m_{\lambda_c}} \frac{(\beta^{[\zeta]}([\lambda_k]_1^{j-2},\lambda_{k_{j-1}},\lambda_{c},\lambda_{k_j},[\lambda_k]_{j+1}^{\zeta}))^{([q]_1^{j-2}[0,q_d,0][q]_{j+1}^{\zeta})}}{ ([q!]_1^{j-2}[0!,q_d!,0!][q!]_{j+1}^{\zeta})  } \nonumber \\
&& \times \left(\prod\limits_{p=1}^{j-2}S_{p}\bm{Y}_p\right)d\bm{E}_{\bm{X}_{j-1}}(\lambda_{j-1})\bm{Y}_{j-1}d\bm{E}_{\bm{C}}(\lambda_{c}) \left(\bm{D} - \lambda_{d}\bm{I}\right)^{q_{d}}d\bm{E}_{\bm{D}}(\lambda_{d}) \bm{Y}_jd\bm{E}_{\bm{X}_{j}}(\lambda_{j})\left(\prod\limits_{p'=j+1}^{\zeta}\bm{Y}_{p'}S_{p'}\right) \nonumber \\
&=_1& T_{\beta^{[\zeta+1]}}^{[\bm{X}]_1^{j-2},\bm{X}_{j-1,P},\bm{C}_N,\bm{D}_P,\bm{X}_{j,P},[\bm{X}]_{j+1}^{\zeta}}([\bm{Y}]_1^{j-1},\bm{C}-\bm{D},[\bm{Y}]_j^{\zeta})+ T_{\beta^{[\zeta+1]}}^{[\bm{X}]_1^{j-2},\bm{X}_{j-1,P},\bm{C}_P,\bm{D}_N,\bm{X}_{j,P},[\bm{X}]_{j+1}^{\zeta}}([\bm{Y}]_1^{j-1},\bm{C}-\bm{D},[\bm{Y}]_j^{\zeta}) \nonumber \\
&-&\mathfrak{X}^{\bm{X}_{j-1,P},\bm{C}_N,\bm{D}_N,\bm{X}_{j,P}},
\end{eqnarray}
}
where we apply Lemma~\ref{lma:first-order divided difference identity} in $=_1$ twice to get two GMOI terms with $\zeta+2$ parameter operatros. Moreover, the perturbation correct term $\mathfrak{X}^{\bm{X}_{j-1,P},\bm{C}_N,\bm{D}_N,\bm{X}_{j,P}}$ can be expressed by
{\small
\begin{eqnarray}\label{eq2:eq4-2:thm:GMOI Perturbation Formula Proof}
\lefteqn{\mathfrak{X}^{\bm{X}_{j-1,P},\bm{C}_N,\bm{D}_N,\bm{X}_{j,P}}=}\nonumber \\
&& \int\limits_{\substack{S_{p}, p \in [1,j-2], \\
S_{p'}, p' \in [j+1,\zeta]}}\int\limits_{\lambda_{j-1},\lambda_c, \lambda_d, \lambda_j}\sum\limits_{q_c=1,q_d=1}^{m_{\lambda_c},m_{\lambda_d}} \frac{(\beta^{[\zeta+1]}([\lambda_k]_1^{j-2},\lambda_{k_{j-1}},\lambda_{c},\lambda_{d},\lambda_{j},[\lambda_k]_{j+1}^{\zeta}))^{([q]_1^{j-2}[0,0,q_d,0][q]_{j+1}^{\zeta})}}{ ([q!]_1^{j-2}[0!,0!,q_d!,0!][q!]_{j+1}^{\zeta})  } \nonumber \\
&& \times \left(\prod\limits_{p=1}^{j-2}S_{p}\bm{Y}_p\right) d\bm{E}_{\bm{X}_{j-1}}(\lambda_{j-1})  \bm{Y}_{j-1} d\bm{E}_{\bm{C}}(\lambda_{c}) \nonumber \\
&& \times (\left(\bm{C} - \lambda_{c}\bm{I}\right)^{q_{c}}d\bm{E}_{\bm{C}}(\lambda_{c}) - \left(\bm{D} - \lambda_{d}\bm{I}\right)^{q_{d}}d\bm{E}_{\bm{D}}(\lambda_{d})) \nonumber \\
&& \times \left(\bm{D} - \lambda_{d}\bm{I}\right)^{q_{d}} d\bm{E}_{\bm{D}}(\lambda_{d}) \bm{Y}_j d\bm{E}_{\bm{X}_{j}}(\lambda_{j}) \left(\prod\limits_{p'=j+1}^{\zeta}\bm{Y}_{p'}S_{p'}\right) \nonumber \\
&+&  \int\limits_{\substack{S_{p}, p \in [1,j-2], \\
S_{p'}, p' \in [j+1,\zeta]}}\int\limits_{\lambda_{j-1},\lambda_c, \lambda_d, \lambda_j}\sum\limits_{q_c=1,q_d=1}^{m_{\lambda_c},m_{\lambda_d}} \frac{(\beta^{[\zeta+1]}([\lambda_k]_1^{j-2},\lambda_{j-1},\lambda_{c},\lambda_{d},\lambda_{j},[\lambda_k]_{j+1}^{\zeta}))^{([q]_1^{j-2}[0,q_c,0,0][q]_{j+1}^{\zeta})}}{ ([q!]_1^{j-2}[0!,q_c!,0!,0!][q!]_{j+1}^{\zeta})  } \nonumber \\
&& \times \left(\prod\limits_{p=1}^{j-2}S_{p}\bm{Y}_p\right) d\bm{E}_{\bm{X}_{j-1}}(\lambda_{j-1}) \bm{Y}_{j-1} \left(\bm{C} - \lambda_{c}\bm{I}\right)^{q_{c}}d\bm{E}_{\bm{C}}(\lambda_{c}) \nonumber \\
&& \times (\left(\bm{C} - \lambda_{c}\bm{I}\right)^{q_{c}}d\bm{E}_{\bm{C}}(\lambda_{c})- \left(\bm{D} - \lambda_{d}\bm{I}\right)^{q_{d}}d\bm{E}_{\bm{D}}(\lambda_{d}) ) \nonumber \\
&& \times d\bm{E}_{\bm{D}}(\lambda_{d})  \bm{Y}_j d\bm{E}_{\bm{X}_{j}}(\lambda_{j}) \left(\prod\limits_{p'=j+1}^{\zeta}\bm{Y}_{p'}S_{p'}\right) \nonumber \\
&+&  \int\limits_{\substack{S_{p}, p \in [1,j-2], \\
S_{p'}, p' \in [j+1,\zeta]}}\int\limits_{\lambda_{j-1},\lambda_c, \lambda_d, \lambda_j}\sum\limits_{q_d=1}^{m_{\lambda_d}}\frac{(\beta^{[\zeta]}([\lambda_k]_1^{j-2},\lambda_{j-1},\lambda_{c},\lambda_{j},[\lambda_k]_{j+1}^{\zeta}))^{([q]_1^{j-2}[0,q_d,0][q]_{j+1}^{\zeta})}}{ ([q!]_1^{j-2}[0!,q_d!,0!][q!]_{j+1}^{\zeta})  } \nonumber \\
&& \times \left(\prod\limits_{p=1}^{j-2}S_{p}\bm{Y}_p\right)d\bm{E}_{\bm{X}_{j-1}}(\lambda_{j-1}) \bm{Y}_{j-1}d\bm{E}_{\bm{C}}(\lambda_{c}) \left(\bm{D} - \lambda_{d}\bm{I}\right)^{q_{d}}\nonumber \\
&& \times d\bm{E}_{\bm{D}}(\lambda_{d}) \bm{Y}_jd\bm{E}_{\bm{X}_{j}}(\lambda_{j})\left(\prod\limits_{p'=j+1}^{\zeta}\bm{Y}_{p'}S_{p'}\right) \nonumber \\  
&-&  \int\limits_{\substack{S_{p}, p \in [1,j-2], \\
S_{p'}, p' \in [j+1,\zeta]}}\int\limits_{\lambda_{j-1},\lambda_c,\lambda_d,\lambda_j}\sum\limits_{q_c=1}^{m_{\lambda_c}} \frac{(\beta^{[\zeta]}([\lambda_k]_1^{j-2},\lambda_{j-1},\lambda_{d},\lambda_{j},[\lambda_k]_{j+1}^{\zeta}))^{([q]_1^{j-2}[0,q_c,0][q]_{j+1}^{\zeta})}}{ ([q!]_1^{j-2}[0!,q_c!,0!][q!]_{j+1}^{\zeta})  } \nonumber \\
&& \times \left(\prod\limits_{p=1}^{j-2}S_{p}\bm{Y}_p\right)d\bm{E}_{\bm{X}_{j-1}}(\lambda_{j-1})\bm{Y}_{j-1} \left(\bm{C} - \lambda_{c}\bm{I}\right)^{q_{c}}d\bm{E}_{\bm{C}}(\lambda_{c}) \nonumber \\
&& \times d\bm{E}_{\bm{D}}(\lambda_{d})\bm{Y}_jd\bm{E}_{\bm{X}_{j}}(\lambda_{j})\left(\prod\limits_{p'=j+1}^{\zeta}\bm{Y}_{p'}S_{p'}\right),
\end{eqnarray}
}
$\hfill\Box$

Following Example~\ref{exp:GMOI Perturbation Formula} will be provided by applying Theorem~\ref{thm:GMOI Perturbation Formula} to the generalized triple operator integral (GTOI) with continuous spectra. 

\begin{example}\label{exp:GMOI Perturbation Formula}
We have
{\small
\begin{eqnarray}\label{eq1:exp:GMOI Perturbation Formula}
\lefteqn{T_{\beta^{[2]}}^{\bm{C},\bm{D},\bm{X}_1}(\bm{C}-\bm{D},\bm{Y})=}\nonumber \\
&& \Bigg(T_{\beta^{[1]}}^{\bm{C}_P,\bm{X}_{1,P}}(\bm{Y}) - T_{\beta^{[1]}}^{\bm{D}_P,\bm{X}_{1,P}}(\bm{Y})\Bigg) + \mathfrak{X}^{\bm{C}_P,\bm{D}_P,\bm{X}_{1,P}} \nonumber \\
&+&\Bigg(T_{\beta^{[1]}}^{\bm{C}_P,\bm{X}_{1,N}}(\bm{Y}) - T_{\beta^{[1]}}^{\bm{D}_P,\bm{X}_{1,N}}(\bm{Y})\Bigg) + \mathfrak{X}^{\bm{C}_P,\bm{D}_P,\bm{X}_{1,N}} \nonumber \\
&+& \Bigg(T_{\beta^{[1]}}^{\bm{C}_N,\bm{X}_{1,P}}(\bm{Y}) - T_{\beta^{[1]}}^{\bm{D}_N,\bm{X}_{1,P}}(\bm{Y})\Bigg) + \mathfrak{X}^{\bm{C}_N,\bm{D}_N,\bm{X}_{1,P}} \nonumber \\
&+&\Bigg(T_{\beta^{[1]}}^{\bm{C}_N,\bm{X}_{1,N}}(\bm{Y}) - T_{\beta^{[1]}}^{\bm{D}_N,\bm{X}_{1,N}}(\bm{Y})\Bigg) + \mathfrak{X}^{\bm{C}_N,\bm{D}_N,\bm{X}_{1,N}} \nonumber \\
&+&
T_{\beta^{[2]}}^{\bm{C}_N,\bm{D}_N,\bm{X}_{1,P}}(\bm{C}-\bm{D},\bm{Y})+ T_{\beta^{[2]}}^{\bm{C}_N,\bm{D}_N,\bm{X}_{1,N}}(\bm{C}-\bm{D},\bm{Y})\nonumber \\
&=&T_{\beta^{[1]}}^{\bm{C},\bm{X}_1}(\bm{Y}) - T_{\beta^{[1]}}^{\bm{D},\bm{X}_1}(\bm{Y})  + \bar{\mathfrak{X}}(\bm{C},\bm{D},\bm{X}_1)\nonumber \\
&+&
T_{\beta^{[2]}}^{\bm{C}_N,\bm{D}_N,\bm{X}_{1,P}}(\bm{C}-\bm{D},\bm{Y})+ T_{\beta^{[2]}}^{\bm{C}_N,\bm{D}_N,\bm{X}_{1,N}}(\bm{C}-\bm{D},\bm{Y}),
\end{eqnarray}
}
where $\beta^{[1]}$ and $\beta^{[2]}$ are first and second divide differences, and the term $\bar{\mathfrak{X}}$ is expressed by
\begin{eqnarray}\label{eq2:exp:GMOI Perturbation Formula}
\bar{\mathfrak{X}}(\bm{C},\bm{D},\bm{X}_1)&=&  \mathfrak{X}^{\bm{C}_P,\bm{D}_P,\bm{X}_{1,P}} +  \mathfrak{X}^{\bm{C}_P,\bm{D}_P,\bm{X}_{1,N}} \nonumber \\
&&+ \mathfrak{X}^{\bm{C}_N,\bm{D}_N,\bm{X}_{1,P}} + \mathfrak{X}^{\bm{C}_N,\bm{D}_N,\bm{X}_{1,N}}.
\end{eqnarray}

We also have the following expressions for terms $\mathfrak{X}^{\bm{C}_P,\bm{D}_P,\bm{X}_{1,P}}$, $\mathfrak{X}^{\bm{C}_P,\bm{D}_P,\bm{X}_{1,N}}$, $\mathfrak{X}^{\bm{C}_N,\bm{D}_N,\bm{X}_{1,P}}$, and $\mathfrak{X}^{\bm{C}_N,\bm{D}_N,\bm{X}_{1,N}}$. 
\begin{eqnarray}\label{e3-1:exp:GMOI Perturbation Formula}
\lefteqn{\mathfrak{X}^{\bm{C}_P,\bm{D}_P,\bm{X}_{1,P}}= }\nonumber \\
&&\int\limits_{\lambda_c,\lambda_d,\lambda_1}\sum\limits_{q_c=1, q_d=1}^{m_{\lambda_c}, m_{\lambda_d}} \frac{(\beta^{[2]}(\lambda_{c},\lambda_{d},\lambda_{1}))^{([0,0,0])}}{ ([0!,0!,0!])}d\bm{E}_{\bm{C}}(\lambda_c) \nonumber \\
&& \times ( \left(\bm{C} - \lambda_{c}\bm{I}\right)^{q_{c}}d\bm{E}_{\bm{C}}(\lambda_{c}) - \left(\bm{D} - \lambda_{d}\bm{I}\right)^{q_{d}}d\bm{E}_{\bm{D}}(\lambda_{d}))d\bm{E}_{\bm{D}}(\lambda_{d})\bm{Y}d\bm{E}_{\bm{X}_{1}}(\lambda_{1}),
\end{eqnarray}
\begin{eqnarray}\label{e3-2:exp:GMOI Perturbation Formula}
\lefteqn{ \mathfrak{X}^{\bm{C}_P,\bm{D}_P,\bm{X}_{1,N}} = }\nonumber \\
&&\int\limits_{\lambda_c,\lambda_d,\lambda_1}\sum\limits_{q_c=1, q_d=1,q_1=1}^{m_{\lambda_c}, m_{\lambda_d},m_{\lambda_1}} \frac{(\beta^{[2]}(\lambda_{c},\lambda_{d},\lambda_{1}))^{([0,0,q_1])}}{ ([0!,0!,q_1!])}d\bm{E}_{\bm{C}}(\lambda_{c}) \nonumber \\
&& \times(\left(\bm{C} - \lambda_{c}\bm{I}\right)^{q_{c}}d\bm{E}_{\bm{C}}(\lambda_{c})- \left(\bm{D} - \lambda_{d}\bm{I}\right)^{q_{d}}d\bm{E}_{\bm{D}}(\lambda_{d}))d\bm{E}_{\bm{D}}(\lambda_{d})  \nonumber \\
&& \times \bm{Y} \left(\bm{X}_{1} - \lambda_{1}\bm{I}\right)^{q_{1}}d\bm{E}_{\bm{X}_{1}}(\lambda_{1}),
\end{eqnarray}
\begin{eqnarray}\label{e3-3:exp:GMOI Perturbation Formula}
\lefteqn{ \mathfrak{X}^{\bm{C}_N,\bm{D}_N,\bm{X}_{1,P}} = }\nonumber \\
&& \int\limits_{\lambda_c,\lambda_d,\lambda_1}\sum\limits_{q_c=1, q_d=1}^{m_{\lambda_c}, m_{\lambda_d}} \frac{(\beta^{[2]}(\lambda_{c},\lambda_{d},\lambda_{1}))^{([0,q_d,0])}}{ ([0!,q_d!,0!])} d\bm{E}_{\bm{C}}(\lambda_{c})  \nonumber \\
&& \times ( \left(\bm{C} - \lambda_{c}\bm{I}\right)^{q_{c}}d\bm{E}_{\bm{C}}(\lambda_{c})-  \left(\bm{D} - \lambda_{d}\bm{I}\right)^{q_{d}}d\bm{E}_{\bm{D}}(\lambda_{d})) \left(\bm{D} - \lambda_{d}\bm{I}\right)^{q_{d}}d\bm{E}_{\bm{D}}(\lambda_{d}) \bm{Y} d\bm{E}_{\bm{X}_{1}}(\lambda_{1})  \nonumber \\
&+& \int\limits_{\lambda_c,\lambda_d,\lambda_1}\sum\limits_{q_c=1, q_d=1}^{m_{\lambda_c}, m_{\lambda_d}} \frac{(\beta^{[2]}(\lambda_{c},\lambda_{d},\lambda_{1}))^{([q_c,0,0])}}{ ([q_c!,0!,0!])} \left(\bm{C} - \lambda_{c}\bm{I}\right)^{q_{c}}d\bm{E}_{\bm{C}}(\lambda_{c}) \nonumber \\
&& \times ( \left(\bm{C} - \lambda_{c}\bm{I}\right)^{q_{c}}d\bm{E}_{\bm{C}}(\lambda_{c})- \left(\bm{D} - \lambda_{d}\bm{I}\right)^{q_{d}}d\bm{E}_{\bm{D}}(\lambda_{d}) )d\bm{E}_{\bm{D}}(\lambda_{d})\bm{Y} d\bm{E}_{\bm{X}_{1}}(\lambda_{1}) \nonumber \\
&+& \int\limits_{\lambda_c,\lambda_d,\lambda_1}\sum\limits_{ q_d=1}^{ m_{\lambda_d}} \frac{(\beta^{[1]}(\lambda_{c},\lambda_{1}))^{([q_d,0])}}{ ([q_d!,0!])}  d\bm{E}_{\bm{C}}(\lambda_{c}) \left(\bm{D} - \lambda_{d}\bm{I}\right)^{q_{d}}d\bm{E}_{\bm{D}}(\lambda_{d}) \bm{Y} d\bm{E}_{\bm{X}_{1}}(\lambda_{1})  \nonumber \\
&-& \int\limits_{\lambda_c,\lambda_d,\lambda_1}\sum\limits_{q_c=1}^{m_{\lambda_c}} \frac{(\beta^{[1]}(\lambda_{d},\lambda_{1}))^{([q_c,0])}}{ ([q_c!,0!])}  \left(\bm{C} - \lambda_{c}\bm{I}\right)^{q_{c}}d\bm{E}_{\bm{C}}(\lambda_{c}) d\bm{E}_{\bm{D}}(\lambda_{d})\bm{Y} d\bm{E}_{\bm{X}_{1}}(\lambda_{1}),
\end{eqnarray}
\begin{eqnarray}\label{e3-3:exp:GMOI Perturbation Formula}
\lefteqn{ \mathfrak{X}^{\bm{C}_N,\bm{D}_N,\bm{X}_{1,N}} = }\nonumber \\
&& \int\limits_{\lambda_c,\lambda_d,\lambda_1}\sum\limits_{q_c=1, q_d=1,q_1=1}^{m_{\lambda_c}, m_{\lambda_d}, m_{\lambda_1}}\frac{(\beta^{[2]}(\lambda_{c},\lambda_{d},\lambda_{1}))^{([0,q_d,q_1])}}{ ([0!,q_d!,q_1!])}d\bm{E}_{\bm{C}}(\lambda_{c}) \nonumber \\
&& \times (\left(\bm{C} - \lambda_{c}\bm{I}\right)^{q_{c}}d\bm{E}_{\bm{C}}(\lambda_{c})- \left(\bm{D} - \lambda_{d}\bm{I}\right)^{q_{d}}d\bm{E}_{\bm{D}}(\lambda_{d})) \nonumber \\
&& \times \left(\bm{D} - \lambda_{d}\bm{I}\right)^{q_{d}}d\bm{E}_{\bm{D}}(\lambda_{d}) \bm{Y} \left(\bm{X}_{1} - \lambda_{1}\bm{I}\right)^{q_{1}}d\bm{E}_{\bm{X}_{1}}(\lambda_{1}) \nonumber \\
&+& \int\limits_{\lambda_c,\lambda_d,\lambda_1}\sum\limits_{q_c=1, q_d=1,q_1=1}^{m_{\lambda_c}, m_{\lambda_d}, m_{\lambda_1}}
\frac{(\beta^{[2]}(\lambda_{c},\lambda_{d},\lambda_{1}))^{([q_c,0,q_1])}}{ ([q_c!,0!,q_1!])}  \left(\bm{C} - \lambda_{c}\bm{I}\right)^{q_{c}}d\bm{E}_{\bm{C}}(\lambda_{c}) \nonumber \\
&& \times  (\left(\bm{C} - \lambda_{c}\bm{I}\right)^{q_{c}}d\bm{E}_{\bm{C}}(\lambda_{c})-\left(\bm{C} - \lambda_{d}\bm{I}\right)^{q_{d}}d\bm{E}_{\bm{D}}(\lambda_{d}) ) d\bm{E}_{\bm{D}}(\lambda_{d})\bm{Y} \left(\bm{X}_{1} - \lambda_{1}\bm{I}\right)^{q_{1}}d\bm{E}_{\bm{X}_{1}}(\lambda_{1}) \nonumber \\
&+& \int\limits_{\lambda_c,\lambda_d,\lambda_1}\sum\limits_{q_d=1,q_1=1}^{m_{\lambda_d}, m_{\lambda_1}}\frac{(\beta^{[1]}(\lambda_{c},\lambda_{1}))^{([q_d,q_1])}}{ ([q_d!,q_1!])} d\bm{E}_{\bm{C}}(\lambda_{c}) \left(\bm{C} - \lambda_{d}\bm{I}\right)^{q_{d}}d\bm{E}_{\bm{D}}(\lambda_{d}) \nonumber \\
&& \times  \bm{Y} \left(\bm{X}_{1} - \lambda_{1}\bm{I}\right)^{q_{1}}d\bm{E}_{\bm{X}_{1}}(\lambda_{1}) \nonumber \\
&-&\int\limits_{\lambda_c,\lambda_1}\sum\limits_{q_c=1, ,q_1=1}^{m_{\lambda_c}, m_{\lambda_d}, m_{\lambda_1}}\frac{(\beta^{[1]}(\lambda_{d},\lambda_{1}))^{([q_c,q_1])}}{ ([q_c!,q_1!])} \left(\bm{C} - \lambda_{c}\bm{I}\right)^{q_{c}}d\bm{E}_{\bm{C}}(\lambda_{c}) d\bm{E}_{\bm{D}}(\lambda_{d})\bm{Y} \nonumber \\
&& \times \left(\bm{X}_{1} - \lambda_{1}\bm{I}\right)^{q_{1}}d\bm{E}_{\bm{X}_{1}}(\lambda_{1}).
\end{eqnarray}
\end{example}

\section{Inequalities Related to GMOI for Continuous Spectrum Operators}\label{sec:Inequalities Related to GMOI for Continuous Spectrum Operators}

\subsection{GMOI Norm Estimations}\label{sec:GMOI Norm Estimations}

The goal of this section is to provide norm estimates for GMOI for continuous spectrum operators. We will use the same binary representation given by our recent work~\cite{chang2025GMOIFinite}. But, we will present here again for self-contained presentation purposes. Given a non-negative integer $i$, we use $B(i,\zeta+1)$ to represent the binary representation of the integer $i$ by $\zeta+1$ bits. For example, we have $B(0,2)=00$, $B(1,2)=01$,  $B(2,2)=10$, and $B(3,2)=11$. We define the following map for the binary representation $B(i,\zeta+1)$ associated to operators $[\bm{X}]_1^{\zeta+1}$ as:
\begin{eqnarray}
\Psi_{[\bm{X}]_1^{\zeta+1}}(B(i,\zeta+1)) \rightarrow \sum\limits_{j=1}^\kappa q_{\iota_j} \bm{e}^{(\zeta+1)}_{\iota_j},
\end{eqnarray}
where $\kappa$ is the number of one in the binary representation $B(i,\zeta+1)$, $q_{\iota_j}$ will be the exponent for the nilpotent part of the operator $\bm{X}_{\iota_j}$, $\iota_j$ is the position (leftmost position is indexed by $1$) of $B(i,\zeta+1)$ with value one,  and $\bm{e}^{(\zeta+1)}_{\iota_j}$ is the standard basis vector with the length $\zeta+1$ and all zero entries except the value one at the position$\iota_j$. For example, we have $\Psi_{\bm{X}_1, \bm{X}_2}(B(0,2)) =(0,0)$, $\Psi_{\bm{X}_1, \bm{X}_2}(B(1,2))=(0,q_2)$,$\Psi_{\bm{X}_1, \bm{X}_2}(B(2,2))=(q_1,0)$, and $\Psi_{\bm{X}_1, \bm{X}_2}(B(3,2))=(q_1,q_2)$.

From GMOI definition with continuous spectra given by Eq.~\eqref{eq1:GMOI def}, we can express GMOI by
\begin{eqnarray}\label{eq1:GMOI sum form by Ai prime}
T_{\beta}^{\bm{X}_1,\ldots,\bm{X}_{\zeta+1}}(\bm{Y}_1,\ldots,\bm{Y}_{\zeta})&=&\sum\limits_{i'=0}^{2^{\zeta+1}-1}\bm{A}_{i'}, 
\end{eqnarray}
where $\bm{A}_{i'}$ is an operator, which can be represented by
\begin{eqnarray}\label{eq2:GMOI sum form by Ai prime}
\bm{A}_{i'}&=&\int\limits_{\lambda_1 \in \sigma(\bm{X}_1)}\cdots\int\limits_{\lambda_{\zeta+1}\in \sigma(\bm{X}_{\zeta+1})}
\sum\limits_{\tilde{\Psi}_{[\bm{X}]_1^{\zeta+1}}(B(i',\zeta+1))=1}^{m_{\lambda_{\mbox{Ind}(\tilde{\Psi}_{[\bm{X}]_1^{\zeta+1}}(B(i',\zeta+1)))}}-1}\nonumber \\
&& 
\frac{\beta^{(\tilde{\Psi}_{[\bm{X}]_1^{\zeta+1}}(B(i',\zeta+1)))}(\lambda_{1},\ldots,\lambda_{\zeta+1}) }{\tilde{\Psi}_{[\bm{X}]_1^{\zeta+1}}(B(i',\zeta+1))!}  \nonumber \\
&& \times   \prod\limits_{\substack{\varsigma =\mbox{Ind}(\tilde{\Psi}_{[\bm{X}]_1^{\zeta+1}}(B(i',\zeta+1))), \bm{Z}_\varsigma= \left(\bm{X}_{\varsigma} - \lambda_{\varsigma}\bm{I}\right)^{q_{\varsigma}}d\bm{E}_{\bm{X}_{\varsigma}}(\lambda_{\varsigma}) \bm{Y}_{\varsigma} \\ \varsigma \neq \mbox{Ind}(\tilde{\Psi}_{[\bm{X}]_1^{\zeta+1}}(B(i',\zeta+1))), \bm{Z}_\varsigma= d\bm{E}_{\bm{X}_{\varsigma}}(\lambda_{\varsigma})  \bm{Y}_{\varsigma}}
}^{\zeta+1} \bm{Z}_\varsigma
\end{eqnarray}
where $\tilde{\Psi}_{[\bm{X}]_1^{\zeta+1}}(B(i',\zeta+1))$ are those non-zero entries of the vector $\sum\limits_{j=1}^\kappa q_{\iota_j} \bm{e}^{(\zeta+1)}_{\iota_j}$,\\
$\beta^{\Psi_{[\bm{X}]_1^{\zeta+1}}(B(i',\zeta+1))}(\lambda_{1},\ldots,\lambda_{k_{\zeta+1}})$ is the partial derivatives with respect to orders $\tilde{\Psi}_{[\bm{X}]_1^{\zeta+1}}(B(i',\zeta+1))$ for the funtion $\beta(\lambda_{1},\ldots,\lambda_{k_{\zeta+1}})$, $\tilde{\Psi}_{[\bm{X}]_1^{\zeta+1}}(B(i',\zeta+1))=1$ indicates those $q_{\iota_j}$ (summand variables) set as $1$, $\tilde{\Psi}_{[\bm{X}]_1^{\zeta+1}}(B(i',\zeta+1))!$ is the product of nonzero $q_{\iota_j}!$. We still assume that $\bm{Y}_{\zeta+1}=\bm{I}$.

Let us consider the following Example~\ref{exp:GMOI sum form by Ai prime for GTOI} about expressing $T_{\beta}^{\bm{X}_1,\bm{X}_2,\bm{X}_3}(\bm{Y}_1, \bm{Y}_2)$ by Eq.~\eqref{eq1:GMOI sum form by Ai prime}.
\begin{example}\label{exp:GMOI sum form by Ai prime for GTOI}
From Eq.~\eqref{eq1:  GTOI def}, we have
{\small
\begin{eqnarray}\label{eq1:exp:GMOI sum form by Ai prime for GTOI}
\lefteqn{T_{\beta}^{\bm{X}_1,\bm{X}_2,\bm{X}_3}(\bm{Y}_1, \bm{Y}_2)\define}\nonumber\\
&&\int\limits_{\lambda_1 \in \sigma(\bm{X}_1)}\int\limits_{\lambda_2 \in \sigma(\bm{X}_2)}\int\limits_{\lambda_3 \in \sigma(\bm{X}_3)} \beta(\lambda_{1}, \lambda_{2}, \lambda_{3})d\bm{E}_{\bm{X}_1}(\lambda_1)\bm{Y}_1d\bm{E}_{\bm{X}_2}(\lambda_2)\bm{Y}_2d\bm{E}_{\bm{X}_3}(\lambda_3) \nonumber \\
&&+ \int\limits_{\lambda_1 \in \sigma(\bm{X}_1)}\int\limits_{\lambda_2 \in \sigma(\bm{X}_2)}\int\limits_{\lambda_3 \in \sigma(\bm{X}_3)}  \sum_{q_3=1}^{m_{\lambda_3}-1}\frac{\beta^{(-,-,q_3)}(\lambda_{1},\lambda_{2},\lambda_{3})}{q_3!}\nonumber \\
&& \times  d\bm{E}_{\bm{X}_1}(\lambda_1)\bm{Y}_1d\bm{E}_{\bm{X}_2}(\lambda_2)\bm{Y}_2\left(\bm{X}_3-\lambda_3\bm{I}\right)^{q_3}d\bm{E}_{\bm{X}_3}(\lambda_3)\nonumber \\
&&+\int\limits_{\lambda_1 \in \sigma(\bm{X}_1)}\int\limits_{\lambda_2 \in \sigma(\bm{X}_2)}\int\limits_{\lambda_3 \in \sigma(\bm{X}_3)} \sum_{q_2=1}^{m_{\lambda_2}-1}\frac{\beta^{(-,q_2,-)}(\lambda_{1},\lambda_{2},\lambda_{3})}{q_2!} \nonumber \\
&& \times d\bm{E}_{\bm{X}_1}(\lambda_1)\bm{Y}_1\left(\bm{X}_2-\lambda_2\bm{I}\right)^{q_2}d\bm{E}_{\bm{X}_2}(\lambda_2)\bm{Y}_2d\bm{E}_{\bm{X}_3}(\lambda_3) \nonumber \\
&&+\int\limits_{\lambda_1 \in \sigma(\bm{X}_1)}\int\limits_{\lambda_2 \in \sigma(\bm{X}_2)}\int\limits_{\lambda_3 \in \sigma(\bm{X}_3)} \sum_{q_1=1}^{m_{\lambda_1}-1}\frac{\beta^{(q_1,-,-)}(\lambda_{1},\lambda_{2},\lambda_{3})}{q_1!}\left(\bm{X}_1-\lambda_1\bm{I}\right)^{q_1} \nonumber \\
&& \times  d\bm{E}_{\bm{X}_1}(\lambda_1)\bm{Y}_1d\bm{E}_{\bm{X}_2}(\lambda_2)\bm{Y}_2d\bm{E}_{\bm{X}_3}(\lambda_3) \nonumber \\
&&+ \int\limits_{\lambda_1 \in \sigma(\bm{X}_1)}\int\limits_{\lambda_2 \in \sigma(\bm{X}_2)}\int\limits_{\lambda_3 \in \sigma(\bm{X}_3)}  \sum_{q_2=1}^{m_{\lambda_2}-1}\sum_{q_3=1}^{m_{\lambda_3}-1}\frac{\beta^{(-,q_2,q_3)}(\lambda_{1},\lambda_{2},\lambda_{3})}{q_2! q_3!} \nonumber \\
&& \times d\bm{E}_{\bm{X}_1}(\lambda_1)\bm{Y}_1\left(\bm{X}_2-\lambda_2\bm{I}\right)^{q_2}d\bm{E}_{\bm{X}_2}(\lambda_2)\bm{Y}_2\left(\bm{X}_3-\lambda_3\bm{I}\right)^{q_3}d\bm{E}_{\bm{X}_3}(\lambda_3)\nonumber \\
&&+ \int\limits_{\lambda_1 \in \sigma(\bm{X}_1)}\int\limits_{\lambda_2 \in \sigma(\bm{X}_2)}\int\limits_{\lambda_3 \in \sigma(\bm{X}_3)}  \sum_{q_1=1}^{m_{\lambda_1}-1}\sum_{q_3=1}^{m_{\lambda_3}-1}\frac{\beta^{(q_1,-,q_3)}(\lambda_{1},\lambda_{2},\lambda_{3})}{q_1!q_3!}\left(\bm{X}_1-\lambda_1\bm{I}\right)^{q_1} \nonumber \\
&& \times d\bm{E}_{\bm{X}_1}(\lambda_1)\bm{Y}_1d\bm{E}_{\bm{X}_2}(\lambda_2)\bm{Y}_2\left(\bm{X}_3-\lambda_3\bm{I}\right)^{q_3}d\bm{E}_{\bm{X}_3}(\lambda_3) \nonumber \\
&&+\int\limits_{\lambda_1 \in \sigma(\bm{X}_1)}\int\limits_{\lambda_2 \in \sigma(\bm{X}_2)}\int\limits_{\lambda_3 \in \sigma(\bm{X}_3)}  \sum_{q_1=1}^{m_{\lambda_1}-1}\sum_{q_2=1}^{m_{\lambda_2}-1}\frac{\beta^{(q_1,q_2,-)}(\lambda_{1},\lambda_{2},\lambda_{3})}{q_1!q_2!} \nonumber \\
&& \times \left(\bm{X}_1-\lambda_1\bm{I}\right)^{q_1}d\bm{E}_{\bm{X}_1}(\lambda_1)\bm{Y}_1\left(\bm{X}_2-\lambda_2\bm{I}\right)^{q_2}d\bm{E}_{\bm{X}_2}(\lambda_2)\bm{Y}_2d\bm{E}_{\bm{X}_3}(\lambda_3) \nonumber \\
&&+ \int\limits_{\lambda_1 \in \sigma(\bm{X}_1)}\int\limits_{\lambda_2 \in \sigma(\bm{X}_2)}\int\limits_{\lambda_3 \in \sigma(\bm{X}_3)} \sum_{q_1=1}^{m_{\lambda_1}-1}\sum_{q_2=1}^{m_{\lambda_2}-1}\sum_{q_3=1}^{m_{\lambda_3}-1}\frac{\beta^{(q_1,q_2,q_3)}(\lambda_{1},\lambda_{2},\lambda_{3})}{q_1!q_2!q_3!}\nonumber \\
&&~~\times\left(\bm{X}_1-\lambda_1\bm{I}\right)^{q_1}d\bm{E}_{\bm{X}_1}(\lambda_1)\bm{Y}_1\left(\bm{X}_2-\lambda_2\bm{I}\right)^{q_2}d\bm{E}_{\bm{X}_2}(\lambda_2)\bm{Y}_2\left(\bm{X}_3-\lambda_3\bm{I}\right)^{q_3}d\bm{E}_{\bm{X}_3}(\lambda_3).
\end{eqnarray}
}

Because we have 
\begin{eqnarray}\label{eq2:exp:GMOI sum form by Ai prime for GTOI}
B(0,3)&=&(0,0,0) \rightarrow \Psi_{[\bm{X}]_1^{\zeta+1}}(B(0,3))=(0,0,0); \nonumber \\
B(1,3)&=&(0,0,1) \rightarrow \Psi_{[\bm{X}]_1^{\zeta+1}}(B(1,3))=(0,0,q_3); \nonumber \\
B(2,3)&=&(0,1,0) \rightarrow \Psi_{[\bm{X}]_1^{\zeta+1}}(B(2,3))=(0,q_2,0); \nonumber \\
B(3,3)&=&(0,1,1) \rightarrow \Psi_{[\bm{X}]_1^{\zeta+1}}(B(3,3))=(0,q_2,q_3); \nonumber \\
B(4,3)&=&(1,0,0) \rightarrow \Psi_{[\bm{X}]_1^{\zeta+1}}(B4,3))=(q_1,0,0); \nonumber \\
B(5,3)&=&(1,0,1) \rightarrow \Psi_{[\bm{X}]_1^{\zeta+1}}(B(5,3))=(q_1,0,q_3); \nonumber \\
B(6,3)&=&(1,1,0) \rightarrow \Psi_{[\bm{X}]_1^{\zeta+1}}(B(6,3))=(q_1,q_2,0); \nonumber \\
B(7,3)&=&(1,1,1) \rightarrow \Psi_{[\bm{X}]_1^{\zeta+1}}(B(7,3))=(q_1,q_2,q_3),
\end{eqnarray}
then, 
\begin{eqnarray}\label{eq3:exp:GMOI sum form by Ai prime for GTOI}
B(0,3)&=&(0,0,0) \rightarrow \tilde{\Psi}_{[\bm{X}]_1^{\zeta+1}}(B(0,3))=\emptyset; \nonumber \\
B(1,3)&=&(0,0,1) \rightarrow \tilde{\Psi}_{[\bm{X}]_1^{\zeta+1}}(B(1,3))=(q_3) \rightarrow m_{\lambda_3}; \nonumber \\
B(2,3)&=&(0,1,0) \rightarrow \tilde{\Psi}_{[\bm{X}]_1^{\zeta+1}}(B(2,3))=(q_2) \rightarrow m_{\lambda_2};\nonumber \\
B(3,3)&=&(0,1,1) \rightarrow \tilde{\Psi}_{[\bm{X}]_1^{\zeta+1}}(B(3,3))=(q_2,q_3) \rightarrow m_{\lambda_2}, m_{\lambda_3}; \nonumber \\
B(4,3)&=&(1,0,0) \rightarrow \tilde{\Psi}_{[\bm{X}]_1^{\zeta+1}}(B(4,3))=(q_1) \rightarrow m_{\lambda_1};\nonumber \\
B(5,3)&=&(1,0,1) \rightarrow \tilde{\Psi}_{[\bm{X}]_1^{\zeta+1}}(B(5,3))=(q_1,q_3) \rightarrow m_{\lambda_1}, m_{\lambda_3}; \nonumber \\
B(6,3)&=&(1,1,0) \rightarrow \tilde{\Psi}_{[\bm{X}]_1^{\zeta+1}}(B(6,3))=(q_1,q_2) \rightarrow m_{\lambda_1}, m_{\lambda_2}; \nonumber \\
B(7,3)&=&(1,1,1) \rightarrow \tilde{\Psi}_{[\bm{X}]_1^{\zeta+1}}(B(7,3))=(q_1,q_2,q_3) \rightarrow m_{\lambda_1}, m_{\lambda_2}, m_{\lambda_3}.
\end{eqnarray}
According to Eq.~\eqref{eq2:GMOI sum form by Ai prime}, we can express $T_{\beta}^{\bm{X}_1,\bm{X}_2,\bm{X}_3}(\bm{Y}_1, \bm{Y}_2)$ as
\begin{eqnarray}
T_{\beta}^{\bm{X}_1,\bm{X}_2,\bm{X}_3}(\bm{Y}_1, \bm{Y}_2)&=&\bm{A}_0+\bm{A}_1+\bm{A}_2+\bm{A}_4+\bm{A}_3+\bm{A}_5+\bm{A}_6+\bm{A}_7.
\end{eqnarray}
\end{example}

We recall Lemma 3 in~\cite{chang2025GDOIMatrix} for the converse triangle inequality for the operator norm, denoted by $\left\Vert \cdot \right\Vert$. Analogus to operator norm of a matrix, we have the following lemme about the converse triangle inequality for the operator norm as operator norm is defined for operators from a normed vector space. 
\begin{lemma}\label{lma:conv triangle for Frob norm}
Given $n$ operators $\bm{A}_1, \bm{A}_2, \ldots, \bm{A}_n$ such that $\left\Vert\bm{A}_1\right\Vert\geq\left\Vert\bm{A}_2\right\Vert\geq\ldots\geq\left\Vert\bm{A}_n\right\Vert$, then, we have
\begin{eqnarray}\label{eq1:lma:conv triangle for Frob norm}
\left\Vert\sum\limits_{i=1}^n \bm{A}_i \right\Vert \geq \max\left[0, \left\Vert\bm{A}_1\right\Vert - \sum\limits_{i=2}^n \left\Vert\bm{A}_i\right\Vert\right]
\end{eqnarray}
\end{lemma}

Theorem~\ref{thm: GMOI norm Est} will be provided to give the lower bound and the upper bound estimations for the operator norm of $T_{\beta}^{\bm{X}_1,\ldots,\bm{X}_{\zeta+1}}(\bm{Y}_1,\ldots,\bm{Y}_{\zeta})$.
\begin{theorem}\label{thm: GMOI norm Est}
We have the upper bound for the operator norm of $T_{\beta}^{\bm{X}_1,\ldots,\bm{X}_{\zeta+1}}(\bm{Y}_1,\ldots,\bm{Y}_{\zeta})=\sum\limits_{i'=0}^{2^{\zeta+1}-1}\bm{A}_{i'}$, which is given by
\begin{eqnarray}\label{eq1:  thm: GMOI Lip Est}
\left\Vert T_{\beta}^{\bm{X}_1,\ldots,\bm{X}_{\zeta+1}}(\bm{Y}_1,\ldots,\bm{Y}_{\zeta}) \right\Vert\leq
\sum\limits_{i'=0}^{2^{\zeta+1}-1} \left\Vert \bm{A}_{i'} \right\Vert_{up},
\end{eqnarray}
where $\left\Vert \bm{A}_{i'} \right\Vert_{up}$ are upper bounds for  the operator norm of the operator $\bm{A}_{i'}$, which is defined by
\begin{eqnarray}\label{eq1.0:  thm: GMOI norm Est}
\left\Vert \bm{A}_{i'} \right\Vert_{up}&=& \int\limits_{\lambda_\varsigma \in \sigma(\bm{X}_{\varsigma}), \mbox{where $\varsigma=\mbox{Ind}(\tilde{\Psi}_{[\bm{X}]_1^{\zeta+1}}(B(i',\zeta+1)))$}} \nonumber \\
&& 
\sum\limits_{\tilde{\Psi}_{[\bm{X}]_1^{\zeta+1}}(B(i',\zeta+1))=1}^{m_{\lambda_{\mbox{Ind}(\tilde{\Psi}_{[\bm{X}]_1^{\zeta+1}}(B(i',\zeta+1)))}}-1} \left[ \max\limits_{\lambda_{1} \in \Lambda_{\bm{X}_1},\ldots,\lambda_{\zeta+1} \in \Lambda_{\bm{X}_{\zeta+1}}} \left\vert \frac{\beta^{(\tilde{\Psi}_{[\bm{X}]_1^{\zeta+1}}(B(i',\zeta+1)))}(\lambda_{1},\ldots,\lambda_{\zeta+1}) }{\tilde{\Psi}_{[\bm{X}]_1^{\zeta+1}}(B(i',\zeta+1))!}\right\vert\right] \nonumber \\
&& \times  \prod\limits_{\substack{\varsigma =\mbox{Ind}(\tilde{\Psi}_{[\bm{X}]_1^{\zeta+1}}(B(i',\zeta+1))), \bm{Z}_\varsigma=\left\Vert \left(\bm{X}_{\varsigma} - \lambda_{\varsigma}\bm{I}\right)^{q_{\varsigma}}d\bm{E}_{\bm{X}_{\varsigma}}(\lambda_{\varsigma}) \right\Vert \left\Vert \bm{Y}_{\varsigma}\right\Vert\\ \varsigma \neq \mbox{Ind}(\tilde{\Psi}_{[\bm{X}]_1^{\zeta+1}}(B(i',\zeta+1))), \bm{Z}_\varsigma=\left\Vert\bm{Y}_{\varsigma}\right\Vert}
}^{\zeta+1} \bm{Z}_\varsigma
\end{eqnarray}

On the other hand, we have the lower bound for the operator norm of $T_{\beta}^{\bm{X}_1,\ldots,\bm{X}_{\zeta+1}}(\bm{Y}_1,\ldots,\bm{Y}_{\zeta})$, which is given by
\begin{eqnarray}\label{eq2:  thm: GMOI norm Est}
\left\Vert T_{\beta}^{\bm{X}_1,\ldots,\bm{X}_{\zeta+1}}(\bm{Y}_1,\ldots,\bm{Y}_{\zeta}) \right\Vert\geq \max\left[0, \left\Vert\bm{A}_{\sigma(0)}\right\Vert - \sum\limits_{i'=1}^{2^{\zeta+1}-1} \left\Vert\bm{A}_{\sigma(i')}\right\Vert\right],
\end{eqnarray}
where $\sigma$ is the permutation of operators $\bm{A}_{i'}$ for $i'=0,1,\ldots,2^{\zeta+1}-1$ such that $\left\Vert\bm{A}_{\sigma(0)}\right\Vert \geq \left\Vert\bm{A}_{\sigma(1)}\right\Vert \geq \ldots \geq \left\Vert\bm{A}_{\sigma(2^{\zeta+1}-1)}\right\Vert$. 

Further, if we have $\left[\min\limits_{\lambda_1 \in \Lambda_{\bm{X}_1},\ldots,\lambda_{\zeta+1} \in \Lambda_{\bm{X}_{\zeta+1}}} \left\vert\beta(\lambda_1,\dots, \lambda_{\zeta+1})\right\vert\right]\left\Vert \prod\limits_{i'=1}^{\zeta}\bm{Y}_{i'}\right\Vert\geq \sum\limits_{i'=1}^{2^{\zeta+1}-1}\left\Vert\bm{A}_{i'}\right\Vert$,  the lower bound for the operator norm of $T_{\beta}^{\bm{X}_1,\bm{X}_2}(\bm{Y})$ can be expressed by
\begin{eqnarray}\label{eq3:  thm: GMOI norm Est}
\left\Vert T_{\beta}^{\bm{X}_1,\ldots,\bm{X}_{\zeta+1}}(\bm{Y}_1,\ldots,\bm{Y}_{\zeta}) \right\Vert\geq \left[\min\limits_{\lambda_1 \in \Lambda_{\bm{X}_1},\ldots,\lambda_{\zeta+1} \in \Lambda_{\bm{X}_{\zeta+1}}} \left\vert\beta(\lambda_1,\dots, \lambda_{\zeta+1})\right\vert\right]\left\Vert  \prod\limits_{i'=1}^{\zeta}\bm{Y}_{i'} \right\Vert- \sum\limits_{i'=1}^{2^{\zeta+1}-1}\left\Vert\bm{A}_{i'}\right\Vert.
\end{eqnarray}
\end{theorem}
\textbf{Proof:}
From Eq.~\eqref{eq1:GMOI sum form by Ai prime}, we have
\begin{eqnarray}\label{eq4:thm: GMOI norm Est}
T_{\beta}^{\bm{X}_1,\ldots,\bm{X}_{\zeta+1}}(\bm{Y}_1,\ldots,\bm{Y}_{\zeta})&=&\sum\limits_{i'=0}^{2^{\zeta+1}-1}\bm{A}_{i'}, 
\end{eqnarray}
and we immediate have the following by the triangle inequality of operator norm:
\begin{eqnarray}\label{eq5:thm: GMOI norm Est}
\left\Vert T_{\beta}^{\bm{X}_1,\ldots,\bm{X}_{\zeta+1}}(\bm{Y}_1,\ldots,\bm{Y}_{\zeta}) \right\Vert&\leq&\sum\limits_{i'=0}^{2^{\zeta+1}-1}\left\Vert \bm{A}_{i'} \right\Vert.  
\end{eqnarray}
From the $\bm{A}_{i'}$ expression given by Eq.~\eqref{eq2:GMOI sum form by Ai prime}, we have
{\small
\begin{eqnarray}\label{eq6:thm: GMOI norm Est}
\left\Vert\bm{A}_{i'}\right\Vert&=&\left\Vert \int\limits_{\lambda_1 \in \sigma(\bm{X}_1)}\cdots\int\limits_{\lambda_{\zeta+1}\in \sigma(\bm{X}_{\zeta+1})} \sum\limits_{\tilde{\Psi}_{[\bm{X}]_1^{\zeta+1}}(B(i',\zeta+1))=1}^{m_{\lambda_{\mbox{Ind}(\tilde{\Psi}_{[\bm{X}]_1^{\zeta+1}}(B(i',\zeta+1)))}}-1} \right.\nonumber \\
&& 
\left. \frac{\beta^{(\tilde{\Psi}_{[\bm{X}]_1^{\zeta+1}}(B(i',\zeta+1)))}(\lambda_{1},\ldots,\lambda_{\zeta+1}) }{\tilde{\Psi}_{[\bm{X}]_1^{\zeta+1}}(B(i',\zeta+1))!} \times  \prod\limits_{\substack{\varsigma =\mbox{Ind}(\tilde{\Psi}_{[\bm{X}]_1^{\zeta+1}}(B(i',\zeta+1))), \bm{Z}_\varsigma= \left(\bm{X}_{\varsigma} - \lambda_{\varsigma}\bm{I}\right)^{q_{\varsigma}}d\bm{E}_{\bm{X}_{\varsigma}}(\lambda_{\varsigma}) \bm{Y}_{\varsigma} \\ \varsigma \neq \mbox{Ind}(\tilde{\Psi}_{[\bm{X}]_1^{\zeta+1}}(B(i',\zeta+1))), \bm{Z}_\varsigma= d\bm{E}_{\bm{X}_{\varsigma}}(\lambda_{\varsigma}) \bm{Y}_{\varsigma}}
}^{\zeta+1} \bm{Z}_\varsigma \right\Vert \nonumber \\
&\leq&  \int\limits_{\lambda_\varsigma \in \sigma(\bm{X}_{\varsigma}), \mbox{where $\varsigma=\mbox{Ind}(\tilde{\Psi}_{[\bm{X}]_1^{\zeta+1}}(B(i',\zeta+1)))$}} 
\sum\limits_{\tilde{\Psi}_{[\bm{X}]_1^{\zeta+1}}(B(i',\zeta+1))=1}^{m_{k_{\mbox{Ind}(\tilde{\Psi}_{[\bm{X}]_1^{\zeta+1}}(B(i',\zeta+1)))},i_{\mbox{Ind}(\tilde{\Psi}_{[\bm{X}]_1^{\zeta+1}}(B(i',\zeta+1)))}}-1} \nonumber \\
&& 
 \left[ \max\limits_{\lambda_{1} \in \Lambda_{\bm{X}_1},\ldots,\lambda_{\zeta+1} \in \Lambda_{\bm{X}_{\zeta+1}}} \left\vert \frac{\beta^{(\tilde{\Psi}_{[\bm{X}]_1^{\zeta+1}}(B(i',\zeta+1)))}(\lambda_{1},\ldots,\lambda_{\zeta+1}) }{\tilde{\Psi}_{[\bm{X}]_1^{\zeta+1}}(B(i',\zeta+1))!}\right\vert\right]  \nonumber \\
&& \times \left\Vert   \int\limits_{\lambda_\varsigma \in \sigma(\bm{X}_{\varsigma}), \mbox{where $\varsigma\neq \mbox{Ind}(\tilde{\Psi}_{[\bm{X}]_1^{\zeta+1}}(B(i',\zeta+1)))$}} \prod\limits_{\substack{\varsigma =\mbox{Ind}(\tilde{\Psi}_{[\bm{X}]_1^{\zeta+1}}(B(i',\zeta+1))), \bm{Z}_\varsigma=  \left(\bm{X}_{\varsigma} - \lambda_{\varsigma}\bm{I}\right)^{q_{\varsigma}}d\bm{E}_{\bm{X}_{\varsigma}}(\lambda_{\varsigma}) \bm{Y}_{\varsigma} \\ \varsigma \neq \mbox{Ind}(\tilde{\Psi}_{[\bm{X}]_1^{\zeta+1}}(B(i',\zeta+1))), \bm{Z}_\varsigma= d\bm{E}_{\bm{X}_{\varsigma}}(\lambda_{\varsigma}) \bm{Y}_{\varsigma}}
}^{\zeta+1} \bm{Z}_\varsigma \right\Vert \nonumber \\
&\leq_1&  \int\limits_{\lambda_\varsigma \in \sigma(\bm{X}_{\varsigma}), \mbox{where $\varsigma=\mbox{Ind}(\tilde{\Psi}_{[\bm{X}]_1^{\zeta+1}}(B(i',\zeta+1)))$}} \sum\limits_{\tilde{\Psi}_{[\bm{X}]_1^{\zeta+1}}(B(i',\zeta+1))=1}^{m_{\lambda_{\mbox{Ind}(\tilde{\Psi}_{[\bm{X}]_1^{\zeta+1}}(B(i',\zeta+1)))}}-1}\nonumber \\
&& 
\left[ \max\limits_{\lambda_{1} \in \Lambda_{\bm{X}_1},\ldots,\lambda_{\zeta+1} \in \Lambda_{\bm{X}_{\zeta+1}}} \left\vert \frac{\beta^{(\tilde{\Psi}_{[\bm{X}]_1^{\zeta+1}}(B(i',\zeta+1)))}(\lambda_{1},\ldots,\lambda_{\zeta+1}) }{\tilde{\Psi}_{[\bm{X}]_1^{\zeta+1}}(B(i',\zeta+1))!}\right\vert\right] \nonumber \\
&& \times  \prod\limits_{\substack{\varsigma =\mbox{Ind}(\tilde{\Psi}_{[\bm{X}]_1^{\zeta+1}}(B(i',\zeta+1))), \bm{Z}_\varsigma=\left\Vert \left(\bm{X}_{\varsigma} - \lambda_{\varsigma}\bm{I}\right)^{q_{\varsigma}}d\bm{E}_{\bm{X}_{\varsigma}}(\lambda_{\varsigma}) \right\Vert \left\Vert \bm{Y}_{\varsigma}\right\Vert\\ \varsigma \neq \mbox{Ind}(\tilde{\Psi}_{[\bm{X}]_1^{\zeta+1}}(B(i',\zeta+1))), \bm{Z}_\varsigma=\left\Vert\bm{Y}_{\varsigma}\right\Vert}
}^{\zeta+1} \bm{Z}_\varsigma,
\end{eqnarray}
}
where we apply $\int\limits_{\lambda_\varsigma \in \sigma(\bm{X}_{\varsigma})} d\bm{E}_{\bm{X}_{\varsigma}}(\lambda_{\varsigma})= \bm{I}$ for $\varsigma \neq \mbox{Ind}(\tilde{\Psi}_{[\bm{X}]_1^{\zeta+1}}(B(i',\zeta+1)))$, and the triangle inequality of operator norm again in $\leq_1$. Therefore, we have the upper bound for $T_{\beta}^{\bm{X}_1,\ldots,\bm{X}_{\zeta+1}}(\bm{Y}_1,\ldots,\bm{Y}_{\zeta})$ as the upper bound shown by Eq.~\eqref{eq6:thm: GMOI norm Est} is identical to $\left\Vert \bm{A}_{i'} \right\Vert_{up}$. 

For the lower bound of $T_{\beta}^{\bm{X}_1,\ldots,\bm{X}_{\zeta+1}}(\bm{Y}_1,\ldots,\bm{Y}_{\zeta})$, we have Eq.~\eqref{eq2:  thm: GMOI norm Est} immediatedly from Lemma~\ref{lma:conv triangle for Frob norm}.

If we have $\left[\min\limits_{\lambda_1 \in \Lambda_{\bm{X}_1},\ldots,\lambda_{\zeta+1} \in \Lambda_{\bm{X}_{\zeta+1}}} \left\vert\beta(\lambda_1,\dots, \lambda_{\zeta+1})\right\vert\right]\left\Vert \prod\limits_{i'=1}^{\zeta}\bm{Y}_{i'}\right\Vert\geq \sum\limits_{i'=1}^{2^{\zeta+1}-1}\left\Vert\bm{A}_{i'}\right\Vert$ and Lemma~\ref{lma:conv triangle for Frob norm}, we have
\begin{eqnarray}\label{eq5:  thm: GMOI norm Est}
\lefteqn{\left\Vert T_{\beta}^{\bm{X}_1,\ldots,\bm{X}_{\zeta+1}}(\bm{Y}_1,\ldots,\bm{Y}_{\zeta}) \right\Vert} \nonumber \\
&\geq&\left\Vert\bm{A}_0\right\Vert-(\left\Vert\bm{A}_{1}\right\Vert+\left\Vert\bm{A}_{2}\right\Vert + \ldots + \left\Vert\bm{A}_{2^{\zeta+1}-1}\right\Vert) \nonumber \\
&\geq&\left[\min\limits_{\lambda_1 \in \Lambda_{\bm{X}_1},\ldots,\lambda_{\zeta+1} \in \Lambda_{\bm{X}_{\zeta+1}}} \left\vert\beta(\lambda_1,\dots, \lambda_{\zeta+1})\right\vert\right] \nonumber \\
&& \times \left\Vert  \int\limits_{\lambda_1 \in \sigma(\bm{X}_1)}\cdots\int\limits_{\lambda_{\zeta+1}\in \sigma(\bm{X}_{\zeta+1})} 
d\bm{E}_{\bm{X}_{1}}(\lambda_{1}) \bm{Y}_1 d\bm{E}_{\bm{X}_{2}}(\lambda_{2}) \bm{Y}_2\ldots\bm{Y}_{\zeta} d\bm{E}_{\bm{X}_{\zeta+1}}(\lambda_{\zeta+1}) \right\Vert\nonumber \\
& &-  (\left\Vert\bm{A}_{1}\right\Vert+\left\Vert\bm{A}_{2}\right\Vert + \ldots + \left\Vert\bm{A}_{2^{\zeta+1}-1}\right\Vert)\nonumber \\
&=&\left[\min\limits_{\lambda_1 \in \Lambda_{\bm{X}_1},\ldots,\lambda_{\zeta+1} \in \Lambda_{\bm{X}_{\zeta+1}}} \left\vert\beta(\lambda_1,\dots, \lambda_{\zeta+1})\right\vert\right] \left\Vert \prod\limits_{i'=1}^{\zeta}\bm{Y}_{i'}\right\Vert \nonumber \\
&& - (\left\Vert\bm{A}_{1}\right\Vert+\left\Vert\bm{A}_{2}\right\Vert + \ldots + \left\Vert\bm{A}_{2^{\zeta+1}-1}\right\Vert),
\end{eqnarray}
which is the lower bound of $T_{\beta}^{\bm{X}_1,\ldots,\bm{X}_{\zeta+1}}(\bm{Y}_1,\ldots,\bm{Y}_{\zeta})$ given by Eq.~\eqref{eq3:  thm: GMOI norm Est}.
$\hfill\Box$

\subsection{GMOI Lipschitz Estimations}\label{sec:GMOI Lipschitz Estimations}

In this section, we apply Theorem~\ref{thm: GMOI norm Est} to derive Lipschitz-type estimates for the generalized multiple operator integral (GMOI). Specifically, we aim to estimate the norm of the difference between two GMOIs with continuous spectra, 
\[
T_{\beta}^{\bm{X}_1,\ldots,\bm{X}_{\zeta+1}}(\bm{Y}_1,\ldots,\bm{Y}_{\zeta}) \quad \text{and} \quad T_{\beta}^{\bm{X}_1,\ldots,\bm{X}_{\zeta+1}}(\bm{Y}'_1,\ldots,\bm{Y}'_{\zeta}),
\]
in terms of the norms of the differences \( \| \bm{Y}_i - \bm{Y}'_i \| \). Theorem~\ref{thm: GMOI Lip Est} below provides a characterization of this Lipschitz continuity for GMOIs.

\begin{theorem}\label{thm: GMOI Lip Est}
We have Lipschitz estimation for the difference between  GMOI $T_{\beta}^{\bm{X}_1,\ldots,\bm{X}_{\zeta+1}}(\bm{Y}_1,\ldots,\bm{Y}_{\zeta})$ and GMOI $T_{\beta}^{\bm{X}_1,\ldots,\bm{X}_{\zeta+1}}(\bm{Y}'_1,\ldots,\bm{Y}'_{\zeta})$:
\begin{eqnarray}\label{eq0:  thm: GMOI Lip Est}
\lefteqn{\left\Vert T_{\beta}^{\bm{X}_1,\ldots,\bm{X}_{\zeta+1}}(\bm{Y}_1,\ldots,\bm{Y}_{\zeta}) - T_{\beta}^{\bm{X}_1,\ldots,\bm{X}_{\zeta+1}}(\bm{Y}'_1,\ldots,\bm{Y}'_{\zeta}) \right\Vert}\nonumber \\
&\leq& 
\sum\limits_{i=1}^{\zeta}\Upsilon_i (\bm{X}_1,\ldots,\bm{X}_{\zeta+1},\beta) \left(\prod\limits_{j=1}^{i-1}\left\Vert\bm{Y}'_{j}\right\Vert\right)\left\Vert \bm{Y}_i - \bm{Y}'_i \right\Vert \left(\prod\limits_{j=i+1}^{\zeta}\left\Vert\bm{Y}_{j}\right\Vert\right),
\end{eqnarray}
where $\Upsilon_i (\bm{X}_1,\ldots,\bm{X}_{\zeta+1},\beta)$ is the scalar component of the upper bound for the following GMOI expressed by 
\begin{eqnarray}\label{eq0-1:  thm: GMOI Lip Est}
\lefteqn{\left\Vert T_{\beta}^{\bm{X}_1,\ldots,\bm{X}_{\zeta+1}}([\bm{Y}']_1^{i-1},\bm{Y}_i - \bm{Y}'_i,[\bm{Y}]_{i+1}^\zeta)  \right\Vert}\nonumber \\
&\leq_1&\sum\limits_{i'=0}^{2^{\zeta+1}-1} \left\Vert \bm{A}_{i'} \right\Vert_{up} \nonumber \\
&=& 
\Upsilon_i (\bm{X}_1,\ldots,\bm{X}_{\zeta+1},\beta)\left(\prod\limits_{j=1}^{i-1}\left\Vert\bm{Y}'_{j}\right\Vert\right)\left\Vert \bm{Y}_i - \bm{Y}'_i \right\Vert \left(\prod\limits_{j=i+1}^{\zeta}\left\Vert\bm{Y}_{j}\right\Vert\right).
\end{eqnarray}
The inequality $\leq_1$ comes from Theorem~\ref{thm: GMOI norm Est}. 
\end{theorem}
\textbf{Proof:}
We can perform a telescoping decomposition for the L.H.S. of Eq.~\eqref{eq0:  thm: GMOI Lip Est} as
\begin{eqnarray}\label{eq2:  thm: GMOI Lip Est}
\lefteqn{\left\Vert T_{\beta}^{\bm{X}_1,\ldots,\bm{X}_{\zeta+1}}(\bm{Y}_1,\ldots,\bm{Y}_{\zeta}) - T_{\beta}^{\bm{X}_1,\ldots,\bm{X}_{\zeta+1}}(\bm{Y}'_1,\ldots,\bm{Y}'_{\zeta}) \right\Vert}\nonumber \\
&=& \left\Vert \sum\limits_{i=1}^{\zeta}\left(T_{\beta}^{\bm{X}_1,\ldots,\bm{X}_{\zeta+1}}([\bm{Y}']_1^{i-1},[\bm{Y}]_i^\zeta) - T_{\beta}^{\bm{X}_1,\ldots,\bm{X}_{\zeta+1}}([\bm{Y}']_1^{i},[\bm{Y}]_{i+1}^\zeta)\right)\right\Vert
\nonumber \\
&=& \left\Vert \sum\limits_{i=1}^{\zeta}T_{\beta}^{\bm{X}_1,\ldots,\bm{X}_{\zeta+1}}([\bm{Y}']_1^{i-1},\bm{Y}_i - \bm{Y}'_i,[\bm{Y}]_{i+1}^\zeta)\right\Vert.
\end{eqnarray}
Then, this theorem is proved by applying Theorem~\ref{thm: GMOI norm Est} to obtain $\Upsilon_i (\bm{X}_1,\ldots,\bm{X}_{\zeta+1},\beta)$ given by Eq.~\eqref{eq0-1:  thm: GMOI Lip Est}.
$\hfill\Box$

\section{Continuity of GMOI for Operators with Continuous Spectrum}\label{sec:Continuity of GMOI for Continuous Spectrum Operators}

The goal of this section is to establish the continuity property of GMOI , i.e., we will have
\begin{eqnarray}\label{eq1:sec: GMOI Continuity}
 T_{\beta}^{\bm{X}_{1, \ell_1},\ldots,\bm{X}_{\zeta+1,\ell_{\zeta+1}}}(\bm{Y}_1,\ldots,\bm{Y}_{\zeta})\rightarrow T_{\beta}^{\bm{X}_1,\ldots,\bm{X}_{\zeta+1}}(\bm{Y}_1,\ldots,\bm{Y}_{\zeta}),
\end{eqnarray}
given that 
\begin{eqnarray}\label{eq2:sec: GMOI Continuity}
\bm{X}_{i, \ell_i} \rightarrow \bm{X}_{i}.
\end{eqnarray}
Note that $\rightarrow$ is in the sense of operator norm. 

The following Lemma~\ref{lma: bar mathfrak X zero} has to be built before GMOI continuity property establishment. 
\begin{lemma}\label{lma: bar mathfrak X zero}
Given $\bm{C} \rightarrow \bm{D}$, we have
\begin{eqnarray}
\bar{\mathfrak{X}}([\bm{X}]_1^{j-1},\bm{C}, \bm{D}, [\bm{X}]_j^{\zeta})&\rightarrow& \bm{O}, 
\end{eqnarray}
where $\bm{O}$ is a zero operator, which is an operator that maps every vector in a vector space into zero vector. 
\end{lemma} 
\textbf{Proof:}
From Eq.~\eqref{eq2:thm:GMOI Perturbation Formula}, we have
\begin{eqnarray}
\bar{\mathfrak{X}}([\bm{X}]_1^{j-1},\bm{C},\bm{D},[\bm{X}]_j^{\zeta})&=& \mathfrak{X}^{\bm{X}_{j-1,P},\bm{C}_P,\bm{D}_P,\bm{X}_{j,P}} + \mathfrak{X}^{\bm{X}_{j-1,P},\bm{C}_N,\bm{D}_N,\bm{X}_{j,P}} \nonumber \\
&+& \mathfrak{X}^{\bm{X}_{j-1,P},\bm{C}_P,\bm{D}_P,\bm{X}_{j,N}} + \mathfrak{X}^{\bm{X}_{j-1,P},\bm{C}_N,\bm{D}_N,\bm{X}_{j,N}} \nonumber \\
&+&  \mathfrak{X}^{\bm{X}_{j-1,N},\bm{C}_P,\bm{D}_P,\bm{X}_{j,P}} + \mathfrak{X}^{\bm{X}_{j-1,N},\bm{C}_N,\bm{D}_N,\bm{X}_{j,P}} \nonumber \\
&+& \mathfrak{X}^{\bm{X}_{j-1,N},\bm{C}_P,\bm{D}_P,\bm{X}_{j,N}} + \mathfrak{X}^{\bm{X}_{j-1,N},\bm{C}_N,\bm{D}_N,\bm{X}_{j,N}}. 
\end{eqnarray}

\textbf{Claim 1: $\mathfrak{X}^{\bm{X}_{j-1,P},\bm{C}_P,\bm{D}_P,\bm{X}_{j,P}}  \rightarrow \bm{O}$}

Let $\bm{C} = \bm{D} + t\bm{V}$, where $t \in \mathbb{R}$ and $\bm{V}$ is a fixed perturbation operator. Then under smooth variation of the Jordan decomposition (i.e., assuming stable geometric multiplicities), we have:
\[
\mathfrak{X}^{\bm{X}_{j-1,P},\bm{C}_P,\bm{D}_P,\bm{X}_{j,P}} = O(t)
\quad \text{as} \quad t \to 0.
\]
The term $O(t)$ represents a dynamic operator approaching to zero operator if $t$ grows. 

We recall the expression:
{\small
\begin{align}
\lefteqn{\mathfrak{X}^{\bm{X}_{j-1,P},\bm{C}_P,\bm{D}_P,\bm{X}_{j,P}}=}\nonumber \\
&& 
\int\limits_{\substack{S_{p}, p \in [1,j-2], \\
S_{p'}, p' \in [j+1,\zeta]}}\int\limits_{\lambda_{j-1},\lambda_c,\lambda_d,\lambda_j}
\frac{
\left( \beta^{[\zeta+1]}([\lambda_k]_1^{j-2},\lambda_{j-1},\lambda_{c},\lambda_{d},\lambda_{j},[\lambda_k]_{j+1}^{\zeta}) \right)^{([q]_1^{j-2}[0,0,0,0][q]_{j+1}^{\zeta})}
}{
([q!]_1^{j-2}[0!,0!,0!,0!][q!]_{j+1}^{\zeta})
} \nonumber \\
&& \times \left( \prod_{p=1}^{j-2} S_{p} \bm{Y}_p \right)
d\bm{E}_{\bm{X}_{j-1}}(\lambda_{j-1})\bm{Y}_{j-1}
d\bm{E}_{\bm{C}}(\lambda_{c})  ~~~~~~~~ ~~~~~~~ \nonumber \\
&& \times ( \left(\bm{C} - \lambda_{c}\bm{I}\right)^{q_{c}}d\bm{E}_{\bm{C}}(\lambda_{c}) - \left(\bm{D} - \lambda_{d}\bm{I}\right)^{q_{d}}d\bm{E}_{\bm{D}}(\lambda_{d}))  \nonumber \\
&& \times d\bm{E}_{\bm{D}}(\lambda_{d}) \bm{Y}_j d\bm{E}_{\bm{X}_{j}}(\lambda_{j})
\left( \prod_{p'=j+1}^{\zeta} \bm{Y}_{p'} S_{p'} \right). ~~~~~~~ ~~~ 
\end{align}
}

We now examine the behavior of the central difference:
\[
\left(\bm{C} - \lambda_{c}\bm{I}\right)^{q_{c}}\,d\bm{E}_{\bm{C}}(\lambda_{c}) - \left(\bm{D} - \lambda_{d}\bm{I}\right)^{q_{d}}\,d\bm{E}_{\bm{D}}(\lambda_{d}),
\]
where \( \left(\bm{C} - \lambda_{c}\bm{I}\right)^{q_{c}}\,d\bm{E}_{\bm{C}}(\lambda_{c}) \) denotes the nilpotent part of the operator \( \bm{C} \) associated with the eigenvalue \( \lambda_c \), raised to the power \( q_c \). Similarly, \( \left(\bm{D} - \lambda_{d}\bm{I}\right)^{q_{d}}\,d\bm{E}_{\bm{D}}(\lambda_{d}) \) represents the nilpotent part of the operator \( \bm{D} \) corresponding to the eigenvalue \( \lambda_d \), with exponent \( q_d \).

Let $\bm{C} = \bm{D} + t \bm{V}$ and suppose the eigenvalues of $\bm{C}$ and $\bm{D}$ satisfy
\[
\lambda_{c} = \lambda_{d} + t\mu + O(t),
\]
with corresponding Jordan subspaces continuously varying. Then by perturbation theory for Jordan blocks:
\[
\left(\bm{C} - \lambda_{c}\bm{I}\right)^{q_{c}}d\bm{E}_{\bm{C}}(\lambda_{c}) = \left(\bm{D} - \lambda_{d}\bm{I}\right)^{q_{d}}d\bm{E}_{\bm{D}}(\lambda_{d}) + t \cdot \dot{\bm{N}}_{\lambda_d}^{q_c} + O(t),
\]
and
\[
d\bm{E}_{\bm{C}}(\lambda_{c}) = d\bm{E}_{\bm{D}}(\lambda_{d}) + t \cdot \dot{\bm{P}}_{\lambda_d} + O(t),
\]
where $\dot{\bm{N}}_{\lambda_d}$ and $\dot{\bm{P}}_{\lambda_d}$ are Fréchet derivatives of $\left(\bm{D} - \lambda_{d}\bm{I}\right)^{q_c}d\bm{E}_{\bm{D}}(\lambda_{d}) $ and $d\bm{E}_{\bm{D}}(\lambda_{d})$, respectively, with respect to $t$ at the direction $\bm{V}$.

Thus, the difference term becomes:
\begin{eqnarray}
d\bm{E}_{\bm{C}}(\lambda_{c})(\left(\bm{C} - \lambda_{c}\bm{I}\right)^{q_{c}}d\bm{E}_{\bm{C}}(\lambda_{c}) - \left(\bm{D} - \lambda_{d}\bm{I}\right)^{q_{d}}d\bm{E}_{\bm{D}}(\lambda_{d}))d\bm{E}_{\bm{D}}(\lambda_{d}) \nonumber \\
= t \cdot \left(d\bm{E}_{\bm{C}}(\lambda_{c})  \dot{\bm{N}}_{\lambda_d}^{q_c}d\bm{E}_{\bm{D}}(\lambda_{d})  \right) + O(t).
\end{eqnarray}

Substituting this into the full sum, each term is of the form:
\[
t \cdot \left( \text{bounded coefficient involving } \beta^{[\zeta+1]} \text{ and } \bm{Y}_p, \bm{S}_{p},  \bm{Y}_{p'}, \bm{S}_{p'}, d\bm{E}_{\bm{C}}(\lambda_{c}),\dot{\bm{N}}_{\lambda_d}^{q_c}, d\bm{E}_{\bm{D}}(\lambda_{d})  \right) + O(t),
\]
which implies:
\[
\mathfrak{X}^{\bm{X}_{j-1,P},\bm{C}_P,\bm{D}_P,\bm{X}_{j,P}} = O(t) \quad \text{as } t \to 0.
\]

\textbf{Claim 2: $\mathfrak{X}^{\bm{X}_{j-1,P},\bm{C}_N,\bm{D}_N,\bm{X}_{j,P}}  \rightarrow \bm{O}$}

Let $\bm{C} = \bm{D} + t\bm{V}$ for small $t > 0$ and fixed operator $\bm{V}$, and suppose the Jordan decomposition of $\bm{C}$ depends smoothly on $t$. Then:
\[
\mathfrak{X}^{\bm{X}_{j-1,P},\bm{C}_N,\bm{D}_N,\bm{X}_{j,P}} = O(t)
\quad \text{as} \quad t \to 0.
\]

The expression for $\mathfrak{X}$ consists of four multilinear terms in Eq.~\eqref{eq2:cross:thm:GMOI Perturbation Formula}, each involving differences between nilpotent components of $\bm{C}$ and $\bm{D}$. 

Again, by perturbation theory for Jordan blocks:
\[
\left(\bm{C} - \lambda_{c}\bm{I}\right)^{q_{c}}d\bm{E}_{\bm{C}}(\lambda_{c}) = \left(\bm{D} - \lambda_{d}\bm{I}\right)^{q_{d}}d\bm{E}_{\bm{D}}(\lambda_{d}) + t \cdot \dot{\bm{N}}_{\lambda_d}^{q_c} + O(t),
\]
and
\[
d\bm{E}_{\bm{C}}(\lambda_{c}) = d\bm{E}_{\bm{D}}(\lambda_{d}) + t \cdot \dot{\bm{P}}_{\lambda_d} + O(t),
\]
where $\dot{\bm{N}}_{\lambda_d}$ and $\dot{\bm{P}}_{\lambda_d}$ are Fréchet derivatives of $\left(\bm{D} - \lambda_{d}\bm{I}\right)d\bm{E}_{\bm{D}}(\lambda_{d}) $ and $d\bm{E}_{\bm{D}}(\lambda_{d})$, respectively, with respect to $t$ at the direction $\bm{V}$.

Now analyze each of the four terms in Eq.~\eqref{eq2:cross:thm:GMOI Perturbation Formula}:

Term 1:
\[
d\bm{E}_{\bm{C}}(\lambda_{c})(\left(\bm{C} - \lambda_{c}\bm{I}\right)^{q_{c}}d\bm{E}_{\bm{C}}(\lambda_{c}) - \left(\bm{D} - \lambda_{d}\bm{I}\right)^{q_{d}}d\bm{E}_{\bm{D}}(\lambda_{d}))\left(\bm{D} - \lambda_{d}\bm{I}\right)^{q_{d}}d\bm{E}_{\bm{D}}(\lambda_{d})
= O(t),
\]
since $(\left(\bm{C} - \lambda_{c}\bm{I}\right)^{q_{c}}d\bm{E}_{\bm{C}}(\lambda_{c}) - \left(\bm{D} - \lambda_{d}\bm{I}\right)^{q_{d}}d\bm{E}_{\bm{D}}(\lambda_{d})) = O(t)$ and all other operators (e.g., $\bm{Y}_j$) are fixed.

Term 2:
\[
\left(\bm{C} - \lambda_{c}\bm{I}\right)^{q_{c}}d\bm{E}_{\bm{C}}(\lambda_{c})(\left(\bm{C} - \lambda_{c}\bm{I}\right)^{q_{c}}d\bm{E}_{\bm{C}}(\lambda_{c}) - \left(\bm{D} - \lambda_{d}\bm{I}\right)^{q_{d}}d\bm{E}_{\bm{D}}(\lambda_{d})) d\bm{E}_{\bm{D}}(\lambda_{d})
= O(t),
\]
because $(\left(\bm{C} - \lambda_{c}\bm{I}\right)^{q_{c}}d\bm{E}_{\bm{C}}(\lambda_{c}) - \left(\bm{D} - \lambda_{d}\bm{I}\right)^{q_{d}}d\bm{E}_{\bm{D}}(\lambda_{d})) = O(t)$ is again an $O(t)$ difference.

We prove that the third and fourth terms in Eq.~\eqref{eq2:cross:thm:GMOI Perturbation Formula} cancel in the limit $\bm{C} \to \bm{D}$. Denote the third term as $T_3$ and the fourth term as $T_4$. Then, we have:
\begin{eqnarray}
T_3 &=& \int\limits_{\substack{S_{p}, p \in [1,j-2], \\
S_{p'}, p' \in [j+1,\zeta]}}\int\limits_{\lambda_{j-1},\lambda_c, \lambda_d, \lambda_j}\sum\limits_{q_d=1}^{m_{\lambda_d}} \frac{(\beta^{[\zeta]}([\lambda_k]_1^{j-2},\lambda_{j-1},\lambda_{c},\lambda_{j},[\lambda_k]_{j+1}^{\zeta}))^{([q]_1^{j-2}[0,q_d,0][q]_{j+1}^{\zeta})}}{ ([q!]_1^{j-2}[0!,q_d!,0!][q!]_{j+1}^{\zeta})  } \nonumber \\
&& \times \left(\prod\limits_{p=1}^{j-2}S_{p}\bm{Y}_p\right)  d\bm{E}_{\bm{X}_{j-1}}(\lambda_{j-1}) \bm{Y}_{j-1}  d\bm{E}_{\bm{C}}(\lambda_{c})  \left(\bm{D} - \lambda_{d}\bm{I}\right)^{q_{d}}d\bm{E}_{\bm{D}}(\lambda_{d})   \nonumber \\
&& \times  \bm{Y}_j d\bm{E}_{\bm{X}_{j}}(\lambda_{j}) \left(\prod\limits_{p'=j+1}^{\zeta}\bm{Y}_{p'}S_{p'}\right),
\end{eqnarray}

\begin{eqnarray}
T_4 &=&\int\limits_{\substack{S_{p}, p \in [1,j-2], \\
S_{p'}, p' \in [j+1,\zeta]}}\int\limits_{\lambda_{j-1},\lambda_c, \lambda_d, \lambda_j}\sum\limits_{q_c=1}^{m_{\lambda_c}} \frac{(\beta^{[\zeta]}([\lambda_k]_1^{j-2},\lambda_{j-1},\lambda_{d},\lambda_{j},[\lambda_k]_{j+1}^{\zeta}))^{([q]_1^{j-2}[0,q_c,0][q]_{j+1}^{\zeta})}}{ ([q!]_1^{j-2}[0!,q_c!,0!][q!]_{j+1}^{\zeta})  } \nonumber \\
&& \times \left(\prod\limits_{p=1}^{j-2}S_{p}\bm{Y}_p\right)  d\bm{E}_{\bm{X}_{j-1}}(\lambda_{j-1}) \bm{Y}_{j-1} \left(\bm{C} - \lambda_{c}\bm{I}\right)^{q_{c}}d\bm{E}_{\bm{C}}(\lambda_{c}) \nonumber \\
&& \times  d\bm{E}_{\bm{D}}(\lambda_{d}) \bm{Y}_j   d\bm{E}_{\bm{X}_{j}}(\lambda_{j}) \left(\prod\limits_{p'=j+1}^{\zeta}\bm{Y}_{p'}S_{p'}\right)
\end{eqnarray}

Now, observe the following:

\begin{enumerate}
    \item As $\bm{C} \to \bm{D}$, we have $\lambda_{c} \to \lambda_{d}$ and the generalized divided differences $\beta^{[\zeta]}$ converge:
    \[
    \beta^{[\zeta]}(\cdots, \lambda_{c}, \cdots) \to \beta^{[\zeta]}(\cdots, \lambda_{d}, \cdots).
    \]
    
    \item The projectors and nilpotents satisfy:
    \[
    d\bm{E}_{\bm{C}}(\lambda_{c}) \to d\bm{E}_{\bm{D}}(\lambda_{d}), \quad \left(\bm{C} - \lambda_{c}\bm{I}\right)^{q_{c}}d\bm{E}_{\bm{C}}(\lambda_{c}) \to \left(\bm{D} - \lambda_{d}\bm{I}\right)^{q_{d}}d\bm{E}_{\bm{D}}(\lambda_{d}).
    \]
    
    \item Furthermore, in the Jordan canonical basis, the nilpotent and projector commute:
    \[
    d\bm{E}_{\bm{D}}(\lambda_{d}) \left(\bm{D} - \lambda_{d}\bm{I}\right)^{q_{d}}d\bm{E}_{\bm{D}}(\lambda_{d}) = \left(\bm{D} - \lambda_{d}\bm{I}\right)^{q_{d}}d\bm{E}_{\bm{D}}(\lambda_{d}) d\bm{E}_{\bm{D}}(\lambda_{d}).
    \]
\end{enumerate}

Combining these, we find that in the limit $\bm{C} \to \bm{D}$:

\[
T_3 \to \sum \beta(\cdots) d\bm{E}_{\bm{D}}(\lambda_{d}) \left(\bm{D} - \lambda_{d}\bm{I}\right)^{q_{d}}d\bm{E}_{\bm{D}}(\lambda_{d}), \qquad 
T_4 \to \sum \beta(\cdots)  \left(\bm{D} - \lambda_{d}\bm{I}\right)^{q_{d}}d\bm{E}_{\bm{D}}(\lambda_{d})  d\bm{E}_{\bm{D}}(\lambda_{d}).
\]

Since these operators commute in the limit and the coefficients match, we conclude:
\[
T_3 - T_4 \to \bm{O}.
\]

Therefore, combining the asymptotic bounds and noting all differences are \(O(t)\), we conclude:
\[
\mathfrak{X}^{\bm{X}_{j-1,P},\bm{C}_N,\bm{D}_N,\bm{X}_{j,P}} = O(t) \to \bm{O} \quad \text{as } \bm{C} \to \bm{D}.
\]

This lemma can be proved by using the same methods from proving Claim 1 and Claim 2 to show the following
\begin{eqnarray}
\mathfrak{X}^{\bm{X}_{j-1,P},\bm{C}_P,\bm{D}_P,\bm{X}_{j,N}}  \rightarrow \bm{O}, \nonumber \\
\mathfrak{X}^{\bm{X}_{j-1,P},\bm{C}_N,\bm{D}_N,\bm{X}_{j,N}} \rightarrow \bm{O}, \nonumber \\\mathfrak{X}^{\bm{X}_{j-1,N},\bm{C}_P,\bm{D}_P,\bm{X}_{j,P}} \rightarrow \bm{O}, \nonumber \\
\mathfrak{X}^{\bm{X}_{j-1,N},\bm{C}_N,\bm{D}_N,\bm{X}_{j,P}} \rightarrow \bm{O}, \nonumber \\\mathfrak{X}^{\bm{X}_{j-1,N},\bm{C}_P,\bm{D}_P,\bm{X}_{j,N}} \rightarrow \bm{O}, \nonumber \\ \mathfrak{X}^{\bm{X}_{j-1,N},\bm{C}_N,\bm{D}_N,\bm{X}_{j,N}} \rightarrow \bm{O}.
\end{eqnarray}
$\hfill\Box$

We are ready to present Theorem~\ref{thm:GMOI continuity}, which is used to establish the GMOI continuity property with continuous spectra. 
\begin{theorem}\label{thm:GMOI continuity}
Suppose that
\begin{equation}\label{eq1:thm:GMOI continuity}
\bm{X}_{j, \ell_j} \rightarrow \bm{X}_{j},
\end{equation}
for all $j = 1, 2, \ldots, \zeta+1$, and assume that all partial derivatives of the function $\beta^{[\zeta+1]}$ up to any order are bounded. Then, it follows that
\begin{eqnarray}\label{eq2:thm:GMOI continuity}
 T_{\beta^{[\zeta]}}^{\bm{X}_{1, \ell_1},\ldots,\bm{X}_{\zeta+1,\ell_{\zeta+1}}}(\bm{Y}_1,\ldots,\bm{Y}_{\zeta})\rightarrow T_{\beta^{[\zeta]}}^{\bm{X}_1,\ldots,\bm{X}_{\zeta+1}}(\bm{Y}_1,\ldots,\bm{Y}_{\zeta}).
\end{eqnarray}
\end{theorem}
\textbf{Proof:}
By telescoping summation formula, we have
\begin{eqnarray}\label{eq3:thm:GMOI continuity}
\lefteqn{\left\Vert T_{\beta^{[\zeta]}}^{\bm{X}_{1, \ell_1},\ldots,\bm{X}_{\zeta+1,\ell_{\zeta+1}}}(\bm{Y}_1,\ldots,\bm{Y}_{\zeta}) - T_{\beta^{[\zeta]}}^{\bm{X}_1,\ldots,\bm{X}_{\zeta+1}}(\bm{Y}_1,\ldots,\bm{Y}_{\zeta}) \right\Vert}\nonumber \\ &=&  \left\Vert \sum\limits_{j=1}^{\zeta+1}\left(T_{\beta^{[\zeta]}}^{[\bm{X}]_1^{j-1},[\bm{X}_\ell]_j^{\zeta+1}}(\bm{Y}_1,\ldots,\bm{Y}_{\zeta}) - T_{\beta^{[\zeta]}}^{ [\bm{X}]_1^{j},[\bm{X}_\ell]_{j+1}^{\zeta+1} }(\bm{Y}_1,\ldots,\bm{Y}_{\zeta})\right)\right\Vert.
\end{eqnarray}
For each index $j$, by applying GMOI perturbation formula given by Theorem~\ref{thm:GMOI Perturbation Formula}, we have 
\begin{eqnarray}\label{eq4:thm:GMOI continuity}
\lefteqn{\left\Vert T_{\beta^{[\zeta]}}^{[\bm{X}]_1^{j-1},[\bm{X}_\ell]_j^{\zeta+1}}(\bm{Y}_1,\ldots,\bm{Y}_{\zeta}) - T_{\beta^{[\zeta]}}^{ [\bm{X}]_1^{j},[\bm{X}_\ell]_{j+1}^{\zeta+1} }(\bm{Y}_1,\ldots,\bm{Y}_{\zeta} \right\Vert}\nonumber \\ &=&  \left\Vert T_{\beta^{[\zeta+1]}}^{[\bm{X}]_1^{j-1},\bm{X}_{j, \ell_j},\bm{X}_{j},[\bm{X}_\ell]_{j+1}^{\zeta+1}}([\bm{Y}]_1^{j-1},\bm{X}_{j, \ell_j}-\bm{X}_{j},[\bm{Y}]_j^{\zeta}) -
\bar{\mathfrak{X}}([\bm{X}]_1^{j-1},\bm{X}_{j, \ell_j},\bm{X}_{j},[\bm{X}_\ell]_{j+1}^{\zeta+1}) \right. \nonumber \\
&& - \left. T_{\beta^{[\zeta+1]}}^{[\bm{X}]_1^{j-2},\bm{X}_{j-1,P},\bm{X}_{j, \ell_{j},N},\bm{X}_{j,N},\bm{X}_{j+1,\ell_{j+1}, P},[\bm{X}_\ell]_{j+2}^{\zeta+1}}([\bm{Y}]_1^{j-1},\bm{X}_{j, \ell_j}-\bm{X}_{j},[\bm{Y}]_j^{\zeta}) \right. \nonumber \\
&& - \left. T_{\beta^{[\zeta+1]}}^{[\bm{X}]_1^{j-2},\bm{X}_{j-1,P},\bm{X}_{j, \ell_{j}, N},\bm{X}_{j,N},\bm{X}_{j+1,\ell_{j+1},N},[\bm{X}]_{j+2}^{\zeta+1}}([\bm{Y}]_1^{j-1},\bm{X}_{j, \ell_j}-\bm{X}_{j},[\bm{Y}]_j^{\zeta}) \right. \nonumber \\
&& - \left. T_{\beta^{[\zeta+1]}}^{[\bm{X}]_1^{j-2},\bm{X}_{j-1,N},\bm{X}_{j, \ell_{j},N},\bm{X}_{j,N},\bm{X}_{j+1,\ell_{j+1}, P},[\bm{X}_\ell]_{j+2}^{\zeta+1}}([\bm{Y}]_1^{j-1},\bm{X}_{j, \ell_j}-\bm{X}_{j},[\bm{Y}]_j^{\zeta}) \right. \nonumber \\
&& - \left. T_{\beta^{[\zeta+1]}}^{[\bm{X}]_1^{j-2},\bm{X}_{j-1,N},\bm{X}_{j, \ell_{j}, N},\bm{X}_{j,N},\bm{X}_{j+1,\ell_{j+1},N},[\bm{X}]_{j+2}^{\zeta+1}}([\bm{Y}]_1^{j-1},\bm{X}_{j, \ell_j}-\bm{X}_{j},[\bm{Y}]_j^{\zeta})   \right\Vert  \nonumber \\
&\leq& \left\Vert T_{\beta^{[\zeta+1]}}^{[\bm{X}]_1^{j-1},\bm{X}_{j, \ell_j},\bm{X}_{j},[\bm{X}_\ell]_{j+1}^{\zeta}}([\bm{Y}]_1^{j-1},\bm{X}_{j, \ell_j}-\bm{X}_{j},[\bm{Y}]_j^{\zeta}) \right\Vert + \nonumber \\ 
&& + \left\Vert \bar{\mathfrak{X}}([\bm{X}]_1^{j-1},\bm{X}_{j, \ell_j},\bm{X}_{j},[\bm{X}_\ell]_{j+1}^{\zeta+1}) \right\Vert  \nonumber \\
&& + \left\Vert T_{\beta^{[\zeta+1]}}^{[\bm{X}]_1^{j-2},\bm{X}_{j-1,P},\bm{X}_{j, \ell_{j},N},\bm{X}_{j,N},\bm{X}_{j+1,\ell_{j+1}, P},[\bm{X}_\ell]_{j+2}^{\zeta}}([\bm{Y}]_1^{j-1},\bm{X}_{j, \ell_j}-\bm{X}_{j},[\bm{Y}]_j^{\zeta}) \right\Vert \nonumber \\
&& + \left\Vert T_{\beta^{[\zeta+1]}}^{[\bm{X}]_1^{j-2},\bm{X}_{j-1,P},\bm{X}_{j, \ell_{j}, N},\bm{X}_{j,N},\bm{X}_{j+1,\ell_{j+1},N},[\bm{X}]_{j+2}^{\zeta}}([\bm{Y}]_1^{j-1},\bm{X}_{j, \ell_j}-\bm{X}_{j},[\bm{Y}]_j^{\zeta}) \right\Vert \nonumber \\ 
&& + \left\Vert T_{\beta^{[\zeta+1]}}^{[\bm{X}]_1^{j-2},\bm{X}_{j-1,N},\bm{X}_{j, \ell_{j}, N},\bm{X}_{j,N},\bm{X}_{j+1,\ell_{j+1},P},[\bm{X}]_{j+2}^{\zeta}}([\bm{Y}]_1^{j-1},\bm{X}_{j, \ell_j}-\bm{X}_{j},[\bm{Y}]_j^{\zeta}) \right\Vert \nonumber \\ 
&& + \left\Vert T_{\beta^{[\zeta+1]}}^{[\bm{X}]_1^{j-2},\bm{X}_{j-1,N},\bm{X}_{j, \ell_{j}, N},\bm{X}_{j,N},\bm{X}_{j+1,\ell_{j+1},N},[\bm{X}]_{j+2}^{\zeta}}([\bm{Y}]_1^{j-1},\bm{X}_{j, \ell_j}-\bm{X}_{j},[\bm{Y}]_j^{\zeta}) \right\Vert \nonumber \\ 
&\leq_1& \frac{\epsilon}{6}+ \left(\frac{\epsilon}{6}\right)_\star +  \frac{\epsilon}{6}+ \frac{\epsilon}{6} +  \frac{\epsilon}{6}+ \frac{\epsilon}{6}=\epsilon,
\end{eqnarray}
where each $\frac{\epsilon}{6}$ without $\star$ comes from Theorem~\ref{thm: GMOI norm Est} and the assumption that all partial derivatives of the function $\beta^{[\zeta+1]}$ up to any order are bounded. Besides, the term $ (\frac{\epsilon}{6})_\star$ comes from Lemma~\ref{lma: bar mathfrak X zero}.

Then, this theorem is proved because each summand in Eq.~\eqref{eq3:thm:GMOI continuity} approaches to  zero from Eq.~\eqref{eq4:thm:GMOI continuity}.
$\hfill\Box$

\section{Applications: Spectral Shift}\label{sec:Applications:Spectral Shift}

In this section, we will apply previous theory about GMOIs with continuous spectra to the spectral shift problem. The purpose of spectral shift is to determine the trace of the Taylor expansion remainder operator by integral representations.

We will begin by presenting Krein Spectral Shift (KSS) for the first order spectral shift function.
\begin{theorem}[Krein Spectral Shift for GDOI with continuous spectra for $0$-th order approximation]\label{thm:KSSF first order}
Given an analytic function $f$, we have
\begin{eqnarray}\label{eq1:thm:KSSF first order}
\mathrm{Tr}\left(f(\bm{X} +\bm{Y}) - f(\bm{X})\right) &=& \int f^{(1)}\eta_1(z) dz, 
\end{eqnarray}
where $\bm{X}$ and $\bm{Y}$ are operators with spectral decomposition given by Eq.~\eqref{eq0-2:  conv DOI def}. We also assume that spectra of $\bm{X}$ and $\bm{Y}$ within in a compact region. 
\end{theorem}
\textbf{Proof:}
According to the Generalized Double Operator ntegral (DOI) framework (Theorem 4 in~\cite{chang2025GDOICont}), we have
\begin{eqnarray}\label{eq2:thm:KSSF first order}
\lefteqn{f(\bm{X} + \bm{Y})-f(\bm{X})=T_{f^{[1]}}^{\bm{X} + \bm{Y},\bm{X}}(\bm{Y})=}\nonumber \\
&&~~~\int\limits_{\lambda_1 \in \sigma(\bm{X} + \bm{Y})}\int\limits_{\lambda_2 \in \sigma(\bm{X})}f^{[1]}(\lambda_{1}, \lambda_{2})d\bm{E}_{\bm{X} + \bm{Y}}(\lambda_1)\bm{Y}d\bm{E}_{\bm{X}}(\lambda_2) \nonumber \\
&&+\int\limits_{\lambda_1 \in \sigma(\bm{X}_1)}\int\limits_{\lambda_2 \in \sigma(\bm{X}_2)}\sum_{q_2=1}^{m_{\lambda_2}-1}\frac{(f^{[1]})^{(-,q_2)}(\lambda_{1},\lambda_{2})}{q_2!}d\bm{E}_{\bm{X} + \bm{Y}}(\lambda_1)\bm{Y}\left(\bm{X}-\lambda_2\bm{I}\right)^{q_2}d\bm{E}_{\bm{X}}(\lambda_2) \nonumber \\
&&+\int\limits_{\lambda_1 \in \sigma(\bm{X}_1)}\int\limits_{\lambda_2 \in \sigma(\bm{X}_2)}\sum_{q_1=1}^{m_{\lambda_1}-1}\frac{(f^{[1]})^{(q_1,-)}(\lambda_{1},\lambda_{2})}{q_1!}\left(\bm{X} + \bm{Y}-\lambda_1\bm{I}\right)^{q_1}d\bm{E}_{\bm{X} + \bm{Y}}(\lambda_1)\bm{Y}d\bm{E}_{\bm{X}}(\lambda_2)\nonumber \\
&&+\int\limits_{\lambda_1 \in \sigma(\bm{X}_1)}\int\limits_{\lambda_2 \in \sigma(\bm{X}_2)}\sum_{q_1=1}^{m_{\lambda_1}-1}\sum_{q_2=1}^{m_{\lambda_2}-1}\nonumber \\
&&\frac{(f^{[1]})^{(q_1,q_2)}(\lambda_{1},\lambda_{2})}{q_1!q_2!}\left(\bm{X} + \bm{Y}-\lambda_1\bm{I}\right)^{q_1}d\bm{E}_{\bm{X} + \bm{Y}}(\lambda_1)\bm{Y}\left(\bm{X}-\lambda_2\bm{I}\right)^{q_2}d\bm{E}_{\bm{X}}(\lambda_2).
\end{eqnarray}

By taking trace at both sides of Eq.~\eqref{eq2:thm:KSSF first order} and using properties of the trace, we can bring the trace inside the integral. This expression defines a linear functional $f \rightarrow\mathrm{Tr}\left(f(\bm{X} +\bm{Y}) - f(\bm{X})\right)$ on the space of derivatives $f$, and by the Riesz representation theorem, there exists a function $\eta_1(z)$ such that 
\begin{eqnarray}\label{eq3:thm:KSSF first order}
\mathrm{Tr}\left(f(\bm{X} +\bm{Y}) - f(\bm{X})\right) &=& \int f^{(1)}(z)\eta_1(z) dz.
\end{eqnarray}
This $\eta_1(z)$ is called the first-order spectral shift function, and it is determined uniquely by the function $f$, and operators $\bm{X}, \bm{Y}$.
$\hfill\Box$

Before considering higher order spectral shift function, we briefly review Taylor expansion of operator-valued functions and remainder term. 

Let \( \bm{X} \) and \( \bm{Y} \) be self-adjoint (or bounded) operators on a Hilbert space, and let \( f \) be a function that is sufficiently smooth (e.g., analytic or \( C^n \)) on a domain containing the spectrum of \( \bm{X} + t\bm{Y} \) for \( t \in [0,1] \). Then the operator-valued function \( f(\bm{X} + \bm{Y}) \) admits a Taylor-type approximation of order \( n-1 \), given by:
\begin{eqnarray}\label{eq1:Taylor exp}
f(\bm{X} + \bm{Y}) \approx \sum_{k=0}^{n-1} \frac{1}{k!} \left. \frac{d^k}{dt^k} f(\bm{X} + t\bm{Y}) \right|_{t=0}.
\end{eqnarray}
The remainder term \( \bm{R}_n(f, \bm{X}, \bm{Y}) \) is defined as the difference between the actual value and the Taylor polynomial:
\begin{eqnarray}\label{eq2:Taylor exp}
\bm{R}_n(f, \bm{X}, \bm{Y}) \define f(\bm{X} + \bm{Y}) -  \sum_{k=0}^{n-1} \frac{1}{k!} \left. \frac{d^k}{dt^k} f(\bm{X} + t\bm{Y}) \right|_{t=0}.
\end{eqnarray}

From Eq.~\eqref{eq2:Taylor exp}, the evaluation of $\bm{R}_n(f, \bm{X}, \bm{Y}) $ requires to determine each term $\frac{1}{k!} \left. \frac{d^k}{dt^k} f(\bm{X} + t\bm{Y}) \right|_{t=0}$. Then, we have to recall Theorem 7 from~\cite{chang2025GMOIFinite} (n-th derivative by GMOI) as below~\footnote{In~\cite{chang2025GMOIFinite}, the generalized multiple operator integrals (GMOIs) are developed in the matrix setting. However, Theorem~7 from~\cite{chang2025GMOIFinite} can be extended to the operator with continous spectra setting, resulting in the perturbation Theorem~\ref{thm:GMOI Perturbation Formula}, which preserves the same structural form as its matrix counterpart.
}.
\begin{theorem}\label{thm:n-th der by GMOI}
Given any operators $\bm{X}$, $\bm{Y}$ and $\tilde{\bm{X}} \define \bm{X}+t\bm{Y}$, we have
\begin{eqnarray}\label{eq1:thm:n-th der by GMOI}
\left.\frac{d^n f (\bm{X}+t \bm{Y})}{d t^n}\right\vert_{t=0}&=&n! T^{[\bm{X}]_1^{n+1}}_{f^{[n]}}([\bm{Y}]_1^n)-
(n-1)!\sum\limits_{j=1}^n T_{f^{[n]}}^{[\bm{X}]_1^{j-1},\bm{X}_N,\bm{X}_N,[\bm{X}]_{j+1}^{n}}([\bm{Y}]_1^n) \nonumber \\
&&-\sum\limits_{i=1}^{n-1}(n-i)!\left.\sum\limits_{j=1}^{n-i+1} \overline{\mathfrak{X}}^{(i)}([\bm{X}]_1^{j-1},\tilde{\bm{X}}, \bm{X},[\tilde{\bm{X}}]_{j+1}^{n-i+1})\right\vert_{t=0} \nonumber \\
&&- \sum\limits_{\ell_{n-1}=1}^{\left(\sum\limits_{k=1}^{\lceil\frac{n-1}{2}\rceil}\frac{(n-k)!}{k!(n-2k)!}\right)- 1(\mbox{~if $n$ is even})} a_{\ell_{n-1}}\left.\frac{d T_{f^{[n-1]}}^{\underline{[\tilde{\bm{X}},\tilde{\bm{X}}_N]}_{\ell_{n-1}}}([\bm{Y}]_1^{n-1})}{dt}\right\vert_{t=0} \nonumber \\
&&+ \left[ \sum\limits_{\rho=4}^n\sum\limits_{i_\rho=1}^{\gamma(\rho)} b_{i_\rho}\left.\overline{\mathfrak{X}}^{(n+2-\rho)}([\tilde{\bm{X}},\bm{X},\bm{X}_N]_{i_\rho})\right\vert_{t=0}\right]\bm{1}(\mbox{if $n \geq 4$}),
\end{eqnarray}
where $n \geq 2$, $a_{\ell_{n-1}}$ are coefficients for those GMOIs with parameters operators containing both $\bm{X}$ and $\bm{X}_N$ from $\frac{d^{n-1} f (\bm{X}+t \bm{Y})}{d t^{n-1}}$, $b_{i_\rho}$ are coefficients for those $\overline{\mathfrak{X}}^{(k)}([\tilde{\bm{X}},\bm{X},\bm{X}_N]_{i_\rho})$ from  $\frac{d^{n-1} f (\bm{X}+t \bm{Y})}{d t^{n-1}}$ with $k=1,2,\ldots,n-3$, and $[\tilde{\bm{X}},\bm{X},\bm{X}_N]_{i_\rho}$ is the array indexed by $i_\rho$ composed by $\tilde{\bm{X}}, \bm{X}$ and $\bm{X}_N$, respectively with length $\rho$. The term $\gamma(\rho)$ is determined by
\begin{eqnarray}\label{eq2:thm:n-th der by GMOI}
\gamma(\rho)&=&  \left[\frac{(\rho-2)!}{1!(\rho-3)!},  \frac{(\rho-3)!}{2!(\rho-5)!}, \ldots, 1 \right] \cdot [\rho-3, \rho-5, \ldots,0],
\end{eqnarray}
where $\cdot$ is the inner product operation. Note that $\bm{1}(\mbox{if $n \geq 4$})$ is the indicator condition function.
\end{theorem}

We have to present the following Lemma~\ref{lma:terms are linear functional} in order to determine the trace of the Taylor expansion remainder operator $\bm{R}_n(f, \bm{X}, \bm{Y})$ by integral representations. 
\begin{lemma}\label{lma:terms are linear functional} 
From Theorem~\ref{thm:n-th der by GMOI}, we have the following linear functional map categories:
\begin{eqnarray}\label{eq1:lma:terms are linear functional} 
\Upsilon_1:  f^{(n)} \rightarrow  \mathrm{Tr}\left(T^{[\bm{X}]_1^{n+1}}_{f^{[n]}}([\bm{Y}]_1^n)\right),
\end{eqnarray}
\begin{eqnarray}\label{eq2:lma:terms are linear functional} 
\Upsilon_2: f^{(n)} \rightarrow  \mathrm{Tr}\left(T_{f^{[n]}}^{[\bm{X}]_1^{j-1},\bm{X}_N,\bm{X}_N,[\bm{X}]_{j+1}^{n}}([\bm{Y}]_1^n)\right).
\end{eqnarray}

Moreover, we have the following bilinear functional map categories:
\begin{eqnarray}\label{eq3:lma:terms are linear functional} 
\Upsilon_3: f^{(n-i+1)}, f^{(n-i)}  \rightarrow  \mathrm{Tr}\left( \left.\overline{\mathfrak{X}}^{(i)}([\bm{X}]_1^{j-1},\tilde{\bm{X}}, \bm{X},[\tilde{\bm{X}}]_{j+1}^{n-i+1})\right\vert_{t=0} \right),
\end{eqnarray}
where $i=1,2,\ldots,n-1$.
\begin{eqnarray}\label{eq4:lma:terms are linear functional} 
\Upsilon_4: f^{(n)},  f^{(n-1)}  \rightarrow  \mathrm{Tr}\left(\left.\frac{d T_{f^{[n-1]}}^{\underline{[\tilde{\bm{X}},\tilde{\bm{X}}_N]}_{\ell_{n-1}}}([\bm{Y}]_1^{n-1})}{dt}\right\vert_{t=0} \right),
\end{eqnarray}
%
\begin{eqnarray}\label{eq3:lma:terms are linear functional} 
\Upsilon_5: f^{(\rho-1)},  f^{(\rho-2)} \rightarrow  \mathrm{Tr}\left(\left.\overline{\mathfrak{X}}^{(n+2-\rho)}([\tilde{\bm{X}},\bm{X},\bm{X}_N]_{i_\rho})\right\vert_{t=0}\right),
\end{eqnarray}
where $\rho=4,5,\ldots,n$.
\end{lemma}
\textbf{Proof:}
From the trace property and GMOI definitions, it is clear that $\Upsilon_1$ and $\Upsilon_2$ are linear functionals. From Theorem 3 in~\cite{chang2025GMOIFinite}, we also know that 
\begin{eqnarray}\label{eq4:lma:terms are linear functional} 
\lefteqn{T_{\beta^{[\zeta+1]}}^{[\bm{X}]_1^{j-1},\bm{C},\bm{D},[\bm{X}]_j^{\zeta}}([\bm{Y}]_1^{j-1},\bm{C}-\bm{D},[\bm{Y}]_j^{\zeta})=}\nonumber \\
&&T_{\beta^{[\zeta]}}^{[\bm{X}]_1^{j-1},\bm{C},[\bm{X}]_{j}^{\zeta}}([\bm{Y}]_1^{j-1},[\bm{Y}]_j^{\zeta}) - T_{\beta^{[\zeta]}}^{[\bm{X}]_1^{j-1},\bm{D},[\bm{X}]_{j}^{\zeta}}([\bm{Y}]_1^{j-1},[\bm{Y}]_j^{\zeta}) + \bar{\mathfrak{X}}([\bm{X}]_1^{j-1},\bm{C},\bm{D},[\bm{X}]_j^{\zeta}) \nonumber \\
&+& T_{\beta^{[\zeta+1]}}^{[\bm{X}]_1^{j-2},\bm{X}_{j-1,P},\bm{C}_N,\bm{D}_N,\bm{X}_{j,P},[\bm{X}]_{j+1}^{\zeta}}([\bm{Y}]_1^{j-1},\bm{C}-\bm{D},[\bm{Y}]_j^{\zeta})\nonumber \\
&+& T_{\beta^{[\zeta+1]}}^{[\bm{X}]_1^{j-2},\bm{X}_{j-1,P},\bm{C}_N,\bm{D}_N,\bm{X}_{j,N},[\bm{X}]_{j+1}^{\zeta}}([\bm{Y}]_1^{j-1},\bm{C}-\bm{D},[\bm{Y}]_j^{\zeta})\nonumber \\
&+& T_{\beta^{[\zeta+1]}}^{[\bm{X}]_1^{j-2},\bm{X}_{j-1,N},\bm{C}_N,\bm{D}_N,\bm{X}_{j,P},[\bm{X}]_{j+1}^{\zeta}}([\bm{Y}]_1^{j-1},\bm{C}-\bm{D},[\bm{Y}]_j^{\zeta})\nonumber \\
&+& T_{\beta^{[\zeta+1]}}^{[\bm{X}]_1^{j-2},\bm{X}_{j-1,N},\bm{C}_N,\bm{D}_N,\bm{X}_{j,N},[\bm{X}]_{j+1}^{\zeta}}([\bm{Y}]_1^{j-1},\bm{C}-\bm{D},[\bm{Y}]_j^{\zeta}).
\end{eqnarray}

From $\bar{\mathfrak{X}}([\bm{X}]_1^{j-1},\bm{C},\bm{D},[\bm{X}]_j^{\zeta})$ definition given by Theorem 3 in~\cite{chang2025GMOIFinite}, we know that $\mathrm{Tr}(\bar{\mathfrak{X}}([\bm{X}]_1^{j-1},\bm{C},\bm{D},[\bm{X}]_j^{\zeta}))$ is a bilinnear functional with respect to the function $\beta^{[\zeta+1]}$ and the function $\beta^{[\zeta]}$. Therefore,  $\Upsilon_3$,  $\Upsilon_4$, and $\Upsilon_5$ are bilinear functionals also.
$\hfill\Box$

We are ready to present the main Theorem~\ref{thm:remainder SS} in this section.
\begin{theorem}[Krein Spectral Shift for GDOI with continuous spectra for $(n-1)$-th order approximation]\label{thm:remainder SS}
By applying Theorem~\ref{thm:n-th der by GMOI} to Eq.~\eqref{eq2:Taylor exp} and arranging terms, we will have
\begin{eqnarray}\label{eq0:thm:remainder SS}
\bm{R}_n(f, \bm{X}, \bm{Y}) &=& \sum\limits_{i=1}^{n-1}\sum_j A_{i,j} T_{f^{[i]}}^{\mathfrak{P}_i(j)}([\bm{Z}_1,\ldots,\bm{Z}_i]_{j}) \nonumber \\
&& + 
\sum\limits_{i'=1}^{n-2}\sum_{j'} B_{i',j'} \overline{\mathfrak{X}}^{(\mathfrak{D}_{i'}( j'))}_{j'}([\bm{Z}'_1,\ldots, \bm{Z}'_{i'+2}]_{j'})\Big|_{t=0},
\end{eqnarray}
where 
\begin{itemize}
\item $\bm{Z}_i$ are operators of $\bm{Y}$;
\item $\bm{Z}'_{i'}$ are operators of $\bm{X}$, $\bm{X}_N$, $\tilde{\bm{X}}$, and $\tilde{\bm{X}}_N$;
\item $\mathfrak{P}_i (j)$ is the $j$-th array of parameter operators of a GMOI with length $i+1$; 
\item $\mathfrak{D}{i'}(j')$ denotes the derivative order associated with the $j'$-th array $[\bm{Z}'_1,\ldots, \bm{Z}'_{i'+2}]_{j'}$, which serves as the input to the operator-valued function $\overline{\mathfrak{X}}_{j'}$;
\item $A_{i,j}, B_{i',j'}$ are complex coefficients determined by $f^{[i]}, f^{[i'+1]},  f^{[i']}$with arguments as eigenvalues of operators from $\mathfrak{P}_i(j)$ and$[\bm{Z}'_1,\ldots, \bm{Z}'_{i'+2}]_{j'}$.
\end{itemize}

Then, there exist single-variable functions \( \eta_{i,j}(z) \) and bivariate functions \( \mathfrak{K}_{i', j'}(z_1, z_2) \) such that
\begin{eqnarray}\label{eq1:thm:remainder SS}
\mathrm{Tr}\left( \bm{R}_n(f, \bm{X}, \bm{Y}) \right) 
&=& \sum\limits_{i=1}^{n-1} \sum_j A_{i,j} \int\limits_z f^{(i)}(z) \, \eta_{i,j}(z) \, dz \nonumber \\
&& + \sum\limits_{i'=1}^{n-2} \sum_{j'} B_{i',j'} \int\limits_{z_1}\int\limits_{z_2} f^{(i'+1)}(z_1) \, \mathfrak{K}_{i', j'}(z_1, z_2) \, f^{(i')}(z_2) \, dz_1 dz_2,
\end{eqnarray}
where \( \mathfrak{K}_{i', j'}(z_1, z_2) \) also acts as the Hilbert–Schmidt kernel mapping \( f^{(i'+1)} \rightarrow f^{(i')} \).
\end{theorem}
\textbf{Proof:}
By taking the trace operation $\mathrm{Tr}$ at both sides of Eq.~\eqref{eq0:thm:remainder SS} and apply trace linear properties, we have
\begin{eqnarray}\label{eq2:thm:remainder SS}
\mathrm{Tr}\left(\bm{R}_n(f, \bm{X}, \bm{Y})\right)&=& \sum\limits_{i=1}^{n-1}\sum_j A_{i,j} \mathrm{Tr}\left(T_{f^{[i]}}^{\mathfrak{P}_i(j)}([\bm{Z}_1,\ldots,\bm{Z}_i]_{j})\right) \nonumber \\
&& + 
\sum\limits_{i'=1}^{n-2}\sum_{j'} B_{i',j'}\mathrm{Tr}\left(\overline{\mathfrak{X}}^{\mathfrak{D}_{i'}( j')}_{j'}([\bm{Z}'_1,\ldots, \bm{Z}'_{i'+2}]_{j'})\Big|_{t=0}\right).
\end{eqnarray}

From Lemma~\ref{lma:terms are linear functional}, as $f^{(i)} \rightarrow  \mathrm{Tr}\left(T_{f^{[i]}}^{\mathfrak{P}_i(j)}([\bm{Z}_1,\ldots,\bm{Z}_i]_{j})\right)$ is a linear functional, from Riesz representation theorem, we can find functions $\eta_{i,j}(z)$ to have
\begin{eqnarray}\label{eq2:thm:remainder SS}
 \mathrm{Tr}\left(T_{f^{[i]}}^{\mathfrak{P}_i(j)}([\bm{Z}_1,\ldots,\bm{Z}_i]_{j})\right)&=&\int\limits_z f^{(i)}(z)\eta_{i,j}(z) dz.
\end{eqnarray}
On the other hand, from Lemma~\ref{lma:terms are linear functional} again, as $f^{(i'+1)}, f^{(i')} \rightarrow \mathrm{Tr}\left(\overline{\mathfrak{X}}^{\mathfrak{D}_{i'}( j')}_{j'}([\bm{Z}'_1,\ldots, \bm{Z}'_{i'+2}]_{j'})\Big|_{t=0}\right)$ is a bilinear functional, we also can find functions $\mathfrak{K}_{i', j'}(z_1,z_2)$ by generalized Riesz representation theorem to have
\begin{eqnarray}\label{eq3:thm:remainder SS}
\mathrm{Tr}\left(\overline{\mathfrak{X}}^{\mathfrak{D}_{i'}( j')}_{j'}([\bm{Z}'_1,\ldots, \bm{Z}'_{i'+2}]_{j'})\Big|_{t=0}\right)&=& \int\limits_{z_1}\int\limits_{z_2} f^{(i'+1)}(z_1)\mathfrak{K}_{i', j'}(z_1,z_2) f^{(i')}(z_2)dz_1 dz_2.
\end{eqnarray}
This theorem is proved by combining Eq.~\eqref{eq2:thm:remainder SS} and Eq.~\eqref{eq3:thm:remainder SS}.
$\hfill\Box$

We will prepare Example~\ref{exp:SS R3} below to demonstrate Theorem~\ref{thm:remainder SS}.
\begin{example}\label{exp:SS R3}
Because we have 
\begin{eqnarray}\label{eq1:exp:SS R3}
\bm{R}_3(f, \bm{X}, \bm{Y}) &=& f(\bm{X} + \bm{Y}) -  \sum_{k=0}^{2} \frac{1}{k!} \left. \frac{d^k}{dt^k} f(\bm{X} + t\bm{Y}) \right|_{t=0} \nonumber \\
&=& T_{f^{[1]}}^{\bm{X} + \bm{Y},\bm{X}}(\bm{Y}) - \frac{1}{1!} \left. \frac{d}{dt} f(\bm{X} + t\bm{Y}) \right|_{t=0} - \frac{1}{2!} \left. \frac{d^2}{dt^2} f(\bm{X} + t\bm{Y}) \right|_{t=0} \nonumber \\
&=&  T_{f^{[1]}}^{\bm{X} + \bm{Y},\bm{X}}(\bm{Y})  - T^{[\bm{X}]_1^{2}}_{f^{[1]}}(\bm{Y}) -  \left( T^{[\bm{X}]_1^{3}}_{f^{[2]}}([\bm{Y}]_1^2) - \frac{1}{2}T^{\bm{X},\bm{X}_N,\bm{X}_N}_{f^{[2]}}([\bm{Y}]_1^2)  \right. \nonumber \\
&& \left. - \frac{1}{2}T^{\bm{X}_N,\bm{X}_N,\bm{X}}_{f^{[2]}}([\bm{Y}]_1^2) - \frac{1}{2}\left.\mathfrak{X}^{(1)}(\tilde{\bm{X}}, \bm{X}, \tilde{\bm{X}})\right\vert_{t=0} - \frac{1}{2}  \left.\mathfrak{X}^{(1)}(\bm{X}, \tilde{\bm{X}}, \bm{X})\right\vert_{t=0} \right) \nonumber \\
&=& T_{f^{[1]}}^{\bm{X} + \bm{Y},\bm{X}}(\bm{Y})  - T^{[\bm{X}]_1^{2}}_{f^{[1]}}(\bm{Y}) -  T^{[\bm{X}]_1^{3}}_{f^{[2]}}([\bm{Y}]_1^2)  + \frac{1}{2}  T^{\bm{X},\bm{X}_N,\bm{X}_N}_{f^{[2]}}([\bm{Y}]_1^2)  \nonumber \\
&& + \frac{1}{2} T^{\bm{X}_N,\bm{X}_N,\bm{X}}_{f^{[2]}}([\bm{Y}]_1^2) +\frac{1}{2} \left.\mathfrak{X}^{(1)}(\tilde{\bm{X}}, \bm{X}, \tilde{\bm{X}})\right\vert_{t=0} + \frac{1}{2} \left.\mathfrak{X}^{(1)}(\bm{X}, \tilde{\bm{X}}, \bm{X})\right\vert_{t=0},
\end{eqnarray}
from Theorem~\ref{thm:remainder SS}, we have
\begin{eqnarray}\label{eq2:exp:SS R3}
\mathrm{Tr}\left(\bm{R}_3(f, \bm{X}, \bm{Y})\right)&=& \int\limits_z f^{(1)}(z) \eta_{1,1}(z)dz - \int\limits_z f^{(1)}(z) \eta_{1,2}(z)dz -  \int\limits_z f^{(2)}(z) \eta_{2,1}(z)dz  \nonumber \\
&& + \frac{1}{2}\int\limits_z f^{(2)}(z) \eta_{2,2}(z)dz + \frac{1}{2} \int\limits_z f^{(2)}(z) \eta_{2,3}(z)dz  \nonumber \\
&& + \frac{1}{2} \int\limits_{z_1,z_2} f^{(2)}(z_1)\mathfrak{K}_{1, 1}(z_1,z_2) f^{(1)}(z_2)dz_1 dz_2   \nonumber \\
&& + \frac{1}{2} \int\limits_{z_1,z_2} f^{(2)}(z_1)\mathfrak{K}_{1, 2}(z_1,z_2) f^{(1)}(z_2)dz_1 dz_2.
\end{eqnarray}
Note that the terms on the right-hand side of Eq.\eqref{eq1:exp:SS R3} are in one-to-one correspondence with those on the right-hand side of Eq.\eqref{eq2:exp:SS R3}, with each term precisely aligned to its counterpart.
\end{example}

\bibliographystyle{IEEETran}
\bibliography{Inf_SpecialCase_and_MOI_Bib}

\end{document}